\documentclass[12pt]{amsart}
\newcommand\theoremnumbering{section}   
\newcommand\version{public}
\newcommand\finalized{yes}
\usepackage{fullpage,setspace}
\usepackage{tfrupee}  

\usepackage{mathrsfs} 
\newcommand\choosefont[1]{\usepackage{#1}}
\usepackage[bottom]{footmisc}  
\usepackage{enumerate}   
\usepackage{url,cite,fancyheadings}
\usepackage{enumerate}
\usepackage{asymptote}
\usepackage{amsmath,amssymb,amscd,mathrsfs,latexsym}
		\usepackage{rsfso}  
		\usepackage{mathtools}
\usepackage{mathdots}

\usepackage{tikz}

\usepackage{ifthen}     

\newcommand\pubpri[2]{%
\ifthenelse{\equal{\version}{public}}%
{{#1}}%
{\ifthenelse{\equal{\finalized}{no}}{\marginpar{\scshape\small Pubpri Alert}{#2}}{#2}}{}}
\newcommand\pubprinoalert[2]{%
\ifthenelse{\equal{\version}{public}}%
{{#1}}%
{#2}}

\newcommand\ignore[1]{}

\usepackage{color}
\providecommand\wantcolor{yes}   %
\definecolor{backgroundyellow}{cmyk}{.2,.1,.8,.2}
\definecolor{backgroundblue}{rgb}{0,0,1}
\definecolor{backgroundred}{rgb}{1,0,0}
\definecolor{backgroundmagenta}{cmyk}{0,1,0,0}
\ifthenelse{\equal{\wantcolor}{no}}{\renewcommand\color[1]{}}{}

\usepackage[shortlabels]{enumitem}

\newcommand\mysection{\section}

\newcommand\mysubsection{\subsection}

\newcommand\mysubsubsection[1]{%
		\subsubsection{\sffamily\upshape\mdseries #1}}

\newcommand\mysss{\mysubsubsection}

\usepackage{chngcntr}
\counterwithin{figure}{section}

\usepackage{url,hyperref}  
\hypersetup{colorlinks=true,
linkcolor=blue,
filecolor=magenta,
urlcolor=cyan,}
\usepackage{graphicx,wrapfig,epstopdf}

\providecommand{\theoremnumbering}{document}

\ifthenelse{\equal{\theoremnumbering}{section}}{\newtheorem{annotation}{Annotation}[section]}%
{\ifthenelse{\equal{\theoremnumbering}{subsection}}{\newtheorem{annotation}{Annotation}[subsection]}{\newtheorem{annotation}{Annotation}}}
\newtheorem{theorem}[annotation]{
		Theorem}
\newtheorem{lemma}[annotation]{
		Lemma}
\newtheorem{definition}[annotation]{
		Definition}
\newtheorem{corollary}[annotation]{
		Corollary}
\newtheorem{proposition}[annotation]{
		Proposition}

\newtheorem{example}[annotation]{
		Example}
\newcommand\bexample{\begin{example}\begin{rm}}
\newcommand\eexample{\end{rm}\hfill$\Box$\end{example}}
\newtheorem{examplenobox}[annotation]{
		Example}
\newcommand\bexamplenobox{\begin{examplenobox}\begin{rm}}
\newcommand\eexamplenobox{\end{rm}\end{examplenobox}}
\newtheorem{exercise}[annotation]{
		Exercise}
\newcommand\bexercise{\noindent\begin{exercise}\begin{rm}}
\newcommand\eexercise{\end{rm}\end{exercise}}
\newtheorem{notation}[annotation]{
		Notation}
\newcommand\bnotation{\begin{notation}\begin{rm}}
\newcommand\enotation{\end{rm}\end{notation}}

\newtheorem{remark}[annotation]{
		Remark}
\newcommand\bremark{\begin{remark}
\begin{upshape}}
\newcommand\eremark{\end{upshape}
\end{remark}}
\newenvironment{remark*}{%
\par\noindent{\scshape 
  Remark: }\begin{rm}}{\hfill\end{rm}\newline} 
\newcommand\bremarkstar{\begin{remark*}}
\newcommand\eremarkstar{\end{remark*}}
\newcommand\bdefn{\begin{definition}
\begin{upshape}}
\newcommand\edefn{\end{upshape}
\end{definition}}
\newtheorem{caveat}[annotation]{
		Caveat}
\newcommand\bcaveat{\begin{caveat}
\begin{upshape}}
\newcommand\ecaveat{\end{upshape}
\end{caveat}}
\newenvironment{caveatstar}{
\par\noindent{\scshape\bfseries
  Caveat: }\begin{rm}}{\end{rm}\newline} 
\newcommand\bcaveatstar{\begin{caveatstar}}
\newcommand\ecaveatstar{\end{caveatstar}}
\newenvironment{myproof}{%
\par\noindent{\scshape 
  Proof: }\begin{rm}}{\hfill$\Box$\end{rm}} 
\newcommand\bmyproof{\begin{myproof}}
\newcommand\emyproof{\end{myproof}}
\newenvironment{myproofnobox}{%
\par\noindent{\scshape Proof: }\begin{rm}}{\end{rm}\hfill}
\newcommand\bmyproofnobox{\begin{myproofnobox}}
\newcommand\emyproofnobox{\end{myproofnobox}}
\newenvironment{myproofof}[1]{%
\par\noindent{\scshape 
  Proof (of~#1): }\begin{rm}}{\hfill$\Box$\end{rm}} 
\newcommand\bmyproofof{\begin{myproofof}}
\newcommand\emyproofof{\end{myproofof}}
\newenvironment{myproofofnobox}[1]{%
\par\noindent{\scshape 
  Proof (of~#1): }\begin{rm}}{\end{rm}\hfill\newline} 
\newcommand\bmyproofofnobox{\begin{myproofofnobox}}
\newcommand\emyproofofnobox{\end{myproofofnobox}}
\newenvironment{solution}{%
\par\noindent{\scshape Solution: }\begin{rm}}{\hfill$\Box$\end{rm}\newline}
\newenvironment{solutionnobox}{%
\par\noindent{\scshape Solution: }\begin{rm}}{\end{rm}}
\newcommand\bsolution{\begin{solution}\begin{rm}}
\newcommand\esolution{\end{rm}\end{solution}}
\newcommand\bsolutionnobox{\begin{solutionnobox}\begin{rm}}
\newcommand\esolutionnobox{\end{rm}\end{solutionnobox}}
\newcommand\bthm{\begin{theorem}}
\newcommand\ethm{\end{theorem}}
\newcommand\bcor{\begin{corollary}}
\newcommand\ecor{\end{corollary}}
\newcommand\blemma{\begin{lemma}}
\newcommand\elemma{\end{lemma}}
\newcommand\bprop{\begin{proposition}}
\newcommand\eprop{\end{proposition}}
\newcommand\beqn{\begin{equation}}
\newcommand\eeqn{\end{equation}}
\newcommand\beqnstar{\begin{equation*}}
\newcommand\eeqnstar{\end{equation*}}
\usepackage{amsthm}
\usepackage{textcomp}
\usepackage{calrsfs,mathrsfs}  

\newcommand\mtitle[1]%
{\marginpar{\sffamily\raggedright\upshape\small #1}}

\providecommand\finalized{no}
\ifthenelse{\equal{\finalized}{yes}}{}{\marginparwidth=2cm}
\ifthenelse{\equal{\finalized}{yes}}%
{\newcommand\mylabel[1]{\label{#1}}}%
{\newcommand\mylabel[1]{\label{#1}\marginpar{[{\ttfamily\upshape\tiny #1}]}}}
\ifthenelse{\equal{\finalized}{yes}}%
{\newcommand\checked[1]{}}%
{\newcommand\checked[1]{\marginpar{[{\ttfamily\upshape\tiny CHECKED: #1}]}}}
\ifthenelse{\equal{\finalized}{yes}}%
{\newcommand\spellchecked[1]{}}%
{\newcommand\spellchecked[1]{\marginpar{[{\ttfamily\upshape\tiny SPELLCHECKED: #1}]}}}

\providecommand\version{public}   
\ifthenelse{\equal{\version}{public}}%
{\newcommand\mcomment[1]{}}%
{\newcommand\mcomment[1]{\marginpar{{\raggedright\sffamily\upshape\small
\begin{spacing}{0.75} #1\end{spacing}}}}}
\ifthenelse{\equal{\version}{public}}%
{\newcommand\fcomment[1]{}}%
{\newcommand\fcomment[1]{\footnote{#1}}}
\ifthenelse{\equal{\version}{public}}%
{\newcommand\comment[1]{}}%
{\newcommand\comment[1]{{\small #1}}}
\newcommand\tensor\otimes
\newcommand\wbar{\overline{w}}
\newcommand\wl{W_\lambda}
\newcommand\wm{W_\mu}

\newcommand\wlwwm{\wl\backslash W/\wm}

\newcommand\Domint{\Lambda^+}
\newcommand\barlsm{\overline{\lambda+\sigma\mu}}
\newcommand\wmu{w\mu}

\newcommand\mapphi{\Phi}
\newcommand\tmapphi{\tilde{\mapphi}}

\newcommand\lbr{\leq}
\newcommand\lneqbr{\lneq}
\newcommand\lsbr{<}

\newcommand\join{\vee}

\newcommand\meet{\wedge}

\newcommand\sequence{\mathfrak{s}}
\newcommand\intvl[1]{I(#1)}
\DeclareMathOperator\brmin{\bold{min}}
\DeclareMathOperator\brmax{\bold{max}}

\DeclareMathOperator\character{\bold{char}}
\newcommand\weyl{\mathfrak{w}}
\newcommand\weylj{\weyl_{j,j+1}}
\newcommand\wppprime{\weyl(\pi\concat\pi')}
\newcommand\plcpmw{\paths(\lambda,w,\mu)}
\newcommand\plwm{\paths(\lambda,w,\mu)}
\newcommand\plcpmphi{\paths(\lambda,\varphi,\mu)}
\newcommand\plcpmphip{\paths(\lambda,\varphi',\mu)}

\providecommand\st{\,|\,}
\newcommand\lieg{\mathfrak{g}}
\newcommand\lieh{\mathfrak{h}}
\newcommand\lieb{\mathfrak{b}}
\newcommand\bplus{\mathfrak{b}_+}
\newcommand\bminus{\mathfrak{b}_-}
\newcommand\vlambda{v_\lambda}
\newcommand\Vlambda{V_\lambda}
\newcommand\Vnu{V_\nu}
\newcommand\vmu{v_\mu}
\newcommand\vwmu{v_{w\mu}}
\newcommand\vwnotmu{v_{w_0\mu}}
\newcommand\vumu{v_{u\mu}}
\newcommand\Vmu{V_\mu}

\newcommand\vwm{v_{w\mu}}
\newcommand\lamdba\lambda
\newcommand\kostant{K}
\newcommand\klwm{\kostant(\lambda,w,\mu)}
\newcommand\klsm{\kostant(\lambda,\sigma,\mu)}
\newcommand\klum{\kostant(\lambda,u,\mu)}

\newcommand\koslwm{\klwm}
\newcommand\kmwinvl{\kostant(\mu,w^{-1},\lambda)}
\newcommand\koslum{\klum}

\newcommand\koslonem{\kostant(\lambda,1,\mu)}
\newcommand\koslwnotm{\kostant(\lambda,w_0,\mu)}

\newcommand\klmw\klwm
\newcommand\klwpm{\kostant(\lambda,w',\mu)}

\newcommand\paths{\mathcal{P}}
\newcommand\skewtab{\mathcal{T}}
\newcommand\ssytset{\mathcal{S}}
\newcommand\pathsl{\paths_\lambda}
\newcommand\pathsm{\paths_\mu}
\newcommand\skewtabm{\skewtab_\mu}
\newcommand\concat{\star}
\newcommand\deodmod{\textcircled{$\star$}}
\newcommand\circstar[2]{#1\,\deodmod\,#2}
\newcommand\pathsmlw{\pathsm^\lambda(w)}

\newcommand\pathsml{\pathsm^\lambda}
\newcommand\skewtabml{\skewtabm^\lambda}
\newcommand\ssytsetml{\ssytset_\mu^\lambda}
\newcommand\ssytsetmd{\ssytset_{\mu}[d]}
\newcommand\ssytsetmld{\ssytset_\mu^\lambda[d]}
\newcommand\skewtabmlw{\skewtabm^\lambda(w)}
\newcommand\skewtabmlwd{\skewtabm^\lambda(w)[d]}
\newcommand\skewtabmld{\skewtabm^\lambda[d]}
\newcommand\ssytsetmlw{\ssytsetml(w)}
\newcommand\ssytsetmlwd{\ssytsetml(w)[d]}
\newcommand\pathsmw{\pathsm(w)}
\newcommand\pil{\pi_\lambda}

\newcommand\pathsplp{\paths_{\pi_\lambda\concat\pi}}
\newcommand\pilp{\pil\concat\pi}
\newcommand\ha{h_\alpha}
\newcommand\hapi{\ha^\pi}

\newcommand\demazure{\Lambda}
\newcommand\demwnot{\demazure_{w_0}}
\newcommand\demw{\demazure_w}

\newcommand\expo[1]{e^{#1}}

\newcommand\dcoset{W_1\backslash W/W_2}
\newcommand\dcosetlm{W_\lambda\backslash W/W_\mu}

\newcommand\rv{r^V}

\newcommand\yseq{\underline{y}}

\newcommand\nubar{\overline{\nu}}

\newcommand\jsigmaw{J_{\sigma W_1}(w)}
\newcommand\jsigmav{J_{\sigma W_1}(v)}
\newcommand\jsigmawp{J_{\sigma W_1}(w')}
\newcommand\jsigmapw{J_{\sigma' W_1}(w)}
\newcommand\jsswminsw{J_{s\sigma W_1}(w\meet sw)}
\newcommand\jssvminsv{J_{s\sigma W_1}(v\meet sv)}
\newcommand\jswminsw{J_{\sigma W_1}(w\meet sw)}
\newcommand\jswmaxsw{J_{\sigma W_1}(w\join sw)}

\newcommand\jidvm{J_{W_1}(v_m)}

\newcommand\ksigmaw{K_{\sigma W_1}(w)}
\newcommand\ksigmav{K_{\sigma W_1}(v)}
\newcommand\ksigmawp{K_{\sigma W_1}(w')}

\newcommand\ksswmaxsw{K_{s\sigma W_1}(w\join sw)}
\newcommand\kssvmaxsv{K_{s\sigma W_1}(v\join sv)}
\newcommand\kswminsw{K_{\sigma W_1}(w\meet sw)}
\newcommand\kswmaxsw{K_{\sigma W_1}(w\join sw)}

\newcommand\length[1]{\textup{length}(#1)}

\newcommand\sld{\mathfrak{sl}_{d}}

\newcommand\taut{\tilde{\tau}}

\newcommand\phit{\tilde{\varphi}}
\newcommand\pmconcat{\paths_{1}\concat\cdots\concat\paths_{n}}

\newcommand\hecke{\mathcal{H}}
\newcommand\pidag{\pi^\dag}

\newcommand{\meets}{\meet_s}

\newcommand{\smeet}{\prescript{}{s}{\meet}}
\newcommand{\sjmeet}{\prescript{}{s_j}{\meet}}
\newcommand{\meetsj}{\meet_{s_j}}
\newcommand{\jmeet}{\prescript{}{j}{\meet}}
\newcommand{\meetj}{\meet_{j}}

\newcommand{\meetm}{\meet_{m}}

\newcommand{\meetone}{\meet_{1}}
\newcommand{\tmeet}{{}_t\meet}

\newcommand\pathstd{\paths_\textup{std}}
\newcommand\cryiso{\Gamma}
\newcommand\indir{\iota}
\newcommand\pathslv{\paths_{\lambda,v}}

\newcommand\pathsmv{\paths_{\mu,v}}
\newcommand\pathsmbv{\paths_{\mu,s_\beta v}}
\newcommand\pathslw{\paths_{\lambda,w}}
\newcommand\pathstdv{\paths_{\textup{std},v}}
\newcommand\pathstdw{\paths_{\textup{std},w}}

\newcommand\funf{\mathfrak{F}}
\newcommand\pathsmvf{\pathsmv(\funf)}
\newcommand\pathsmwf{\paths_{\mu,w}(\funf)}
\newcommand\pathsmvi{\pathsmv(\iota)}
\newcommand\pathsmwi{\paths_{\mu,w}(\iota)}
\newcommand\pathsmbvf{\pathsmbv(\funf)}

\newcommand\biject\Psi
\newcommand\cleq\prec
\newcommand\cgeq\succ
\newcommand\boxb{\mathbf{b}}
\newcommand\boxc{\mathbf{c}}
\newcommand\boxbp{\boxb'}
\newcommand\boxbr{\boxb^r}
\newcommand\boxbtint{\boxb^\tint}

\newcommand\boxbtintpone{{\boxb^{\tint+1}}}
\newcommand\boxbrpone{\boxb^{r+1}}
\newcommand\newboxn{\mathbf{n}}

\newcommand\mpart{{\mu_1'}}
\newcommand\ssyt{S}
\newcommand\ssytp{\ssyt'}
\newcommand\ssytr{\ssyt^{(r)}}

\newcommand\ssyttint{\ssyt_{\tint}}
\newcommand\ssyttintmone{\ssyt_{\tint-1}}

\newcommand\imsl{v}
\newcommand\imslsr{\imsl(r)}
\newcommand\imsltruncr{\imsl^{(r)}}

\newcommand\sigmat{\tilde{\sigma}}
\newcommand\procw{u}

\newcommand\procwrpone{\procw_{r+1}}
\newcommand\procwtint{\procw_{\tint}}
\newcommand\tint{p}
\newcommand\vecx{\underline{x}}
\newcommand\vecy{\underline{y}}
\newcommand\vecz{\underline{z}}
\newcommand\veca{\underline{a}}

\newcommand\vecg{\underline{g}}
\newcommand\vech{\underline{h}}

\newcommand\ltrunc[2]{{}_{(#2)}{#1}}
\newcommand\colmax{i_0}
\newcommand\nbem[1]{\underline{\mathbf{#1}}}
\newcommand\boxbt{\tilde{\boxb}}
\newcommand\boxbdelta{\boxb_\delta}
\newcommand\boxbtdelta{\boxbt_\delta}
\newcommand\boxbpdelta{\boxb'_\delta}
\newcommand\word{\mathbf{w}}
\newcommand\wordr{\mathbf{w}_{\textup{row}}}
\newcommand\wordc{\mathbf{w}_{\textup{col}}}
\newcommand\wordp{\mathbf{w}_{\textup{pos}}}
\newcommand\wordk{\mathbf{w}_k}
\newcommand\weight{\textup{wt}}
\def\a{\alpha}

\def\lieg{\mathfrak g}
\def\lieh{\mathfrak h}

\providecommand\tensor\otimes

\def\baselinestretch{1.15}
			\choosefont{mathpazo}
\usepackage{tkz-euclide}
\usepackage{tikz}
\usetikzlibrary{calc,angles,positioning,intersections,quotes,decorations.markings,shapes}
\usetikzlibrary{lindenmayersystems}
\begin{document}
\title[KK-modules via paths]{A study of Kostant-Kumar modules via Littelmann paths}
\author{Mrigendra Singh Kushwaha}
\address{The Institute of Mathematical Sciences, HBNI, Chennai, 600\,113, India}
\email{mrigendra@imsc.res.in, mrigendra154@gmail.com}
\author{K.~N.~Raghavan}
\address{The Institute of Mathematical Sciences, HBNI, Chennai, 600\,113, India}
\email{knr@imsc.res.in}
\author{Sankaran Viswanath}
\address{The Institute of Mathematical Sciences, HBNI, Chennai, 600\,113, India}
\email{svis@imsc.res.in}
\thanks{\noindent The authors acknowledge support under a DAE project grant to IMSc.}

\begin{abstract}

	We study, by means of Littelmann's theory of paths, Kostant-Kumar modules (KK modules for short),  which by definition are certain submodules of the tensor product of two irreducible integrable highest weight representations of a symmetrizable Kac-Moody algebra.   Our main result is an identification of a path model for any KK module as a subset of the well known path model for the tensor product consisting of concatenations of Lakshmibai-Seshadri paths. 
	The technical results about extremal elements in Coxeter groups 
	that we formulate and prove en route and the technique of their proofs should be of independent interest.  We also discuss the existence of PRV components and generalised PRV components in KK modules.


Specialising to the case of the special linear Lie algebra,  we record a decomposition rule for KK modules in terms of Littlewood-Richardson tableaux.
In this connection,  we present a new procedure to determine the permutation that is the initial element of the minimal standard lift of a semi-standard Young tableau.
The appendix, necessitated by the derivation of the tableau decomposition rule, 
	deals with standard concatenations of Lakshmibai-Seshadri paths of arbitrary shapes,  of which semi-standard Young tableaux form a very special case.
\end{abstract}  

\maketitle
\tableofcontents
\renewcommand\baselinestretch{1.16}
\mysection{Description of the results}\mylabel{s:intro}
\noindent
Let $\lieg$ be a symmetrizable Kac Moody Lie algebra.
For $\lambda$ a dominant integral weight,  let $\Vlambda$ denote the irreducible integrable representation of~$\lieg$ with highest weight $\lambda$, and $\vlambda$ a highest weight (non-zero) vector in~$\Vlambda$.

Fix dominant integral weights $\lambda$, $\mu$  and an element~$w$ of the Weyl group~$W$.  Let $\vwm$ denote a non-zero vector in the one-dimensional weight space of weight~$w\mu$ in~$\Vmu$.  The cyclic $\lieg$-submodule $U\lieg(v_\lambda\tensor v_{w\mu})$ of $\Vlambda\tensor\Vmu$, where $U\lieg$ denotes the universal enveloping algebra of~$\lieg$, is called a {\em Kostant-Kumar module\/}, {\em KK module\/} for short, and denoted $\klwm$.  
The following facts are elementary to prove and well known (see~\S\ref{s:kosdef}):
\begin{itemize}
	\item $\koslonem$ is the copy of $V_{\lambda+\mu}$ in the tensor product (where $1$ denotes the identity element of~$W$);  $\koslwnotm$ is the whole tensor product $\Vlambda\tensor\Vmu$, when $\lieg$ is of finite type and $w_0$ denotes the longest element of~$W$.
\item $\klwpm\subseteq\klwm$ for $w'$ in the Weyl group with $w'\leq w$ in the Bruhat order.
\end{itemize}

In this paper we study KK modules by means of Littelmann's theory of paths~\cite{litt:inv,litt:ann}.  
			Let $\pathsl$ and $\pathsm$ denote respectively the sets of Lakshmibai-Seshadri paths of shapes $\lambda$ and~$\mu$.
				Let $\pathsl\concat\pathsm$ be the set of paths~$\pi\concat\pi'$ with $\pi\in\pathsl$ and $\pi'\in\pathsm$, where $\pi\concat\pi'$ denotes the concatenation of $\pi$ and $\pi'$.  By work of Littelmann~\cite{litt:inv,litt:ann},   it is well known that $\pathsl\concat\pathsm$ is a ``path model'' for $V_\lambda\tensor V_\mu$.
		The technical novelty that we introduce to the study of KK modules via paths is the association~(\S\ref{ss:defkosset}) of an element denoted $\weyl(\pi\concat\pi')$ of the Weyl group~$W$ to each element $\pi\concat\pi'$ of $\pathsl\concat\pathsm$. To describe this association, let us assume for the sake of simplicity that both $\lambda$ and $\mu$ are regular (see \S\ref{ss:defkosset} for the general definition). Let $\tau$ and $\varphi$ be the elements of~$W$ representing respectively the final direction of~$\pi$ and the initial direction of~$\pi'$.    Consider the following subset of~$W$:
		\[ \{\tau'\varphi\st \tau'\in W, \textup{$\tau'\leq \tau^{-1}$ in the Bruhat order}\} \]
		This admits a unique minimal element (Corollary~\ref{c:takeaway}\,\eqref{i:takeaway2}), which $\weyl(\pi\concat\pi')$ is defined to be.  The key property of $\weyl(\pi\concat\pi')$ is that it is invariant under the action of root operators on~$\pi\concat\pi'$ (Proposition~\ref{p:wind}).   Towards the proof of this property, we establish in~\S\ref{s:extremal} some general results about extremal elements in Coxeter groups.   These results and especially the technique of their proofs are, we believe, of independent interest.

		Our main result (Theorem~\ref{t:pathkos}) identifies, under the possibly removable hypothesis that $\lieg$ is either of finite type or symmetric, a subset of~$\pathsl\concat\pathsm$ that is a path model for the KK module $\klwm$.  More precisely, we have,  for any element~$w$ in the Weyl group:
		\[ \character{\klwm}=\sum_{\pi\concat\pi'\in\plwm}\exp{(\pi\concat\pi')(1)} 
			\]
			where $\plwm$ denotes the set of those paths $\pi\concat\pi'$ in~$\pathsl\concat\pathsm$ with $\weyl(\pi\concat\pi')\leq w$.
The proof of this theorem combines the invariance of $\weyl(\pi\concat\pi')$ under root operators mentioned above with Joseph's decomposition rule~\cite[Theorems~5.25,~5.22]{joseph:demflag} for KK modules (which we recall below in~\S\ref{s:decompkos} and from which the restrictive hypothesis on~$\lieg$ is inherited) and two fundamental results of Littelmann, namely the path character formula~\cite[page~330]{litt:inv} and the isomorphism theorem~\cite[Theorem~7.1]{litt:ann}.

Let $w'$ be an element of the Weyl group and $\nu'$ denote the unique dominant conjugate of $\lambda+w'\mu$.    A copy of $V_{\nu'}$ in~$\Vlambda\tensor\Vmu$ (or any submodule thereof) is called a ``PRV component'',   following a conjecture of Parthasarathy, Ranga Rao, and Varadarajan~\cite{prv} that was proved independently by Kumar~\cite[Theorem~2.15]{kumar:prv} and Mathieu~\cite[Corollaire~3]{mathieu}.   The use of paths to prove the existence of PRV components is well known:  see e.g.\ \cite[\S7]{litt:inv} and \cite[\S2.7]{joseph:hc}.  The proof of Theorem~\ref{t:rkprvforkos} below, which is about PRV components in KK modules, follows this beaten track.   Theorem~\ref{t:gkprv} is about ``generalised PRV components'' in KK modules and follows Montagard~\cite[Theorem~3.1]{montagard}.

		The second part of the paper (\S\ref{s:tabkk}--\ref{s:proc-proof})  deals with the special case of the special linear Lie algebra.  In~\S\ref{s:tabkk} we deduce, from the general decomposition rule (Theorem~\ref{t:decompkos}),  a ``tableau decomposition rule'' (\S\ref{ss:tabkk}) for KK modules in terms of LR tableaux (LR is short for Littlewood-Richardson),  generalising the classical LR rule (see e.g.~\cite{fulton:yt,mcld}) for decomposing $\Vlambda\tensor\Vmu$.   As is well known and in any case easy to see
		(\S\ref{sss:st-to-ssyt}, \S\ref{sss:isoTS}),  the LR tableaux that figure in the classical LR rule can be identified with semi-standard Young tableaux (SSYT for short) of shape~$\mu$ that are $\lambda$-dominant.  
		  Again as is well known (and recalled with proof in \S\ref{ss:etastd}),  any SSYT has an interpretation as a standard concatenation of LS paths (LS is short for Lakshmibai-Seshadri) and therefore gives rise to a standard tuple and corresponding minimal standard lift (in the sense of~\S\ref{ss:standard}). 

		  Fix a LR tableau~$T$ that by the classical rule makes a contribution to the decomposition of the whole tensor product~$\Vlambda\tensor\Vmu$.  Let $\imsl$ be the permutation that is the initial element of the minimal standard lift of the SSYT corresponding to~$T$.  The tableau decomposition rule---of which we state three variations (\eqref{e:kktab}, \eqref{e:kktabalt}, \eqref{e:kktabpoly}) for ease of ready reference---says that $T$ contributes to the decomposition of the KK module $\klwm$ if and only if $\imsl\leq w$ in the Bruhat order. 

		  The standard concatenation of LS paths attached to $T$ may in turn be interpreted, via the isomorphism theorem of Littelmann mentioned above (in connection with the proof of the main Theorem~\ref{t:pathkos}),  as an LS path~$\pi$ of shape~$\mu$.     The fact that the permutation~$\imsl$ represents the initial direction of~$\pi$ is what is involved in the derivation (in \S\ref{ss:pfkktab}) of the tableau decomposition rule.  
		  While this fact will come as no surprise to experts in the theory of paths---it is hinted at by Littelmann already in~\cite{litt:inv} and later stated in~\cite[\S11]{litt_plactic} with a sketch of proof---we could not find a reference to a complete proof.   The appendix below is an attempt to address this inadequacy in the literature. It develops with complete proofs the necessary results about standard concatenations of LS paths of arbitrary shapes.

		  Given the important role played by the initial element~$\imsl$ of the minimal standard lift of a SSYT in the tableau decomposition rule for KK modules,   we present in~\S\ref{sss:ssyt-to-w} a new procedure to produce it from the SSYT.     The naivest way to produce $\imsl$ would perhaps be to closely follow the proof of its existence in~\S\ref{ss:standard}, or in other words repeatedly apply the construction in Deodhar's lemma (Proposition~\ref{p:deodhar}).     Lascoux-Sch\"utzenberger~\cite{ls:keys} give an algorithm to produce~$v$ via their notion of the ``right key'' associated to a SSYT.  Willis~\cite{willis} gives a more efficient algorithm to produce the right key.  Our procedure stands apart from these other ones.  The justification for it (namely the proof that it produces~$\imsl$) is provided in~\S\ref{s:proc-proof}.

It is a pleasure to thank B.~Ravinder and R.~Venkatesh for pointers to relevant literature and useful discussions.   We cannot thank the anonoymous referee enough for insightful comments on the first version of this paper,   which helped fill a rather large gap in our awareness of the literature,  and,  in particular,  enabled weakening of the hypothesis in our results.

\mysection{Generalities on extremal elements in Coxeter groups}\mylabel{s:extremal}
\noindent
Proposition~\ref{p:wind} is the main technical result of the present paper.   It relies (via Lemma~\ref{l:wind}) on the existence and properties of unique minimal elements in a particular kind of subsets of a Weyl group.    The purpose of this section is to formulate and prove the required result about these minimal elements---Corollary~\ref{c:takeaway} below---in the more natural context of Coxeter groups.   The arguments leading up to the result are all elementary.  
The reader willing to accept it at face value without proof may want to skip this section at a first pass.

The results in~\S\ref{ss:brcoset},~\ref{ss:deodhar}, and~\ref{ss:standard} are well known (e.g., from~\cite{deo:ca:87,lak-litt-mag} as specifically indicated in a few places below),  but we have included them because we need them and it is easier to prove them ab initio in our set up than to refer to  sources.   It not only makes the paper more self-contained but also more readable with these results stated and proved rather than just quoted.

\mysubsection{Notation for this section}\mylabel{ss:ntnextrml}
Let $(W,S)$ be a Coxeter system.  
Let $\lbr$ denote the (strong) Bruhat order on~$W$.
For a subset $K$ of the Coxeter group~$W$,   we let $\brmin{K}$ and $\brmax{K}$ denote respectively the unique minimal and unique maximal elements of~$K$ in the Bruhat order (if they do exist).    For $u$ in~$W$ and $s$ in~$S$, the elements $u$ and $su$ (respectively $u$ and $us$) are comparable.    Thus $\brmin{\{u,su\}}$, $\brmax{\{u,su\}}$, $\brmin{\{u,us\}}$, and $\brmax{\{u,us\}}$ make sense.   We denote these respectively by $u\meet su$, $u\join su$, $u\meet us$, and $u\join us$.


\mysubsection{The results}\mylabel{ss:rsltsextrml}
The following basic fact is repeatedly applied in this section:
\begin{equation} \label{e:parallel}
	\textup{For elements $u\lbr v$ in $W$ and $s$ in~$S$, we have $u\meet su\lbr v\meet sv$ and $u\join su\lbr v\join sv$.}\tag{*}
\end{equation}
The ``right analogue'' of the above fact asserts:  $u\meet us\lbr v\meet vs$ and $u\join us\lbr v\join vs$ (under the same hypothesis).   Only the left analogues of the ``one sided'' results below are explicitly stated.  Their right analogues hold good too.


\bremark    Suppose that $u\leq v$ is a covering relation in~$W$ (that is, $\length{u}=\length{v}-1$ and $u=t v$ for some reflection~$t$ in~$W$).   Then, if for $s$ in~$S$,  we have $sv<v$ and $u<su$, then $t=s$.   Indeed, it follows from~\eqref{e:parallel} that $u\leq sv$, but then equality is forced since $u$ and $sv$ have the same length.
	\eremark

A simple application of \eqref{e:parallel} gives:
\begin{proposition}\mylabel{p:ppar}
	Suppose that a subset $K$ of the Coxeter group~$W$ has a unique minimal element $u$ under~$\lbr$.   Then, for any $s$ in~$S$,  the subset $K\cup sK$ also has a unique minimal element under $\lbr$, namely, $u\meet su$.  Analogously, if $K$ admits a unique maximal element~$v$ under $\lbr$, then $v\join sv$ is the unique maximal element of $K\cup sK$.
\end{proposition}
\bcor\mylabel{c:cppar}  Let $\sequence$: $s_1$, $s_2$, \ldots\  be a (possibly infinite) sequence of simple reflections (elements of $S$).    For $\sequence'$: $s_{i_1}$, $s_{i_2}$, \ldots, $s_{i_m}$ a finite subsequence of $\sequence$,  let $w(\sequence')$ denote the element $s_{i_m}s_{i_{m-1}}\cdots s_{i_1}$ of the Coxeter group (note the order reversal).  Let $K$ be a subset of~$W$ with a unique minimal element~$u$ with respect to~$\lbr$.  Then $\cup_{\sequence'}w(\sequence')K$, where the union runs over all finite subsequences $\sequence'$ of~$\sequence$, has a unique minimal element $u_\infty$, the stable value of $u_j$ as $j\to\infty$,  where $u_j$ is recursively defined:  $u_0=u$, and $u_{j+1}=u_j\meet s_{j+1}u_i$ for $0\leq j$.
\ecor
\bmyproof For $j$ a non-negative integer, let $\sequence_j$ denote the subsequence $s_1$, $s_2$, \ldots, $s_j$ of $\sequence$. By a repeated application of~Proposition~\ref{p:ppar},  we see that $u_j$ is the unique minimal element of $K_j:=\cup_{\sequence'}w(\sequence')K$,  where the union runs over subsequences of~$\sequence_j$.   Since the subsets $K_j$ increase with~$j$,  it follows that $u_{j+1}\lbr u_j$.   Since any decreasing sequence in the Bruhat order stabilizes,  we conclude that $u_j$ is constant for $j$ sufficiently large.
\emyproof
\bremark\mylabel{r:cppar}  What about the maximal analogue of Corollary~\ref{c:cppar}?  Let $K$ be a subset of~$W$ that has a unique maximal element~$v$.    With notation as in the proof just above,  we conclude analogously that $v_k$ is the unique maximal element of~$K_k$,  where $v_k$ is defined recursively as follows:  $v_0=v$, and $v_{i+1}=v\join s_{i+1}v$ for $0\leq i$. Since the $K_k$ increase with~$k$,  we have $v_{k+1}\geq v_k$.    If the $v_k$ stabilize to a stable value $v_\infty$ as $k\to\infty$ (which in general need not happen),  then $v_\infty$ is the unique maximal element of $\cup_{k\geq0}K_k$.  In particular,  the maximal analogue holds if the sequence~$\sequence$ is finite.
\eremark
We now apply Corollary~\ref{c:cppar} (and its right analogue) in two special cases. First, let $\sigma$ be an element of~$W$ and, with notation as in the corollary, choose the sequence $\sequence$: $s_1$, $s_2$, \ldots, $s_m$ of elements of~$S$ to be such that $s_ms_{m-1}\cdots s_1$ is a reduced expression for $\sigma$.   Then $\{w(\sequence')\st \textup{$\sequence'$ is a subsequence of~$\sequence$}\}$ equals $\intvl{\sigma}:=\{\sigma'\in W\st \sigma'\lbr\sigma\}$.  We conclude that $\intvl{\sigma}K$ has a unique minimal element and further that this element is the unique minimal element in $\intvl{\sigma}u$.   Now applying the right analogue of this argument,   we obtain:
\bcor\mylabel{c:cpparintvl} Let $K$ be a subset of the group~$W$ that admits a unique minimal element $u$.   Then,  for any two elements $\sigma_1$ and $\sigma_2$ of~$W$,  the set $\intvl{\sigma_1}K\intvl{\sigma_2}$ has a unique minimal element, and this element is the unique minimal element in $\intvl{\sigma_1}u\intvl{\sigma_2}$.
\ecor
\noindent
The special case of the above result (as also its maximal analogue, namely, Corollary~\ref{c:cpparmax}) when $K$ is a singleton and $\sigma_2$ is the identity element appears in~\cite[Lemma~11\,(i)]{lak-litt-mag}.  
\bcor\mylabel{c:cpparintvl1} Let $\sigma$ and $\varphi$ be elements of~$W$, and $s$ an element of~$S$.   Suppose that $\sigma s\lsbr\sigma$.    Then $\brmin{\intvl{\sigma}\varphi}$ equals either $\brmin{\intvl{\sigma}s\varphi}$ or $\brmin{\intvl{\sigma s}\varphi}$ accordingly as $s\varphi\lsbr\varphi$ or $\varphi\lsbr s\varphi$.
\ecor
\bmyproof Choose a sequence~$\sequence$: $s=s_1$, $s_2$, \ldots, $s_m$ such that $s_ms_{m-1}\cdots s_1$ is a reduced expression for $\sigma$.  Let $\varphi_i$, $\varphi'_i$, and $\varphi''_i$ be sequences defined recursively as follows:
\begin{itemize}
	\item $\varphi_0=\varphi$,  and $\varphi_{i+1}=\varphi_i\meet s_{i+1}\varphi_i$ for $0\leq i<m$
	\item $\varphi'_0=s\varphi$,  and $\varphi'_{i+1}=\varphi'_i\meet s_{i+1}\varphi'_i$ for $0\leq i<m$
	\item $\varphi''_1=\varphi$,  and $\varphi''_{i+1}=\varphi''_i\meet s_{i+1}\varphi''_i$ for $1\leq i<m$
		\end{itemize}
By Corollary~\ref{c:cppar},  $\brmin{\intvl{\sigma}\varphi}$, $\brmin{\intvl{\sigma}s\varphi}$, and $\brmin{\intvl{\sigma s}\varphi}$ are equal respectively to $\varphi_m$, $\varphi'_m$, and $\varphi''_m$.    

First suppose that $s\varphi \lsbr \varphi$.   Then $\varphi_1=\varphi'_1$: indeed, $\varphi_1=\varphi\meet s\varphi=s\varphi$, and $\varphi'_1=s\varphi\meet s(s\varphi))=s\varphi$.   Thus $\varphi_i=\varphi'_i$ for all $1\leq i$, and in particular for $i=m$.

Now suppose that $\varphi \lsbr s\varphi$.   Then $\varphi_1=\varphi''_1$: indeed, $\varphi_1=\varphi\meet s\varphi=\varphi$, and $\varphi''_1=\varphi$ by definition.   Thus $\varphi_i=\varphi''_i$ for all $1\leq i$, and in particular for $i=m$.      
\emyproof

Towards a second application of Corollary~\ref{c:cppar}, let $S_1$ be a subset of~$S$ and $W_1$ the subgroup of~$W$ generated by~$S_1$.  With notation as in Corollary~\ref{c:cppar}, choose the sequence $\sequence$: $s_1$, $s_2$, \ldots\ to consist of elements of~$S_1$ and such that every element of~$W_1$ arises as $w(\sequence')$ for some finite subsequence $\sequence'$ of~$\sequence$.  Then $\{w(\sequence')\st \textup{$\sequence'$ is a subsequence of~$\sequence$}\}$ equals~$W_1$.  We conclude that $W_1K$ has a unique minimal element.   Now applying the right analogue of this argument,  we obtain:
\bcor\mylabel{c:cpparpar} Let $K$ be a subset of the group~$W$ that admits a unique minimal element $u$.   Then,  for any two standard parabolic subgroups $W_1$ and $W_2$ of~$W$,  the set $W_1KW_2$ has a unique minimal element, and this element is the unique minimal element in $W_1uW_2$.  In particular,  any double coset of a pair of standard parabolic subgroups has a unique minimal element.
\ecor
\noindent
Of course,  when the subgroups~$W_1$ and $W_2$ are finite, Corollary~\ref{c:cpparpar} is a special case of~Corollary~\ref{c:cpparintvl}.   Indeed,  letting $w_1$ and $w_2$ be the unique maximal elements of~$W_1$ and $W_2$ respectively,  we have $W_1=\intvl{w_1}$ and $W_2=\intvl{w_2}$.

As the maximal analogues of the above two corollaries,  we have:
\bcor\mylabel{c:cpparmax}   Let $K$ be a subset of~$W$ having a unique maximal element~$v$.   Then,  for any two elements~$\sigma_1$ and $\sigma_2$ of~$W$,   the set $\intvl{\sigma_1}K\intvl{\sigma_2}$ has a unique maximal element, namely, the unique such element in $\intvl{\sigma_1}v\intvl{\sigma_2}$.  In particular,  for any two finite standard parabolic subgroups~$W_1$ and $W_2$ of~$W$, the union $W_1KW_2$ of double cosets has a unique maximal element,  namely, the unique such element in $W_1vW_2$.
\ecor
\bremark\mylabel{r:deodhar1} {\scshape (Relation to Deodhar's $\star$ operation.)}
In~\cite[Lemma~2.4]{deo:ca:87}, Deodhar states:    there exists a unique associative binary operation $\star$ on~$W$ such that $w\star\textup{id}=w$ and $w\star s=w\join ws$ for all $w\in W$ and $s\in S$.    The uniqueness is clear.    For the proof of the existence, we define $w\star x:=\brmax{\intvl{w}\intvl{x}}$ for all $w$ and $x$ in~$W$ ($\brmax{\intvl{w}\intvl{x}}$ exists by Corollary~\ref{c:cpparmax}).   It is easy to verify, using Corollary~\ref{c:cpparmax}, that this operation has the requisite properties:     $\brmax{\intvl{w}\intvl{\textup{id}}}=\brmax{\intvl{w}}=w$;  $\brmax{\intvl{w}\intvl{s}}=\brmax{w\intvl{s}}=w\join ws$;  and
\begin{gather*}\brmax{\intvl{w}\intvl{\brmax{\intvl{x}\intvl{y}}}}=
\brmax{\intvl{w}{\brmax{(\intvl{x}\intvl{y}})}}=
\brmax{\intvl{w}\intvl{x}\intvl{y}}\\=
\brmax{(\brmax{\intvl{w}\intvl{x}})\intvl{y}}=
\brmax{\intvl{\brmax{(\intvl{w}\intvl{x}})}\intvl{y}}, \end{gather*}
 so associativity holds.  

 We have:
 \begin{itemize}
	 \item   The unique maximal element of $\intvl{\sigma_1}K\intvl{\sigma_2}$ in Corollary~\ref{c:cpparmax} is $\sigma_1\star v\star\sigma_2$.
	 \item 
			 $\intvl{w}\intvl{x}=\intvl{w\star x}$ for all $w$ and $x$ in~$W$.
	 \item  Let $K$ be a subset of~$W$ with a unique maximal element $v$.  For any collection $\sigma_1$, \ldots, $\sigma_s$, $\tau_1$, \ldots, $\tau_t$ of elements in~$W$,  the set $\intvl{\sigma_1}\cdots\intvl{\sigma_s}K\intvl{\tau_1}\cdots\intvl{\tau_t}$ equals $\intvl{\sigma_1\star\cdots\star\sigma_s}K\intvl{\tau_1\star\cdots\star\tau_t}$ and admits a unique maximal element, namely, $\sigma_1\star\cdots\star\sigma_s\star v \star\tau_1\star\cdots\star\tau_t$.  \end{itemize}

\eremark
\bremark\mylabel{r:hecke}   Consider the specialized Hecke algebra~$\hecke$ defined as the associative algebra with identity (over say a field $k$) generated by variables $T_s$, $s\in S$, and subject to the relations $T_s^2=T_s$ (for all $s$ in~$S$) and the braid relations.  For $w\in W$,  let $T_w$ be the element $T_{s_{i_1}}T_{s_{i_2}}\cdots T_{s_{i_r}}$ of~$\hecke$ where $s_{i_1}s_{i_2}\cdots s_{i_r}$ is a reduced expression for $w$:   this definition does not depend on the choice of reduced expression because the braid relations are satisfied.  The algebra $\hecke$ is just the semigroup algebra of the semigroup $W$ with respect to the $\star$ operation as in Remark~\ref{r:deodhar1}:   $\{T_w\st w\in W\}$ is a basis for~$\hecke$ and $T_wT_x=T_{w\star x}$ for all $w$ and $x$ in~$W$.

Let $kW$ be the free $k$-vector space with elements of~$W$ as a basis.  We can make $kW$ to be $\hecke$-$\hecke$ bimodule as follows.   For $s$ in $S$, let $\meets$ denote the (left) operator on~$W$ defined by $\meets w:=sw\meet w$ and $\smeet$ the right operator on~$W$ defined by $w\smeet:=w\meet ws$ (for $w$ in~$W$).  The linear extensions of the operators $\smeet$ and $\meets$ to $kW$ are denoted by the same symbols.    We have,  for $s$, $t$ in~$S$ and $w$ in~$W$:
\begin{itemize}
	\item $\meets(\meets w)=
	\meets w$  and $(w\smeet)\smeet=w\smeet$.
	\item $(\meets w)\tmeet=\meets(w\tmeet)$.
	\item Let $s_1$, \ldots, $s_m$ be a sequence of elements of $S$ such that $s_m s_{m-1} \cdots s_1$ is a reduced expression for an element $\sigma$ of~$W$.  Then 
		$\meetm\cdots\meetone(w)=\brmin{\intvl{\sigma}w}$ (see the paragraph preceding Corollary~\ref{c:cpparintvl}) and analogously $w\meetone\cdots\meetm=\brmin{w\intvl{\sigma^{-1}}}$, where $\meetj$ and $\jmeet$ stand for $\meetsj$ and $\sjmeet$ respectively.   Thus the operators $\smeet$ (respectively $\meets$), $s\in S$, satisfy the braid relations.
\end{itemize}
Thus, letting $T_s$, $s\in S$, act on~$kW$ on the left by $\meets$ and on the right by $\smeet$,  we get a bimodule structure on~$kW$.

\eremark
\mysubsection{Bruhat order on double coset spaces}\mylabel{ss:brcoset}
Let $W_1$ and $W_2$ be standard parabolic subgroups of~$W$.     It is convenient to identify the coset space $\dcoset$ as a subset of~$W$ via the association $W_1uW_2\mapsto\brmin{W_1uW_2}$.   The Bruhat order on $\dcoset$ is the restriction to this subset of the Bruhat order on~$W$.
\bcor\mylabel{c:brcoset}  Given two elements $W_1uW_2$, $W_1vW_2$ in~$\dcoset$,   we have $W_1uW_2\leq W_1vW_2$ in Bruhat order if and only if there exist $u'$ in~$W_1uW_2$ and $v'$ in~$W_1vW_2$ such that $u'\leq v'$.   In particular,  if $\brmin{W_1uW_2}\leq v'$ for some $v'$ in $W_1vW_2$,  then $\brmin{W_1uW_2}\leq v''$ for every $v''\in W_1vW_2$.
\ecor
\bmyproof
For the only if part,   just take $u'=\brmin{W_1uW_2}$ and $v'=\brmin{W_1vW_2}$.   For the if part,  apply Corollary~\ref{c:cpparpar} with $K=\{u',v'\}$.     Since $u'\leq v'$,   we conclude that $\brmin{W_1u'W_2}=\brmin{W_1KW_2}$.    But $\brmin{W_1u'W_2}=\brmin{W_1uW_2}$ and $\brmin{W_1KW_2}\leq\brmin{W_1v'W_2}=\brmin{W_1vW_2}$ since $v'\in K$. \emyproof 
\bcor\mylabel{c2:brcoset}   For $u$ an element of $W$ and $s$ an element of $S$,  suppose that $suW_1\lneqbr uW_1$.   Then for every $v$ in $uW_1$,  we have $sv<v$.   Conversely,  if $su<u$ for the unique minimal element $u$ in $uW_1$,  then $suW_1\lneqbr uW_1$,  and $su$ is the unique minimal element in $suW_1$.
\ecor
\bmyproof   For the first statement, observe the following: if $v<sv$,  then, by Corollary~\ref{c:brcoset}, $uW_1=vW_1\lbr svW_1=suW_1$,  a contradiction.      For the first part of the converse,  observe that $suW_1\lbr uW_1$ by Corollary~\ref{c:brcoset},  and that equality cannot hold (if $su$ were to belong to $uW_1$ the minimality of $u$ in $uW_1$ would be contradicted).   For the second part of the converse,   suppose that $x\in suW_1$.  Then $sx\in uW_1$, and so $u\lbr sx$, which means $su=u\meet su\lbr sx\meet x\lbr x$.
\emyproof
\bremark\mylabel{r:c2:brcoset} For $u$ in $W$ and $s$ in~$S$,   it is possible that $usW_1\lneq uW_1$ but there exists $v$ in~$uW_1$ with $vs>v$.  For example,   let 
\[
	W\ =\ \langle s_1, s_2, s_3\st s_1^2=s_2^2=s_3^2=1,\ s_1s_2s_1=s_2s_1s_2,\ s_2s_3s_2=s_3s_2s_3,\ s_1s_3=s_3s_1 \rangle
	\]
$W_1=\langle s_1,s_3\rangle$, $u=s_1s_2$, and $s=s_2$.   Then $uW_1=\{s_1s_2,s_1s_2s_1,s_1s_2s_3,s_1s_2s_1s_3\}$, the minimal element in $uW_1$ is $u$,  and $usW_1=W_1\lneq uW_1$,   but $vs>v$ for $v=s_1s_2s_3$ in $uW_1$.
\eremark
\bprop\mylabel{p:brcoset}  Suppose that $W_1uW_2\leq W_1vW_2$.   Then,  given $u'$ in $W_1uW_2$,  there exists $v'$ in $W_1vW_2$ with $u'\leq v'$.
\eprop
\bmyproof Proceed by induction on the length of~$u'$.    Let $u_0$ be the minimal element in $W_1uW_2$.  We have $u_0\leq u'$.   If $u'=u_0$,   then $u'\leq v'$ for any $v'$ in $W_1vW_2$ (since $u_0\leq v_0$  by definition,  where $v_0$ is the minimal element of~$W_1vW_2$).   Now suppose that $u_0\lneq u'$.  Then there exists either $s\in S\cap W_1$ such that $su'<u'$,  or $t\in S\cap W_2$ such that $u't<u'$.   Let us suppose that the former condition holds (the case when the latter holds is handled analogously).   Observe that $su'$ belongs to $W_1uW_2$.  By induction,  there exists $v'$ in~$W_1vW_2$ such that $su'\leq v'$.    By~\eqref{e:parallel} at the beginning of this section,  we have \( u'\leq su'\join u'\leq v'\join sv'\).  But $v'\join sv'$ belongs to $W_1vW_2$.  
\emyproof
\mysubsection{Deodhar's Lemma}\mylabel{ss:deodhar} \noindent
Let $W_1$ be a standard parabolic subgroup of~$W$ and let $\sigma$, $w$ be elements of~$W$.  
Set $\jsigmaw:=\{v\in\sigma W_1\st w\lbr v\}$.   
\bprop\mylabel{p:deodhar:0} Let $s$ be in $S$.
\begin{enumerate}
	\item\label{i:deodp01} $\jsigmaw\supseteq\jsigmawp$ for $w\leq w'$.
	\item\label{i:deodp02} $\jsigmaw\subseteq s\jsswminsw$
	\item\label{i:deodp03} $\jsigmaw$ is non-empty if and only if $wW_1\leq\sigma W_1$.
\end{enumerate}
\eprop
\bmyproof  Statement~\eqref{i:deodp01} is immediate.    For~\eqref{i:deodp02},  just observe that $w\leq x$ implies $w\meet sw\leq x\meet sx\leq  sx$.   As for~\eqref{i:deodp03}, the only if part is trivial,  and the if part follows from Proposition~\ref{p:brcoset}.~\emyproof

\begin{lemma}\mylabel{l:deodhar}  Suppose that $\sigma$ is the minimal element of~$\sigma W_1$,  and let $s\in S$ be such that $s\sigma<\sigma$. Then
	\begin{enumerate}
		\item\label{i:ldeod1} $sx<x$ for any $x$ in $\jsigmaw$; and $v<sv$ for any $v$ in $\jsswminsw$.
		\item\label{i:ldeod2} $\jsigmaw=s \jsswminsw$.   In particular,  $\jsigmaw=\jswminsw=\jswmaxsw$.
	\item\label{i:ldeod3} If either $\jsigmaw$ or $\jsswminsw$ has a unique minimal element $u$,  then so does the other and $su$ is that unique minimal element.\end{enumerate} \end{lemma}
		\bmyproof
		\eqref{i:ldeod1}:    We have $s\sigma W_1\lneq \sigma W_1$ by the second part of Corollary~\ref{c2:brcoset}.     From the first part of that corollary,   it follows that $sx<x$ for any $x$ in $\sigma W_1$.    Since $\jsigmaw\subseteq\sigma W_1$ by definition,   the first statement follows.   The second statement too follows from the first assertion in Corollary~\ref{c2:brcoset}.

		\eqref{i:ldeod2}:   For the first assertion,  given item~\eqref{i:deodp02} of Proposition~\ref{p:deodhar:0},  it is enough to show that $\jsigmaw\supseteq s\jsswminsw$.  Suppose $x$ belongs to $\jsswminsw$.  Then evidently $sx$ belongs to $\sigma W_1$.   By item~\eqref{i:ldeod1} (of the present lemma),  we have $sx=x\join sx$.    Since $x\geq w\meet sw$ by hypothesis,  we have $sx=x\join sx\geq (w\meet sw)\join s(w\meet sw)\geq w$.

	To see that $\jsigmaw=\jswminsw$ (respectively, $\jsigmaw=\jswmaxsw$), put $v=w\meet sw$ (respectively, $v=w\join sw$).   Then $v\meet sv=w\meet sw$, and so, by the first assertion, $\jsigmav=s\jssvminsv=s\jsswminsw=\jsigmaw$.   


	\eqref{i:ldeod3}~Suppose $u$ is the unique minimal element in~$\jsswminsw$.  Then $su$ belongs to $\jsigmaw$ by~\eqref{i:ldeod2}. Let $x$ be any element in $\jsigmaw$.  We have $x\join sx=x$ by~\eqref{i:ldeod1}, and $sx\in\jsswminsw$ by~\eqref{i:ldeod2}, so $u\lbr sx$.  Thus $su\lbr u\join su\lbr sx\join s(sx)=sx\join x =x$.

	Suppose $u$ is the unique minimal element in~$\jsigmaw$.  Then $su$ belongs to $\jsswminsw$ by~\eqref{i:ldeod2}. Let $y$ be any element in $\jsswminsw$. We have $u\meet su=su$ by~\eqref{i:ldeod1}, and $sy\in\jsigmaw$ by~\eqref{i:ldeod2}, so $u\lbr sy$.  Thus $su = u\meet su\lbr sy\meet s(sy)=sy\meet y \lbr y$.
	\emyproof
	\begin{proposition}\mylabel{p:deodhar}	\textup{ ({\scshape``Deodhar's lemma''}, see e.g.~\cite[Lemma~5.8]{lms3})}
Suppose that $\jsigmaw$ is not empty.  Then it contains a unique minimal element.   Moreover, this element can be constructed recursively as follows:  let $\sigma$ be the minimal element in its coset $\sigma W_1$ and  let $s_1$, \ldots, $s_m$ be a sequence of elements of~$S$ such that $s_1\cdots s_m$ is a reduced expression for~$\sigma$;   put $v_0=w$ and $v_j=v_{j-1}\meet s_{j}v_{j-1}$ for $1\leq j\leq m$;   then $v_m$ belongs $W_1$,  and the minimal element of $\jsigmaw$ is just $\sigma v_m$. \end{proposition} 
	\bmyproof  By repeated application of Lemma~\ref{l:deodhar}~\eqref{i:ldeod2},   we have $\jsigmaw=\sigma \jidvm$.   Since $\jsigmaw$ is non-empty,  it follows that $\jidvm$ is non-empty as well, which means $v_m$ belongs to $W_1$ and is the unique minimal element in $\jidvm$.    That $\sigma v_m$ is the unique minimal element of~$\jsigmaw$ follows by a repeated application of Lemma~\ref{l:deodhar}~\eqref{i:three}.  \emyproof
\ignore{
	\bcor\mylabel{c:deodhar}  Suppose that $\jsigmaw$ is not empty.   Let $\sigma'$ in $W$ be such that $\sigma W_1\leq\sigma' W_1$.   Then $\jsigmapw$ is not empty and $\brmin{\jsigmaw}\leq\brmin{\jsigmapw}$.
	\ecor
	\bmyproof Let $\sigma$, $\sigma'$ be respectively the minimal elements in~ $\sigma W_1$, $\sigma'W_1$.   Proceed by induction on the length of~$\sigma'$.    If this length is $0$,  then $\sigma=\sigma'$ ($=\textup{identity}$) and so $\jsigmaw=\jsigmapw$ and in particular the statement is true.   Now suppose that the length of~$\sigma'$ is at least~$1$.   Let $s$ be in $S$ such that $s\sigma'<\sigma'$.  We consider two cases.

	First suppose that $s\sigma <\sigma$.    Then 

	First suppose that $s\sigma>\sigma$.   Then $\sigma=\sigma\meet s\sigma\leq \sigma'\meet s\sigma'=s\sigma'$.   So $\sigma W_1\leq s\sigma'W_1$.   By the induction hypothesis $\jsigmaw\leq $

	\emyproof
}
	\bremark\mylabel{r:p:deodhar}  We make a few remarks regarding the construction in Proposition~\ref{p:deodhar}.  
	\begin{enumerate}
		\item The element $v_m$ in the statement of Proposition~\ref{p:deodhar} is $\brmin{\intvl{\sigma^{-1}}w}$ (see the third item in the list in Remark~\ref{r:hecke}).  Thus $\jsigmaw$ is non-empty if and only if $\brmin{\intvl{\sigma^{-1}}w}$ belongs to $W_1$ and in this case its unique minimal element is $\sigma(\brmin{\intvl{\sigma^{-1}}w})$.
		\item  
		 Given a double coset $W_1\sigma W_2$,  where $W_2$ is also a standard parabolic subgroup,  there need not be a unique minimal element  among those in the double coset that are $\geq w$.  Consider for instance the following simple example.  Let $W$ be the Weyl group of type $A_2$:  $S=\{s_1,s_2\}$ and $W=\langle s_1, s_2\st s_1^2=s_2^2=1, s_1s_2s_1=s_2s_1s_2\rangle$.   Put $W_1=W_2=\langle s_1\rangle=\{1,s_1\}$,  $\sigma=s_2$, and $w=s_1$.   Observe that among the elements in $W_1s_2W_2=\{s_2,s_1s_2,s_2s_1,s_1s_2s_1\}$ that are $\geq s_1$,  there are two minimal ones, namely, $s_1s_2$ and $s_2s_1$.
	 \item\label{i:r:p:deodhar:ww'} Let $w'$ be an element of $w$ such that $w\leq w'$.  Then, evidently, $\brmin{\jsigmaw}\leq\brmin{\jsigmawp}$ (assuming that both sets are non-empty).
		\item 
			Let $\sigma'$ be an element of~$W$ such that $\sigma W_1\leq\sigma'W_1$.   Suppose that $\jsigmaw$ is non-empty.   Then $\jsigmapw$ is non-empty too:  for any $u$ in $\sigma W_1$ (and in particular for any $u$ in $\jsigmaw$),  there exists, by Proposition~\ref{p:brcoset}, $u'$ in~$\sigma'W_1$ such that $u\leq u'$, and evidently $u'$ belongs to $\jsigmawp$ when $u$ belongs to $\jsigmaw$.   However,  it need not be true that $\brmin{\jsigmaw}\leq\brmin{\jsigmapw}$,  as the following simple example shows.   Let $W$ be the Weyl group of type $A_3$: \[\langle s_1, s_2, s_3\st s_1^2=s_2^2=s_3^2=(s_1s_2)^3=(s_2s_3)^3=(s_1s_3)^2=1\rangle\]  Let $W_1$ be the parabolic subgroup $\langle s_2\rangle=\{1,s_2\}$,  $w=s_2$,  $\sigma=s_1s_3$, and $\sigma'=s_1s_2s_3$.   Then $\jsigmaw$ is non-empty, and $\brmin{\jsigmaw}=s_1s_3s_2\not\leq\brmin{\jsigmapw}=s_1s_2s_3$.
		\item\label{i:deofns}
			{(\scshape Deodhar's functions $f$ and $g$~\cite{deo:ca:87})}   Suppose that $\sigma$ is the least element in the coset $\sigma W_1$.  Let $\sigma'$ be an element of the Weyl group that is the least in its coset $\sigma'W_1$ and suppose that $\sigma W_1\leq\sigma'W_1$ (equivalently $\sigma\leq\sigma'$).

			Given $x$ in $W_1$,  there exists $x'$ in~$W_1$ such that $\sigma x\leq\sigma' x'$ (see Proposition~\ref{p:brcoset}).     Deodhar in~\cite[Lemma~2.2]{deo:ca:87} states that there is a function $f':W_1\to W_1$ (depending upon $\sigma$ and $\sigma'$) such that, for $x$ and $x'$ in~$W_1$, $\sigma x\leq \sigma' x'$ if and only if $f'x\leq x'$. (Deodhar writes $g$ for $f'$.)  In the notation of Remark~\ref{r:hecke},  this function is given by $f'x=\meet_{\sigma'^{-1}}\sigma x$.  Indeed it follows from Proposition~\ref{p:deodhar} that $\meet_{\sigma'^{-1}}\sigma x$ has the required property.

			Deodhar also asserts the existence of a function $f:W_1\to W_1$ (depending upon $\sigma$ and $\sigma'$) such that, for $x$ and $x'$ in~$W_1$, $\sigma x\leq \sigma' x'$ if and only if $x\leq fx'$.  To describe this function,  let $s_1$, \ldots, $s_m$ be a sequence of elements of~$S$ such that $s_1\cdots s_m$ is a reduced expression for $\sigma'$.   Put $v_0=\sigma$, and $v_{i+1}:=v_i\meet s_iv_i$ inductively for $1\leq i\leq m$.    Let $p$ be least, $0\leq p\leq m$, such that $v_p=\textup{id}$:  such a $p$ exists because $\sigma\leq\sigma'$.     Let $s_{i_1}$, \ldots, $s_{i_q}$ be the subsequence of $s_{p+1}$, \ldots, $s_m$ consisting precisely of those elements that belong to $W_1$.   Putting $z:=s_{i_1}\star\cdots\star s_{i_q}$,  Deodhar's function $f$ is given by $fx'=z\star w'$. 
			We omit the justification (which is not difficult) since we have no use in what follows for this function~$f$.
		\item\label{i:ksigmaw} (see~\cite[Lemma~11\,(ii)]{lak-litt-mag}) Put $\ksigmaw:=\{v\in \sigma W_1\st v\leq w\}$ ($=I(w)\cap \sigma W_1$).   If $\ksigmaw$ is non-empty (or, equivalently, $\sigma\leq w$) then it has a unique maximal element.   Indeed, writing $w$ as $\sigma' x'$ where $\sigma'$ is the minimal element in $wW_1$ and $x'$ is in~$W_1$,   this unique maximal element is $\sigma fx'$, where $f$ is Deodhar's function of the previous item.
	\end{enumerate}
	\eremark
\ignore{
	\mysubsection{Maximal elements in cosets of parabolics intersected with lower ideals}\mylabel{ss:deodhar2}  As in the last subsection,
let $W_1$ be a standard parabolic subgroup of~$W$ and let $\sigma$, $w$ be elements of~$W$.  
Set $\ksigmaw:=\{u\in\sigma W_1\st u\lbr w\}$.   
\bprop\mylabel{p:deodhar2:0} Let $s$ be in $S$.
\begin{enumerate}
	\item\label{i:deodp201} $\ksigmaw\subseteq\ksigmawp$ for $w\leq w'$.
	\item\label{i:deodp202} $\ksigmaw\subseteq s\ksswmaxsw$
\end{enumerate}
\eprop
\bmyproof  Statement~\eqref{i:deodp01} is immediate.    For~\eqref{i:deodp02},  just observe that $x\leq w$ implies $sx \leq x\join sx\leq w\join sw$.
\emyproof

\begin{lemma}\mylabel{l:deodhar2}  Suppose that $\sigma$ is the minimal element of~$\sigma W_1$,  and let $s\in S$ be such that $s\sigma<\sigma$. Then
	\begin{enumerate}
		\item\label{i:ldeod21} $sx<x$ for any $x$ in $\ksigmaw$; and $v<sv$ for any $v$ in $\ksswmaxsw$.
		\item\label{i:ldeod22} $\ksigmaw=s \ksswmaxsw$.   In particular,  $\ksigmaw=\kswminsw=\kswmaxsw$.
	\item\label{i:ldeod23}  If either $\ksigmaw$ or $\ksswmaxsw$ has a unique maximal element $v$,  then so does the other and $sv$ is that unique maximal element.\end{enumerate} \end{lemma}
		\bmyproof
		\eqref{i:ldeod21}:    We have $s\sigma W_1\lneq \sigma W_1$ by the second part of Corollary~\ref{c2:brcoset}.     From the first part of that corollary,   it follows that $sx<x$ for any $x$ in $\sigma W_1$.    Since $\ksigmaw\subseteq\sigma W_1$ by definition,   the first statement follows.   The second statement too follows from the first assertion in Corollary~\ref{c2:brcoset}.

		\eqref{i:ldeod22}:   For the first assertion,  given item~\eqref{i:deodp202} of Proposition~\ref{p:deodhar2:0},  it is enough to show that $\ksigmaw\supseteq s\ksswmaxsw$.  Suppose $x$ belongs to $\ksswmaxsw$.  Then evidently $sx$ belongs to $\sigma W_1$.   By item~\eqref{i:ldeod21} (of the present lemma),  we have $sx=x\meet sx$.    Since $x\leq w\join sw$ by hypothesis,  we have $sx=x\meet sx\leq (w\join sw)\meet s(w\join sw)\leq w$.

	To see that $\ksigmaw=\kswminsw$ (respectively, $\ksigmaw=\kswmaxsw$), put $v=w\meet sw$ (respectively, $v=w\join sw$).   Then $v\join sv=w\join sw$, and so, by the first assertion, $\ksigmav=s\kssvmaxsv=s\ksswmaxsw=\ksigmaw$.   


	\eqref{i:ldeod3}~Suppose $u$ is the unique minimal element in~$\jsswminsw$.  For $x$ in $\jsigmaw$, it follows (since $sx\in\jsswminsw$ and so $u\lbr sx$, and $x\join sx=x$), that $su\lbr u\join su\lbr sx\join x=x$.
	Suppose $u$ is the unique minimal element in~$\jsigmaw$.  For $v$ in $\jsswminsw$, it follows (since $sv\in\jsigmaw$ and so $u\lbr sv$, and $su=u\meet su$), that $su= u\meet su\lbr sv\join v\lbr v$.
	\emyproof
	\begin{proposition}\mylabel{p:deodhar}	\textup{\scshape (``Deodhar's lemma'')}
Suppose that $\jsigmaw$ is not empty.  Then it contains a unique minimal element.   Moreover, this element can be constructed recursively as follows:   let $s_1$, \ldots, $s_m$ be a sequence of elements of~$S$ such that $s_1\cdots s_m$ is a reduced expression for the unique minimal element~$\sigma$ in $\sigma W_1$;   put $v_0=w$ and $v_j=v_{j-1}\meet s_{j}v_{j-1}$ for $1\leq j\leq m$;   then $v_m$ belongs $W_1$,  and the minimal element of $\jsigmaw$ is just $\sigma v_m$. \end{proposition} 
	\bmyproof  By repeated application of Lemma~\ref{l:deodhar}~\eqref{i:ldeod2},   we have $\jsigmaw=\sigma \jidvm$.   Since $\jsigmaw$ is non-empty,  it follows that $\jidvm$ is non-empty as well, which means $v_m$ belongs to $W_1$ and is the unique minimal element in $\jidvm$.    That $\sigma v_m$ is the unique minimal element of~$\jsigmaw$ follows by a repeated application of Lemma~\ref{l:deodhar}~\eqref{i:three}.  \emyproof
\ignore{
	\bcor\mylabel{c:deodhar}  Suppose that $\jsigmaw$ is not empty.   Let $\sigma'$ in $W$ be such that $\sigma W_1\leq\sigma' W_1$.   Then $\jsigmapw$ is not empty and $\brmin{\jsigmaw}\leq\brmin{\jsigmapw}$.
	\ecor
	\bmyproof Let $\sigma$, $\sigma'$ be respectively the minimal elements in~ $\sigma W_1$, $\sigma'W_1$.   Proceed by induction on the length of~$\sigma'$.    If this length is $0$,  then $\sigma=\sigma'$ ($=\textup{identity}$) and so $\jsigmaw=\jsigmapw$ and in particular the statement is true.   Now suppose that the length of~$\sigma'$ is at least~$1$.   Let $s$ be in $S$ such that $s\sigma'<\sigma'$.  We consider two cases.

	First suppose that $s\sigma <\sigma$.    Then 

	First suppose that $s\sigma>\sigma$.   Then $\sigma=\sigma\meet s\sigma\leq \sigma'\meet s\sigma'=s\sigma'$.   So $\sigma W_1\leq s\sigma'W_1$.   By the induction hypothesis $\jsigmaw\leq $

	\emyproof
}
	\bremark\mylabel{r:p:deodhar}  We make a few remarks regarding the construction in Proposition~\ref{p:deodhar}.  
	\begin{enumerate}
		\item The element $v_m$ in the statement of Proposition~\ref{p:deodhar} is $\brmin{\intvl{\sigma^{-1}}w}$ (see the third item in the list in Remark~\ref{r:hecke}).  Thus $\jsigmaw$ is non-empty if and only if $\brmin{\intvl{\sigma^{-1}}w}$ belongs to $W_1$ and in this case its unique minimal element is $\sigma(\brmin{\intvl{\sigma^{-1}}w})$.
		\item  
		 Given a double coset $W_1\sigma W_2$,  where $W_2$ is also a standard parabolic subgroup,  there need not be a unique minimal element  among those in the double coset that are $\geq w$.  Consider for instance the following simple example.  Let $W$ be the Weyl group of type $A_2$:  $S=\{s_1,s_2\}$ and $W=\langle s_1, s_2\st s_1^2=s_2^2=1, s_1s_2s_1=s_2s_1s_2\rangle$.   Put $W_1=W_2=\langle s_1\rangle=\{1,s_1\}$,  $\sigma=s_2$, and $w=s_1$.   Observe that among the elements in $W_1s_2W_2=\{s_2,s_1s_2,s_2s_1,s_1s_2s_1\}$ that are $\geq s_1$,  there are two minimal ones, namely, $s_1s_2$ and $s_2s_1$.
	 \item Let $w'$ be an element of $w$ such that $w'\geq w$.  Then, evidently, $\brmin{\jsigmawp}\geq\brmin{\jsigmaw}$ (assuming that both sets are non-empty).
		\item 
			Let $\sigma'$ be an element of~$W$ such that $\sigma'W_1\geq\sigma W_1$.   Suppose that $\jsigmaw$ is non-empty.   Then $\jsigmawp$ is non-empty too:  for any $u$ in $\sigma W_1$ (and in particular for any $u$ in $\jsigmaw$),  there exists, by Proposition~\ref{p:brcoset}, $u'$ in~$\sigma'W_1$ such that $u'\geq u$, and evidently $u'$ belongs to $\jsigmawp$ when $u$ belongs to $\jsigmaw$.   However,  it need not be true that $\brmin{\jsigmawp}\geq\brmin{\jsigmaw}$,  as the following simple example shows.   Let $W$ be the Weyl group of type $A_3$: \[\langle s_1, s_2, s_3\st s_1^2=s_2^2=s_3^2=(s_1s_2)^3=(s_2s_3)^3=(s_1s_3)^2=1\rangle\]  Let $W_1$ be the parabolic subgroup $\langle s_2\rangle=\{1,s_2\}$,  $w=s_2$,  $\sigma=s_1s_3$, and $\sigma'=s_1s_2s_3$.   Then $\jsigmaw$ is non-empty, and $\brmin{\jsigmaw}=s_1s_3s_2\not\leq\brmin{\jsigmawp}=s_1s_2s_3$.
	\end{enumerate}
	\eremark
	\mysubsection{Maximal elements in cosets of parabolics intersected with lower ideals}\mylabel{ss:deodhar2}
} 
\bexample\mylabel{x:deodhar} 
Let $W$ be the group of permutations of~$[n]$ with $S$ as the set of simple transpositions $(1,2)$, \ldots, $(n-1,n)$.  For $r$, $1\leq r< n$,  let $W_r$ be the standard maximal parabolic subgroup generated by all simple transpositions except $(r,r+1)$.    Let $\sigma$ be a permutation of~$[n]$ with one-line notation $\sigma_1\sigma_2\ldots \sigma_n$.    For $\sigma$ to be of minimal length in its coset $\sigma W_r$,   it is necessary and sufficient that the sequences $\sigma_1\ldots\sigma_r$ and $\sigma_{r+1}\ldots\sigma_n$ are both increasing.   Suppose that this is the case.

Let $w$ be another permutation.  For $J_{\sigma W_r}(w)$ to be non-empty,  it is necessary and sufficient that $w^r_1\leq \sigma_1$, \ldots, $w^r_r\leq \sigma_r$, where $w^r_1<\ldots<w^r_r$ are the elements $w_1$, \ldots, $w_r$ arranged in increasing order.   Suppose that this is the case.

Let us suppose further that the $w_1$, \ldots, $w_r$ were themselves in increasing order,  so that $w^r_1=w_1$, \ldots, $w^r_r=w_r$.  In this case,
$\tau:=\brmin J_{\sigma W_{r}}(w)$ is determined as follows.    Put $\tau_j=\sigma_j$ for $1\leq j\leq r$.  For $j>r$,  an induction on~$j$ determines $\tau_j$ as follows.     Let $\tau^{j-1}_1<\ldots<\tau^{j-1}_{j-1}$ be $\tau_1$, \ldots, $\tau_{j-1}$ in increasing order and $w^j_1<\ldots<w^j_j$ be $w_1$, \ldots, $w_j$ in increasing order.    Then $\tau_j=w^j_k$,   where $k$ is the largest, $1\leq k\leq j$,  such that $\tau^{j-1}_{k-1}<w^j_k$ (we put $\tau^{j-1}_0=-\infty$).


As an example, let $n=6$, $r=3$, $w=145362$, and $\sigma=246135$.  Then, by the recipe above, $\tau=246153$.

Justification for the recipe appears later: see~Example~\ref{x:deodharcont} in~\S\ref{ss:deodharcont}.~\hfill\eexample
	\mysubsection{Standard tuples and standard lifts}\mylabel{ss:standard}   Let $W_1$, \ldots, $W_m$ be a sequence of standard parabolic subgroups of~$W$ and let $\theta=(\tau_1,\ldots,\tau_m)\in W/W_1\times \cdots\times W/W_m$.     We call $\theta$ {\em standard\/} if there exists a chain $\taut_1\geq\ldots\geq\taut_m$ of elements in $W$ such that $\taut_jW_j=\tau_j$ for $1\leq j\leq m$.   Such a chain is called a {\em standard lift\/} of~$\theta$.

  	Fix a standard tuple $\theta$ of cosets and a standard lift of it as above.  Put $\sigma_m=\brmin{\tau_m}$.   Observe that $\taut_{m-1}\geq\taut_m\geq\sigma_m$.   This means that $J_{\tau_{m-1}}(\sigma_m)$ is not empty, and so it has a unique minimal element by Proposition~\ref{p:deodhar}.   Put $\sigma_{m-1}=\brmin{J_{\tau_{m-1}}(\sigma_m)}$.   Proceeding this way, choose inductively $\sigma_j=\brmin{J_{\tau_j}(\sigma_{j+1})}$ for $j$ equal to $m-2$, $m-3$, \ldots, $1$.  We call the chain $\sigma_1\geq\ldots\geq\sigma_m$ the {\em minimal standard lift\/} of~$\theta$.
	We denote by $\weyl(\theta)$ the initial element $\sigma_1$ of the minimal standard lift of~$\theta$.  

	Let $\sigma_1\geq\ldots\geq\sigma_n$ be the minimal standard lift of~$\theta$ (for some standard tuple~$\theta$ of cosets).   As is easily observed (by a downward induction on~$j$),   $\sigma_j\leq\taut_j$ for $1\leq j\leq m$ for any standard lift $\taut_1\geq\ldots\geq\taut_m$ of~$\theta$.    Furthermore,  this property characterises the minimal standard lift.
\mysubsection{The takeaway from this section}\mylabel{ss:takeaway}
Finally,  we isolate the takeaway from this section in its (admittedly strange and whimsical) specific form that will be invoked later. 
\bcor\mylabel{c:takeaway}  Let $(W,S)$ be a Coxeter system,  $\tau$, $\varphi$ be elements of~$W$,  and $W_1$, $W_2$ be standard parabolic subgroups of~$W$.    Let $\intvl{\tau^{-1}}$ be the Bruhat interval $\{w\in W\st w\lbr\tau^{-1}\}$.  
Then:
\begin{enumerate}
	\item\label{i:takeaway1} If $\tau'$, $\varphi'$ are elements of~$W$ such that $\tau W_1=\tau' W_1$ and $\varphi W_2=\varphi' W_2$,   
		then:\[{W_1\intvl{{\tau'}^{-1}}\varphi'W_2}={W_1\intvl{\tau^{-1}}\varphi W_2}\]
	\item\label{i:takeaway2}
There exists a unique minimal element in $W_1\intvl{\tau^{-1}}\varphi W_2$, denoted $\brmin{W_1\intvl{\tau^{-1}}\varphi W_2}$. 
\item\label{i:takeaway3}
Let $s$ be an element of $S$ such that $s\tau\lsbr\tau$  and $s\varphi W_2\geq \varphi W_2$.  Then
		\[
		\brmin{W_1\intvl{\tau^{-1}}\varphi W_2}=
		\brmin{W_1\intvl{\tau^{-1}s}\varphi W_2} \]
\item \label{i:takeaway4}
Let $s$ be an element of $S$ such that $s\varphi<\varphi$ and $s\tau W_1\leq\tau W_1$.  Then
		\[
		\brmin{W_1\intvl{\tau^{-1}}\varphi W_2}=
		\brmin{W_1\intvl{\tau^{-1}}s\varphi W_2} \]
\end{enumerate}
\ecor
\bmyproof 
For~\eqref{i:takeaway1},  it being evident that $W_1\intvl{\tau^{-1}}\varphi W_2=W_1\intvl{\tau^{-1}}\varphi'W_2$,   it is enough to show that $W_1\intvl{\tau^{-1}}=W_1\intvl{\tau'^{-1}}$.
By Corollary~\ref{c:cpparpar}, $W_1\tau^{-1}=W_1\tau'^{-1}$ has a unique minimal element, say $\sigma^{-1}$.    It is enough to show that $W_1\intvl{\sigma^{-1}}=W_1\intvl{\tau^{-1}}$.  Since $\sigma^{-1}\lbr \tau^{-1}$, it follows that $\intvl{\sigma^{-1}}\subseteq\intvl{\tau^{-1}}$ and so $W_1\intvl{\sigma^{-1}}\subseteq W_1\intvl{\tau^{-1}}$.  To prove the other way containment,  write $\tau^{-1}=u^{-1}\sigma^{-1}$ with $u\in W_1$ and $\ell(\tau)=\ell(u)+\ell(\sigma)$,  where $\ell$ stands for ``length''.  Suppose that $\rho^{-1}\lbr \tau^{-1}$.  Then $\rho^{-1}=v^{-1}\rho'^{-1}$ with $v\lbr u$ (hence $v\in W_1$) and $\rho'\lbr\sigma$.  Thus $W_1\rho^{-1}=W_1v^{-1}\rho'^{-1}=W_1\rho'^{-1}$,  and we are done.

Assertion \eqref{i:takeaway2} follows from Corollaries~\ref{c:cpparintvl} and~\ref{c:cpparpar}.

Proof of~\eqref{i:takeaway3}: By Corollary~\ref{c:cpparintvl}, $\brmin{\intvl{\tau^{-1}}\varphi}$ exists. First suppose that $s\varphi>\varphi$.    Then, by~Corollary~\ref{c:cpparintvl1}, $\brmin{\intvl{\tau^{-1}}\varphi}=\brmin{\intvl{\tau^{-1}s}\varphi}$, and,  by Corollary~\ref{c:cpparpar}, the desired equality follows.  Next suppose that $s\varphi<\varphi$.   Then, by Corollary~\ref{c:brcoset}, $s\varphi W_2\leq\varphi W_2$ and, given our hypothesis that $s\varphi W_2\geq \varphi W_2$,   we conclude that $s\varphi W_2=\varphi W_2$.    Put $\varphi':=s\varphi$.  Then $\varphi=s\varphi'>\varphi'=s\varphi$, and, by the first case,  we have $\brmin{W_1\intvl{\tau^{-1}}\varphi' W_2}=\brmin{W_1\intvl{\tau^{-1}s}s\varphi'W_2}$.   But, since $\varphi'W_2=\varphi W_2$,  this is exactly the desired equality.

Proof of~\eqref{i:takeaway4}: This is analogous to the proof of~\eqref{i:takeaway3}. By Corollary~\ref{c:cpparintvl}, $\brmin{\intvl{\tau^{-1}}\varphi}$ exists. First suppose that $s\tau<\tau$.    Then, by~Corollary~\ref{c:cpparintvl1}, $\brmin{\intvl{\tau^{-1}}\varphi}=\brmin{\intvl{\tau^{-1}}s\varphi}$, and,  by Corollary~\ref{c:cpparpar}, the desired equality follows.  Next suppose that $s\tau>\tau$.   Then, by Corollary~\ref{c:brcoset}, $s\tau W_1\geq\tau W_1$ and, given our hypothesis that $s\tau W_1\leq \tau W_1$,   we conclude that $s\tau W_1=\tau W_1$.    Put $\tau':=s\tau$.  Then $s\tau'=\tau<s\tau=\tau'$, and, by the first case,  we have $\brmin{W_1\intvl{\tau'^{-1}}\varphi W_2}=\brmin{W_1\intvl{\tau'^{-1}}s\varphi W_2}$.   But, since $\tau'W_1=\tau W_1$,  we have, by part~\eqref{i:takeaway1}, $W_1\intvl{\tau'^{-1}}=W_1\intvl{\tau^{-1}}$, and the desired equality follows.~\hfill\emyproof
%

\mysection{The KK filtration on concatenated LS paths}\mylabel{s:kosfiltpath}
\noindent
In this section,  we introduce the two key elements that underpin this entire paper, namely: \begin{itemize}
	\item the definition of Kostant-Kumar (KK) sets of concatenated Lakshmibai-Seshadri (LS) paths (Equation~\eqref{e:kosfiltpath}).
	\item the result that such a KK set is invariant under the action of root operators (Proposition~\ref{p:wind})
\end{itemize}
\ignore{ 
They are the only novelty in the proofs of our theorems.  The technical work in~\S\ref{s:extremal} is aimed solely at establishing them.   
In this section,  we define KK subsets of concatenated LS paths (Equation~\eqref{e:kosfiltpath}) and show that such a subset is invariant under the root operators on paths (Proposition~\ref{p:wind}).  As mentioned in the introduction,  the definition and the proposition are the key elements underpinning this paper.
}
%
The following notation remains fixed throughout this paper:  $\lieg$ denotes a symmetrizable Kac-Moody algebra;  $\lambda$ and $\mu$ are fixed dominant integral weights;    $W$ is the Weyl group and $W_\lambda$, $W_\mu$ are respectively the stabilizers in~$W$ of $\lambda$, $\mu$.

We assume familiarity with the basic notions and results of Littelmann's theory~\cite{litt:inv,litt:ann} of paths.  Let $\pathsl$, $\pathsm$ be respectively the sets of Lakshmibai-Seshadri (LS) paths of shape $\lambda$,~$\mu$.  Let $\pathsl\concat\pathsm:=\{\pi\concat\pi'\st \pi\in \pathsl,\ \pi'\in\pathsm\}$,  where $\concat$ denotes concatenation.


Recall that a path $\pi$ in~$\pathsl$ consists of a sequence $\tau_1>\tau_2>\ldots>\tau_r$ of elements in~$W/W_\lambda$ and a sequence $0=a_0<a_1<\ldots<a_{r-1}<a_r=1$ of rational numbers (subject to some integrality conditions as in~\cite[\S2]{litt:inv}, the details of which are not so relevant for the moment).   We call $\tau_1$ the {\em initial direction\/} and $\tau_r$ the {\em final direction\/} of~$\pi$.  
\mysubsection{Definition of a KK set in $\pathsl\concat\pathsm$}\mylabel{ss:defkosset}
Given a path $\pi\concat\pi'$ in $\pathsl\concat\pathsm$,  we define, using Corollary~\ref{c:takeaway}, part~\eqref{i:takeaway2}, its associated Weyl group element $\weyl(\pi\concat\pi')$ by:
\begin{equation}\label{e:defmain}
\weyl(\pi\concat\pi'):=\brmin{W_\lambda\intvl{\tau^{-1}}\varphi W_\mu}
\end{equation}
where $\tau$ and $\varphi$ are lifts in~$W$ respectively of the final direction of~$\pi$ and the initial direction of~$\pi'$.   Part~\eqref{i:takeaway1} of Corollary~\ref{c:takeaway} says that $\weyl(\pi\concat\pi')$ is independent of the choice of the lifts~$\tau$ and $\varphi$.

Given an element $\varphi=W_\lambda w W_\mu$ of the double coset space $\dcosetlm$,   we define the associated {\em KK set\/} by:
\begin{equation}\label{e:kosfiltpath}
	\plcpmphi
	:= \{\pi\concat\pi'\in\pathsl\concat\pathsm\st
	\weyl{(\pi\concat\pi')}\leq w\}
\end{equation}
The choice of the lift $w$ in~$W$ of~$\varphi$ does not matter (see Corollary~\ref{c:brcoset}),  and we often write $\plcpmw$ in place of $\plcpmphi$.

Clearly $\plcpmphi\subseteq\plcpmphip$ if $\varphi\leq\varphi'$.  Thus the KK sets form an increasing filtration of the space $\pathsl\concat\pathsl$ of concatenated LS paths, with underlying poset being the double coset space $\dcosetlm$ with its Bruhat order.  We call this the {\em KK filtration\/} on paths.

\bremark\mylabel{r:alternative}   For this remark alone,  we suppose that $W$ is finite.  Let $w_0$ be the longest element of~$W$.    For $w$ and $x$ in~$W$,  let $\circstar{w}{x}:=(ww_0\star xw_0)w_0$,  where $\star$ is Deodhar's operation on the Weyl group discussed in~\S\ref{r:deodhar1}.   
Let $\pi$ in~$\pathsl$ comprise the sequence $\tau_1>\tau_2>\ldots>\tau_r$ of elements in~$W/W_\lambda$ and the sequence $0=a_0<a_1<\ldots<a_{r-1}<a_r=1$ of rational numbers.   Let $\pidag$ be the path in~$\pathsl$ comprising $w_0\tau_r>\ldots>w_0\tau_1$ and $0=1-a_r<1-a_{r-1}<\ldots<1-a_1<1=1-a_0$.   Then
\[
	\weyl{(\pi\concat\pi')}= \circstar{\phi(\pidag)^{-1}}{\phi(\pi')}
	\]
	where $\phi(\eta)$ for an LS path~$\eta$ is the minimal lift in~$W$ of the initial direction of~$\eta$.
\eremark
\mysubsection{Stability of KK sets under root operators}\mylabel{ss:stabkosfilt}
For $\alpha$ a simple root,  let $e_\alpha$ and $f_\alpha$ be the root operators on paths as defined in~\cite{litt:ann}.   Although this definition differs from the earlier one in~\cite{litt:inv},  it is ``backwards-compatible'':  as explained in~\cite[Corollary~2 on page 512]{litt:ann}, the results of~\cite{litt:inv} are unaffected and we can freely quote them.

Let $\pi$ be a path in $\pathsl$.   Recall from~\cite{litt:inv}:
\begin{enumerate}[label=(\alph*)]
	\item\label{i:straight} The straight line path $\pil$ from the origin to $\lambda$ belongs to $\pathsl$.
	\item\label{i:plint} $\pi$ is piece-wise linear and its end point $\pi(1)$ is an integral weight.
	\item\label{i:eandf} For a simple root $\alpha$, if $e_\alpha(\pi)$ (respectively $f_\alpha(\pi)$) does not vanish,  then it belongs to $\pathsl$, its end point is $\pi(1)+\alpha$ (respectively $\pi(1)-\alpha)$,  and $f_\alpha(e_\alpha(\pi))$ (respectively $e_\alpha(f_\alpha(\pi))$) equals~$\pi$.
	\item\label{i:direction} Let $\alpha$ be a simple root and $\tau$ the initial (respectively, final) direction of $\pi$.   If $e_\alpha\pi$ does not vanish, then its initial (respectively, final) direction is either $\tau$ or $s_\alpha\tau$.  The same holds for $f_\alpha$ in place of~$e_\alpha$.
	\item\label{i:pifromtop} $\pi$ is obtained from $\pil$ by applying a suitable finite sequence of the root operators $f_\alpha$,  as $\alpha$ varies.  In particular, the end point of $\pi$ is of the form $\lambda-\kappa$ where $\kappa$ is a non-negative integral linear combination of the simple roots.
	\item\label{i:hslm} 
		Every value that is a local minimum of the function $\hapi(t):=\langle \pi(t),\alpha^\vee\rangle$ on $t\in[0,1]$ is an integer,  for every simple root $\alpha$.  (A value $\hapi(t_0)$ is called a {\em local minimum\/} if $\hapi(t_0)\leq \hapi(t)$ for $0\leq |t-t_0|<\epsilon$ for some $\epsilon>0$.)    This follows from the proof of~\cite[Lemma~4.5, part~(d)]{litt:ann} although the definition of local minimum there is less inclusive.
	\item\label{i:unique}  If $e_\alpha(\pi)$ vanishes for every simple root $\alpha$,  then $\pi=\pil$~\cite[Corollary in~\S3.5]{litt:inv}.    In particular, if (the image of) $\pi$ lies entirely in the dominant Weyl chamber,  then $\pi=\pil$.
\end{enumerate}
\begin{lemma}\mylabel{l:wind}
	Let $\pi\concat\pi'$ be a path in~$\pathsl\concat\pathsm$ and $\alpha$ a simple root.  
Then: \begin{enumerate}
	\item\label{i:hslmconcat}
	\newcounter{myenumikos}     
	\setcounter{myenumikos}{\theenumi}  
		Every local minimum value of the function $\ha^{\pi\concat\pi'}(t):=\langle (\pi\concat\pi')(t), \alpha^\vee\rangle$ is an integer.  
\end{enumerate}
	Suppose that $e_\alpha(\pi\concat\pi')$ does not vanish.  Then:
	\begin{enumerate}
\setcounter{enumi}{\value{myenumikos}}  
	\item\label{i:epipi'} $e_\alpha(\pi\concat\pi')$ equals either $e_\alpha\pi\concat\pi'$ or $\pi\concat e_\alpha\pi'$.
		\item\label{i:winv} $\weyl(\pi\concat\pi')=\weyl(e_\alpha(\pi\concat\pi'))$.
\end{enumerate}
\end{lemma}
\bmyproof
Statement~\eqref{i:hslmconcat} holds because $\pi(1)$ is an integral weight (item~\ref{i:plint} above) and local minima of both functions $\ha^\pi$ and $\ha^{\pi'}$ are integers (item~\ref{i:hslm} above).

Item~\eqref{i:epipi'} appears as~\cite[Lemma~2.7]{litt:ann}.  At any rate, it follows readily from the definition of~$e_\alpha$ once we know that the absolute minima of the functions~$\hapi$ and $\ha^{\pi'}$ are both integers,  which is guaranteed by item~\ref{i:hslm} above.

To prove~\eqref{i:winv}, let $\tau$ be the final direction of~$\pi$ and $\varphi$ the initial direction of~$\pi'$.     First suppose that $e_\alpha(\pi\concat\pi')=e_\alpha\pi\concat\pi'$.   By item~\ref{i:direction} above,  the final direction of $e_\alpha\pi$ is either $\tau$ or $s_\alpha\tau$.  If it is~$\tau$, then there is nothing for us to do.    In case it is $s_\alpha\tau$,   then, from the definition of~$e_\alpha$ and properties of~$\pi$ and $\pi'$, it follows that~$s_\alpha\tau<\tau$ and $s_\alpha\varphi\geq\varphi$.  The assertion now follows from part~\ref{i:takeaway3} of~Corollary~\ref{c:takeaway}.

Now suppose that $e_\alpha(\pi\concat\pi')=\pi\concat e_\alpha\pi'$.   By item~\ref{i:direction} above,  the initial direction of $e_\alpha\pi'$ is either $\varphi$ or~$s_\alpha\varphi$.  If it is $\varphi$, then there is nothing for us to do.    In case it is $s_\alpha\varphi$,   then, from the definition of~$e_\alpha$ and properties of~$\pi$ and $\pi'$, it follows that~$s_\alpha\varphi<\varphi$ and $s_\alpha\tau\leq\tau$.  The assertion now follows from part~\ref{i:takeaway4} of~Corollary~\ref{c:takeaway}.
\emyproof
\ignore{   
\begin{itemize}
	\item   The final direction of~$e_\alpha\pi$ is either~$\tau$ or~$s_\alpha\tau$; the initial direction of $e_\alpha\pi'$ is either $\varphi$ or $s_\alpha\varphi$.
	\item  The final direction of $e_\alpha\pi$ is $s_\alpha\tau$ only if $s_\alpha\tau<\tau$ and $s_\alpha\varphi\geq\varphi$.
	\item  The initial direction of $e_\alpha\pi'$ is $s_\alpha\varphi$ only if $s_\alpha\varphi<\varphi$ and $s_\alpha\tau\leq\tau$.
\end{itemize}
and that the final direction of $e_\alpha\pi$ is $s_\alpha\tau$.  
  Suppose that $e_\alpha(\pi\concat\pi')=\pi\concat e_\alpha\pi'$ and that the initial direction of $e_\alpha\pi'$ is $s_\alpha\varphi$.   Then the result follows from part~\eqref{i:takeaway4} of Corollary~\ref{c:takeaway}.
\emyproof
}   
\mysubsubsection{The equivalence relation on~$\pathsl\concat\pathsm$ defined by root operators}\mylabel{sss:sim}
Given $\pi\concat\pi'$ and $\sigma\concat\sigma'$ paths in $\pathsl\concat\pathsm$,  let us say $\pi\concat\pi'$ {\em is related to\/} $\sigma\concat\sigma'$ if $\pi\concat\pi'$ equals either $e_\alpha(\sigma\concat\sigma')$ or $f_\alpha(\sigma\concat\sigma')$ for some simple root~$\alpha$.  This relation is symmetric since $\pi\concat\pi'=e_\alpha(\sigma\concat\sigma')$ if and only if $f_\alpha(\pi\concat\pi')=\sigma\concat\sigma'$.   Denote by $\sim$ the reflexive and transitive closure of this relation (as we vary over all simple roots).   

As an immediate consequence of the item~\eqref{i:winv} of Lemma~\ref{l:wind}, we have:
	\begin{proposition}\mylabel{p:wind}  The association $\pi\concat\pi'\mapsto\weyl(\pi\concat\pi')$ is constant on equivalence classes of the equivalence relation~$\sim$.    In particular, for any $\varphi\in\dcosetlm$, the KK set $\plcpmphi$ is a union of such equivalence classes.\end{proposition}
\noindent
In other words, each KK set is stable under the root operators.   We will show in~\S\ref{s:pathkos} that a KK set provides a {\em path model\/} for the corresponding KK module.

\mysection{More Preliminaries}\mylabel{s:prelims}
\noindent
Notation is fixed as in~\S\ref{s:kosfiltpath}: $\lieg$ is a symmetrizable Kac-Moody algebra;  $\lambda$, $\mu$ are dominant integral weights; $W$ is the Weyl group; and $W_\lambda$, $W_\mu$ are the stabilisers in $W$ of $\lambda$, $\mu$ respectively.
\mysubsection{Geometric interpretation of minimal representatives in $\wlwwm$}\mylabel{ss:wlwwm}
\noindent
We now give a geometric interpretation of the unique minimal element in a given double coset in $\wlwwm$. 
The association $\wbar\leftrightarrow w\mu$  (for $w\in W$) gives a bijection of the coset space~$W/W_\mu$ with the set $W\mu$ of $W$-conjugates of~$\mu$.  We identify the sets $W/W_\mu$ and $W\mu$ via this bijection.   The double coset space~$\wlwwm$ may then be identified with the set of $\wl$-orbits of $W\mu$.
\bprop\mylabel{p:wlwwm}  Every $\wl$-orbit of the set $W\mu$ of $W$-conjugates of~$\mu$ contains a unique element $w\mu$ such that $\lambda+tw\mu$ is dominant for some real number $t>0$.
\eprop
\bmyproof  
Each such orbit contains a unique $w\mu$ that is $\wl$-dominant.  The $\wl$-dominance means precisely that $\langle w\mu,\alpha^\vee\rangle\geq0$ for every simple root $\alpha$ in~$\wl$.   It is easily verified that 
$w\mu$ has the desired property.    Conversely,  if $w\mu$  is not $\wl$-dominant,  then $\langle w\mu,\alpha^\vee\rangle <0$ for some simple root~$\alpha$ in~$\wl$,  and so $\langle\lambda+tw\mu,\alpha^\vee\rangle=t\langle w\mu,\alpha^\vee\rangle<0$ for $t>0$.
\emyproof

The double coset space $\wlwwm$ may thus be identified with the set of those Weyl conjugates $w\mu$ of~$\mu$ such that $\lambda+tw\mu$ is dominant for some positive $t$.   We illustrate this by means of an example.
\bexample\mylabel{x:b2coset}   Let $\lieg$ be of type~$B_2$.
Let $e_1$ and $e_2$ be the standard basis vectors in~$\mathbb{R}^2$ with its standard inner product.    We may take
$\alpha_1:=e_1$ and $\alpha_2:= e_2-e_1$ to be the simple roots.   Then the set of all positive roots is $\{\alpha_2,\alpha_1,\alpha_2+\alpha_1,\alpha_2+2\alpha_1\}$,  and the fundamental weights are $\varpi_1 = \frac{1}{2}(\epsilon_1+\epsilon_2)$ and $\varpi_2=\epsilon_2$.  The Weyl group consists of~$8$ elements:
\[ W = \{ 1,s_1,s_2,s_1s_2,s_2s_1,s_2s_1s_2,s_1s_2s_1,s_1s_2s_1s_2=s_2s_1s_2s_1  \} \]
where $s_1$ and $s_2$ are the reflections in the hyperplanes perpendicular to $\alpha_1$ and $\alpha_2$ respectively.
The shaded portion in the figure is the dominant Weyl chamber. 

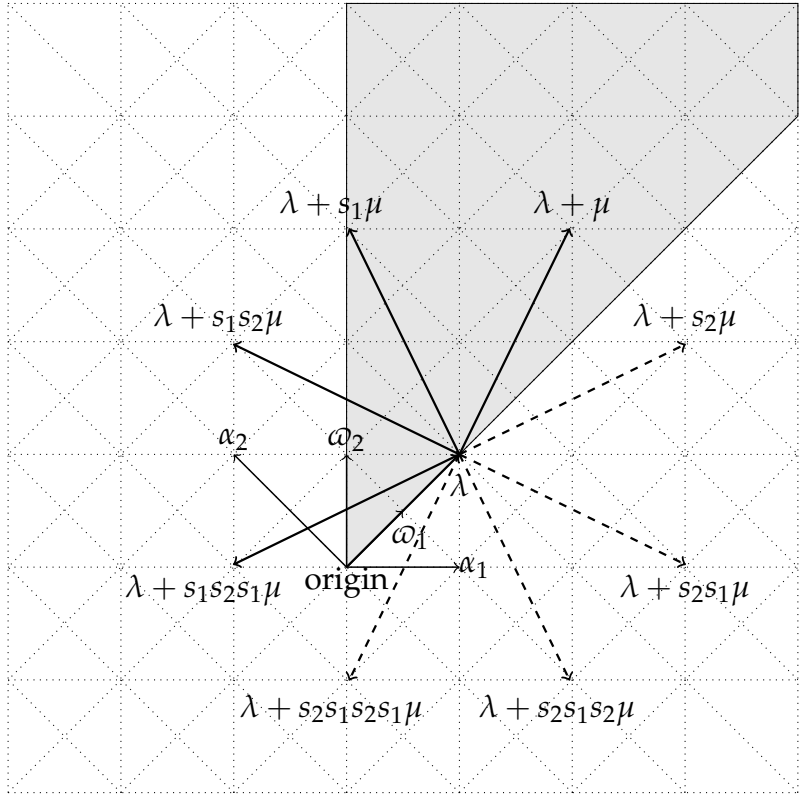
\begin{figure}
\centering
\begin{tikzpicture}

\draw[fill = gray!20]  (6,4.5) -- (12,10.5) -- (12,12) -- (6,12) -- cycle;

\draw[scale = 1.5,style = dotted] (1,1) grid (8,8);
\draw [style = dotted](1.5,1.5)--(12,12);
\draw [style = dotted](3,1.5)--(12,10.5)--(10.5,12)--(1.5,3)--cycle;
\draw [style = dotted](4.5,1.5)--(12,9)--(9,12)--(1.5,4.5)--cycle;
\draw [style = dotted](6,1.5)--(12,7.5)--(7.5,12)--(1.5,6)--cycle;
\draw [style = dotted](7.5,1.5)--(12,6)--(6,12)--(1.5,7.5)--cycle;
\draw [style = dotted](9,1.5)--(12,4.5)--(4.5,12)--(1.5,9)--cycle;

	\draw [style = dotted](10.5,1.5)--(1.5,10.5);
\draw [style = dotted](12,1.5)--(1.5,12);
\draw [style = dotted](12,3)--(3,12);



\draw[->,line width = .2 mm] (6,4.5)--(4.5,6);
\node at (4.5,6.2) {$\alpha_2$};
\draw[->,line width = .2 mm] (6,4.5)--(7.5,4.5);
\node at (7.7,4.5) {$\alpha_1$};
\draw[->,line width = .2 mm] (6,4.5)--(6,6);
\node at (6,6.2) {$\varpi_2$};
\draw[->,line width = .2 mm] (6,4.5)--(6.75,5.25);
\node at (6.85,4.9) {$\varpi_1$};

\draw[->,line width = .3 mm] (6,4.5)--(7.5,6);
\node at (7.5,5.6) {$\lambda$};
	\node at (6,4.5-0.2) {$\textup{origin}$};

\node at (9,9.3) {$\lambda+\mu$};
\node at (10.5,7.8) {$\lambda+s_2\mu$};
\node at (10.5,4.2) {$\lambda+s_2s_1\mu$};
\node at (5.8,9.3) {$\lambda+s_1\mu$};
\node at (5.8,2.6) {$\lambda+s_2s_1s_2s_1\mu$};
\node [rotate = 0]at (8.8,2.6) {$\lambda+s_2s_1s_2\mu$};
\node at (4.3,7.8) {$\lambda+s_1s_2\mu$};
\node at (4.1,4.2) {$\lambda+s_1s_2s_1\mu$};
\draw[->,line width = .3 mm] (7.5,6)--++(64:3.35);
\draw[->,line width = .3 mm,dashed] (7.5,6)--++(-64:3.35);
\draw[->,line width = .3 mm] (7.5,6)--++(116:3.35);
\draw[->,line width = .3 mm,style = dashed] (7.5,6)--++(-26:3.35);
\draw[->,line width = .3 mm,style = dashed] (7.5,6)--++(26:3.35);
\draw[->,line width = .3 mm,style = dashed] (7.5,6)--++(-116:3.35);
\draw[->,line width = .3 mm] (7.5,6)--++(154:3.35);
\draw[->,line width = .3 mm] (7.5,6)--++(-154:3.35);






\end{tikzpicture}
\caption{Illustration of Propositions~\ref{p:wlwwm},~\ref{p:wlwwm2}; see Example~\ref{x:b2coset}}
\label{f:b2coset}
\end{figure}

Take $\lambda=2\omega_1$ and $\mu = 2\omega_1+\omega_2$.
 The stabilizers of $\lambda$ and $\mu$ are respectively:   $W_\lambda=\{1,s_2\}$ and $W_\mu=\{1\}$.
The set of double cosets $\wlwwm$ is:
\[ \{ {\{1,s_2\}},{\{s_1,s_2s_1 \}},{\{ s_1s_2,s_2s_1s_2 \}},{\{ s_1s_2s_1,s_2s_1s_2s_1 \}}  \} \]
 As is clear from Figure~\ref{f:b2coset},   $\mu$, $s_1\mu$, $s_1s_2\mu$, and $s_1s_2s_1\mu$ are all the conjugates of $\mu$ for which the line segment joining $\lambda$ to the conjugate lies for some positive distance in the dominant Weyl chamber.  
\eexample
\bprop\mylabel{p:wlwwm2} Given a double coset in~$\wlwwm$,  let $w$ be the unique minimal element in it with respect to the Bruhat order (as guaranteed by Corollary~\ref{c:cpparpar}).  Then $w\mu$ is such that $\lambda+tw\mu$ is dominant for all small positive~$t$.
\eprop
\bmyproof  From the proof of Proposition~\ref{p:wlwwm},  it is enough to show that $w\mu$ is $\wl$-dominant.    Suppose that this is not so.   Then there exists simple root $\alpha$ with $s_\alpha$ in~$\wl$ such that $\langle w\mu,\alpha^\vee\rangle<0$.   We then have $s_\alpha w<w$, which contradicts the hypothesis that $w$ is the minimal element in its double coset.
\emyproof

\mysubsection{Two key propositions}\mylabel{ss:2keyp}\mylabel{ss:keyp}
\noindent
An LS path $\pi$ of shape~$\mu$ is said to be {\em $\lambda$-dominant\/} if $\lambda+\pi(t)$ belongs to the dominant Weyl chamber for every $t\in[0,1]$.   The set of $\lambda$-dominant paths of shape~$\mu$ is denoted by~$\pathsml$.       For $w$ an element of the Weyl group,  $\pathsmlw$ denotes the elements of $\pathsml$ whose initial direction is $\leq wW_\mu$.
\bprop\mylabel{p:dompath}
Let $\theta$ be a path in $\pathsl\concat\pathsm$.    Then there exists a unique path~$\eta$ in the equivalence class (of the relation~$\sim$ defined in~\S\ref{sss:sim}) containing $\theta$ such that $e_\alpha\eta$ vanishes on $\eta$ for all simple roots~$\alpha$.    Moreover,  $\eta$ has the following properties:
\begin{enumerate}
	\item\label{i:etadom} $\eta$ lies entirely in the dominant Weyl chamber.
	\item\label{i:piinpml} $\eta=\pi_\lambda\concat\pi$ for some $\pi$ in~$\pathsml$.
	\item\label{i:initv} $\weyl(\theta)=\weyl(\eta)=v$ where $v$ is minimal in the Weyl group such that $v\mu$ is the initial direction of~$\pi$.   In particular, if $\theta\in\plcpmw$ (for some $w$ in~$W$),  then $\pi\in\pathsmlw$.
\end{enumerate}
\eprop
\bmyproof  
For the existence of $\eta$, there is the following standard argument. Construct by induction a sequence $\theta_0$, $\theta_1$, \ldots\ of elements in the equivalence class of~$\theta$ as follows.   Choose $\theta_0$ to be $\theta$.  Given $\theta_i$, if $e_\alpha\theta_i$ vanishes for all simple roots $\alpha$,  then set $\eta=\theta_i$ and we are done.   If not,  then choose $\alpha$ simple root arbitrarily such that $e_\alpha\theta_i$ does not vanish and put $\theta_{i+1}=e_\alpha\theta_i$.    By induction $\theta_{i+1}$ belongs to the equivalence class of~$\theta$.   We will eventually find an $\eta$ this way,  for this process must terminate at some point.   In fact,  the length of the sequence is bounded by the sum of the coefficients of~$\kappa$ where $\kappa$ is the non-negative integral linear combination of the simple roots such that the end point of~$\theta$ equals $\lambda+\mu-\kappa$.

Since $e_\alpha\eta$ vanishes for all simple~$\alpha$ and since the absolute minimum of the function $\ha^\eta(t)$ is an integer for every simple~$\alpha$ (see item~\eqref{i:hslmconcat} in~Lemma~\ref{l:wind} above),   it follows from the definition of~$e_\alpha$ that $\eta$ lies entirely in the dominant Weyl chamber.

The uniqueness of $\eta$ now follows from~\cite[Corollary~1\,(b) of~\S7]{litt:ann}.

Write $\eta=\zeta\concat\pi$ with $\zeta\in\pathsl$ and $\pi\in\pathsm$.  Since $\eta$ lies entirely in the dominant Weyl chamber,  clearly so does~$\zeta$.   Thus~$\zeta=\pi_\lambda$ by item~\ref{i:unique} in~\S\ref{ss:stabkosfilt} above, and $\pi$ belongs to~$\pathsml$.

The equality $\weyl(\theta)=\weyl(\eta)$ follows from Proposition~\ref{p:wind}.  
Since $\eta$ lies entirely in the dominant Weyl chamber,   it follows that $\lambda+tv\mu$ is dominant for sufficiently small $t\geq0$.  By Proposition~\ref{p:wlwwm2}, the unique minimal element of $W_\lambda v W_\mu$ lies in $vW_\mu$ and hence equals $v$.   But $w(\eta)=\brmin{W_\lambda vW_\mu}$ by its definition.  \emyproof
\bprop\mylabel{p:dompatheq}
With notation as in Proposition~\ref{p:dompath},  write $\theta=\pi_1\concat\pi_2$.   Then the following conditions are equivalent:
\begin{enumerate}
	\item\label{i:pdeq1} $\eta=\pi_\lambda\concat\pi_\mu$ (that is, $\pi=\pi_\mu$)
	\item\label{i:pdeq2} $\weyl(\theta)=\textup{identity}$
	\item\label{i:pdeq3} there exist $\taut$ and $\phit$ in~$W$ such that $\taut\geq\phit$, $\taut W_\lambda$ is the final direction of $\pi_1$, and $\phit W_\mu$ is the initial direction of~$\pi_2$.
\end{enumerate}
\eprop
\bmyproof
Let $v$ be as defined in item~\eqref{i:initv} of Proposition~\ref{p:dompath}.     Observe that condition~\eqref{i:pdeq1} is equivalent to saying that $v$ is identity.   Since $\weyl(\theta)=v$ by Proposition~\ref{p:dompath}\,\eqref{i:initv},   we have \eqref{i:pdeq1}$\Leftrightarrow$\eqref{i:pdeq2}.

\eqref{i:pdeq2}$\Rightarrow$\eqref{i:pdeq3}:   Let $\tau$ and $\varphi$ be arbitrary elements in~$W$ such that $\tau W_\lambda$ is the final direction of~$\pi_1$ and $\varphi W_\mu$ is the initial direction of~$\pi_2$.   Condition~\eqref{i:pdeq2} says that $\brmin{W_\lambda \intvl{\tau^{-1}}\varphi W_\mu}$ equals identity.    Let $u\in W_\lambda$,  $\sigma\leq\tau$ in $W$, and $v\in W_\mu$ be such that $u^{-1}\sigma^{-1}\varphi v=\textup{identity}$,  or $\varphi v=\sigma u$.    We have $\sigma W_\lambda\leq \tau W_\lambda$ (by Corollary~\ref{c:brcoset}).   By Proposition~\ref{p:brcoset},  there exists $\taut$ in~$\tau W_\lambda$ such that $\sigma u\leq \taut$.   Taking $\phit=\varphi v=\sigma u$, \eqref{i:pdeq3} is proved.

\eqref{i:pdeq3}$\Rightarrow$\eqref{i:pdeq2}:   Since $\phit\leq\taut$,  it follows that $\phit^{-1}$ belongs to $\intvl{\taut^{-1}}$.   This implies that $\brmin{W_\lambda \intvl{\taut^{-1}}\phit W_\mu}$ equals identity.  But $\weyl(\theta)=\brmin{W_\lambda \intvl{\taut^{-1}}\phit W_\mu}$ by definition.  \emyproof
\mysubsection{Extremal paths}\mylabel{ss:extpath}
Let $\theta$ be a path in $\pathsl\concat\pathsm$ and let $\eta$ be as in Proposition~\ref{p:dompath} above.   Following Montagard~\cite{montagard}, we call $\theta$ {\em extremal\/} if the dominant Weyl conjugate~$\overline{\theta(1)}$ of the end point of~$\theta$ equals the end point~$\lambda+\pi(1)$ of~$\eta$.

The following observation~\cite[Theorem~2.2~(i)]{montagard} applied to the path $\pi_\lambda\concat\pi_{u\mu}$ is already used in Littelmann's proof~\cite[\S7]{litt:inv} of the PRV conjecture (here $u$ denotes an element of~$W$ and $\pi_{u\mu}$ the straight line path to the extremal weight $u\mu$ in~$V_\mu$):
\bprop\mylabel{p:extremal}    If a path $\theta\in\pathsl\concat\pathsm$ lies entirely in the dominant Weyl chamber except perhaps for a portion of its last straight line segment,  then $\theta$ is extremal in the above sense.
\eprop

\mysection{The KK (sub)modules of~$V_\lambda\tensor V_\mu$}\mylabel{s:kosdef}
\noindent
In this section we recall the definition of Kostant-Kumar (KK) modules and two basic results about them (Propositions~\ref{p:koswind} and~\ref{p:kosleq}). 

Let $\lieg$ be a symmetrizable Kac-Moody algebra. Let $\lambda$, $\mu$ be dominant integral weights.  Let $\Vlambda$, $\Vmu$ be the irreducible integrable $\lieg$-modules with respective highest weights $\lambda$, $\mu$.  Let $\wl$, $\wm$ be the respective stabilizers of $\lambda$, $\mu$ in the Weyl group~$W$.
\mysubsection{Filtration by KK modules of $\Vlambda\tensor\Vmu$}\mylabel{ss:kosfilt}
Fix an element $w$ of the Weyl group.  Let $\vlambda$ be a highest weight vector in~$\Vlambda$. Let $\vwmu$ be a non-zero vector in the (one-dimensional) weight space $V_{\wmu}$ of weight $w\mu$ in~$\Vmu$.  The {\em Kostant-Kumar module\/}, or simply {\em KK module\/}, $\koslwm$ is defined to be the cyclic submodule of the tensor product $\Vlambda\tensor\Vmu$ generated by $\vlambda\tensor\vwmu$:
\begin{equation}\label{e:kosmoddef}
	\koslwm := U\lieg\left(\vlambda\tensor\vwmu\right)
\end{equation}
where $U\lieg$ denotes the universal enveloping algebra of~$\lieg$.
\begin{proposition}\mylabel{p:koswind}
	Let $u$ and $w$ be elements of the Weyl group such that $W_\lambda uW_\mu=W_\lambda wW_\mu$.   Then $\koslum=\koslwm$.
\end{proposition}
\bmyproof
For the proof,  we will first recall a basic result from~\cite[\S3.8]{kac}. Let $V$ be an integrable representation of~$\lieg$.   For a simple reflection $s$,  there is a corresponding linear automorphism $\rv$ of~$V$ (defined in~\cite[Lemma~3.8]{kac}) such that:  
\begin{enumerate}
	\item\label{i:one} For $v\in V$ a weight vector of weight~$\eta$,  $\rv(v)$ is a weight vector of weight~$s(\eta)$.
	\item\label{i:two} $r^{V\tensor V'}=\rv\tensor r^{V'}$ for $V'$ an integrable representation.
	\item\label{i:three} For $v\in V$,  there exists $x_v\in U\lieg$ such that $\rv(v)=x_v(v)$.
\end{enumerate}

It suffices to show that $\koslum\subseteq\koslwm$, for then the other containment also holds by interchanging the roles of $u$ and $w$.  Write $u=\tau w\varphi$ with $\tau\in W_\lambda$ and $\varphi\in W_\mu$. We have $u\mu=\tau w\varphi\mu=\tau w \mu$.   Let $\tau=s_{i_1}\cdots s_{i_k}$ be a reduced expression for $\tau$.   Note that all $s_{i_j}$ belong to $W_\lambda$.   Consider the operator $r_{i_1}\cdots r_{i_k}$ on $\Vlambda\tensor\Vmu$ where $r_{i_j}=r_{i_j}^{\Vlambda\tensor\Vmu}
$ is the linear automorphism corresponding to $s_{i_j}$ (as recalled above).

On the one hand,  by properties~\eqref{i:one} and \eqref{i:two} above,   we have:
\[ r_{i_k}(\vlambda\tensor\vwmu) = r_{i_k}^{\Vlambda}(\vlambda)\tensor r_{i_k}^{\Vmu}(\vwmu)= c\cdot v_{s_{i_k}\lambda}\tensor v_{s_{i_k}w\mu} = c\cdot \vlambda\tensor v_{s_{i_k}w\mu} \]
where $c$ is a non-zero scalar.    By a chain of similar calculations,  we get
\begin{equation}\label{e:kosdef1} r_{i_1}\cdots r_{i_k}(\vlambda\tensor\vwmu) = c'\cdot \vlambda\tensor v_{\tau w\mu}=c'\cdot \vlambda\tensor v_{u\mu} \end{equation}
where $c'$ is a non-zero scalar.

On the other hand,  by property~\eqref{i:three}, there exist elements $x_{i_1}$, \ldots, $x_{i_k}$ of $U\lieg$ such that 
\begin{equation}\label{e:kosdef2}
r_{i_1}\cdots r_{i_k} (v_\lambda\tensor\vwmu)= x_{i_1}\cdots x_{i_k} (v_\lambda\tensor\vwmu)
\end{equation}
From \eqref{e:kosdef1} and \eqref{e:kosdef2} we get
\[ \vlambda\tensor v_{u\mu} = c'^{-1}\cdot r_{i_1}\cdots r_{i_k}(\vlambda\tensor\vwmu) =
 c'^{-1}\cdot x_{i_1}\cdots x_{i_k}(\vlambda\tensor\vwmu) \]
 and thus $\koslum=U\lieg(\vlambda\tensor\vumu)\subseteq U\lieg(\vlambda\tensor\vwmu)=\koslwm$.
\emyproof
\begin{proposition}\mylabel{p:kosleq} 
	For elements $u$ and $w$ of the Weyl group~$W$ such that $W_\lambda u W_\mu\leq W_\lambda w W_\mu$ in the Bruhat order on $\dcosetlm$ (see~\S\ref{ss:brcoset}),  we have $\koslum\subseteq\koslwm$.
\end{proposition}
\bmyproof
By Proposition~\ref{p:koswind}, we may assume~$u=\brmin{W_\lambda uW_\mu}$ and $w=\brmin{W_\lambda wW_\mu}$, so that $u\leq w$.  Let $U\lieb(\vwmu)$ be the Demazure module generated by $\vwmu$.     Since $u\leq w$,  we have $\vumu\in U\lieb(\vwmu)$.   Thus $\vlambda\tensor\vumu\in U\lieb(\vlambda\tensor\vwmu)\subseteq U\lieg(\vlambda\tensor\vwmu)$,  and $U\lieg(\vlambda\tensor\vumu)\subseteq U\lieg(\vlambda\tensor\vwmu)$.
\emyproof
\bremark  
 The KK module $\koslonem = U\lieg(\vlambda\tensor\vmu)$ corresponding to the identity element~$1$ of the Weyl group is the copy of the irreducible representation $V_{\lambda+\mu}$ in $\Vlambda\tensor\Vmu$.    When $\lieg$ is of finite type,  the KK module $\koslwnotm$ corresponding to the longest element $w_0$ of the Weyl group is the whole tensor product $V_\lambda\tensor V_\mu$.   Indeed,   letting $\bplus$ and $\bminus$  denote respectively the positive and negative Borel subalgebras,   we have 
\begin{align*}
	U\lieg(v_\lambda\tensor\vwnotmu)=U\bminus\cdot U\bplus (\vlambda\tensor\vwnotmu)&=U\bminus (\vlambda\tensor U\bplus \vwnotmu)\\ &=U\bminus(\vlambda\tensor\Vmu)=(U\bminus\vlambda)\tensor\Vmu=\Vlambda\tensor\Vmu
\end{align*}
\eremark

\mysection{Recall of a decomposition rule for Kostant-Kumar (KK) modules}\mylabel{s:decompkos}
\noindent
The decomposition rule (Theorem~\ref{t:decompkos} below) that gives the break up of a KK module into a direct sum of irreducibles is well known.  For example, at least in the case when $\lieg$ is symmetric (i.e., has a symmetric generalized Cartan matrix),  it  follows immediately from Joseph's results~\cite[Theorems~5.25,~5.22]{joseph:demflag}.    Our purpose in this section is to state the theorem and also give, for the sake of readability and completeness, a proof in the case when $\lieg$ is of finite type.  The restrictive hypothesis on~$\lieg$ (namely that it be of finite type or symmetric) is imposed only due to the use of a positivity result of Lusztig~\cite[22.1.7]{lusztig:book93} by Joseph in~\cite{joseph:demflag} and is possibly not required: see~\cite[\S1.4]{joseph:demflag}.
\begin{theorem}\textup{(Joseph~\cite{joseph:demflag})}\mylabel{t:decompkos} 
Let $\lieg$ be a symmetrizable Kac-Moody Lie algebra that is either of finite type or symmetric.   Let $\lambda, \mu$ be dominant integral weights and $w$ an element of the Weyl group.  Then the decomposition of the KK module~$\klwm$ as a direct sum of irreducible $\lieg$-modules is given by
	\begin{equation}\label{e:decompkos}
		\klwm =
		\bigoplus V_{\lambda+\pi(1)}
		\quad\textup{where the sum is over $\pi\in\pathsmlw$}
	\end{equation}
where~$\pathsmlw$ denotes the set of $\lambda$-dominant LS paths of shape~$\mu$ with initial direction $\leq wW_\mu$ ($\lambda$-dominance of a path is defined in~\S\ref{ss:keyp}).
\end{theorem}
\noindent
In the case when $\lieg$ is symmetric,   the theorem follows from Joseph's results as already indicated (see also Naoi \cite[Remark 2.12]{naoi}).  In the finite type case,  a proof is recorded below (\S\ref{ss:pfdecompkos}).  
For this proof,  we need a result of Lakshmibai, Littelmann, and Magyar~\cite{lak-litt-mag},  which is a combinatorial analogue of the existence of ``excellent filtrations'', a la Joseph \cite{joseph:dcf}, Mathieu \cite{mathieu:m1,mathieu:m2}, Polo \cite{polo}, and van der Kallen \cite{vdk}.   We first recall this result.
\mysubsection{A result of Lakshmibai-Littelmann-Magyar}\mylabel{ss:llm}
%
In order to state the result,  we introduce some notation.
The term {\em path} in this section means a piecewise linear path whose endpoint lies in the weight lattice (for instance, a concatenation of LS paths of various shapes). 
Let $\paths$ be a set of paths. We define its {\em character\/}, denoted $\character{\paths}$,  by:
$ \character{\paths} := \sum_{\eta \in \paths} \expo{\eta(1)} $.
If $\pi$ is any path, we let $\pi \concat \paths$ denote the set of paths $\{\pi \concat \eta: \eta \in \paths\}$.
Suppose $\pi$ is a path such that $\pi(t)$ belongs to the dominant Weyl chamber for all $t \in [0,1]$. 
Fix a reduced word
$w = s_{\beta_1} s_{\beta_2} \cdots s_{\beta_k}$:  here $\beta_i$ are simple roots.  Define
\[ C(\pi,w) := \{ f_{\beta_1}^{n_1} f_{\beta_2}^{n_2} \cdots  f_{\beta_k}^{n_k} \,\pi: \, n_i \geq 0 \text{ for all } i\}\]
This set is independent of the reduced word chosen, and has character:
\[ \character{C(\pi,w)}  = \demw(\expo{\pi(1)}) \]
where $\demw$ is the Demazure operator corresponding to $w$ (see, e.g., \cite[\S5.1]{litt:inv}). Further, when $\pi$ is the straight line path $\pi_\mu$, we have $C(\pi_\mu, w) = \pathsmw$, the set of LS paths of shape~$\mu$ with initial direction $\leq wW_\mu$. The following key result appears in \cite[Proposition 12]{lak-litt-mag} (see also~\cite[Theorem in~\S2.11, also~\S3.5]{joseph:crystal}).

\bprop \mylabel{p:llm} 
With notation as above,   there exists a Weyl group valued function $\pi\mapsto w(\pi)$ on~$\pathsmlw$ such that
\begin{equation} \label{e:llm}
\pi_\lambda \concat C(\pi_\mu, w) = \bigsqcup_{\pi \in \pathsmlw} C(\pi_\lambda \concat \pi, w(\pi))
\end{equation}
(The precise form of the function $\pi\mapsto w(\pi)$ is immaterial for our purposes.)
\eprop

Computing characters of both sides in \eqref{e:llm}, we obtain:
\begin{equation}\label{e:llm-char-form}
\character{\pi_\lambda \concat C(\pi_\mu, w)} = \sum_{\pi \in \pathsmlw} \demazure_{w(\pi)} (\expo{\lambda + \pi(1)})
\end{equation}

\mysubsection{Proof of Theorem~\ref{t:decompkos} for $\lieg$ of finite type}\mylabel{ss:pfdecompkos}
By a result of Kumar~\cite[Theorem~2.14]{kumar:prv},  the character of the KK module $\klwm$ is given by
\begin{equation}\label{e:charklwm}
	\character{\klwm}=\demwnot(\expo{\lambda}\cdot\demw(\expo{\mu}))
\end{equation}
where $w_0$ is the longest Weyl group element. 
%
%
Since $\demw(\expo{\mu})$ is the character of ${C(\pi_\mu, w)}$, we obtain \begin{equation}\label{e:charklwm:1}\character{\klwm} =\demwnot(\character{\pi_\lambda \concat C(\pi_\mu, w)})\end{equation} 
	Substituting from~\eqref{e:llm-char-form} into~\eqref{e:charklwm}, we obtain:
	\begin{equation*}
\character{\klwm} =  \sum_{\pi \in \pathsmlw} \demwnot \demazure_{w(\pi)} (\expo{\lambda + \pi(1)}) 
=  \sum_{\pi \in \pathsmlw} \demwnot (\expo{\lambda + \pi(1)})
\end{equation*}
since $\demwnot \demazure_\sigma = \demwnot$ for all $\sigma \in W$. This latter fact follows from the following well-known property of the Demazure operators: if $\alpha$ is a simple root, then $\demazure_{s_\alpha} \demw$ equals $\demazure_{s_\alpha w}$ or $\demw$ according as $s_\alpha w$ is  $>w$ or $<w$.

But now $\demwnot (\expo{\lambda + \pi(1)})$ is the character of the $\lieg$-module $V_{\lambda+\pi(1)}$ (by the Demazure character formula applied to $w_0$). Thus the modules on both sides of \eqref{e:decompkos} have the same character, and the proof is complete.\hfill$\Box$
\bexamplenobox\mylabel{x:b2}    Consider the situation of Example~\ref{x:b2coset} and Figure~\ref{f:b2coset}.   Figure~\ref{f:b2} depicts all $\lambda$-dominant paths of shape~$\mu$, colour coded by their initial directions as follows:
\[
	\begin{array}{c|c}
		\textup{Initial direction} & \textup{Colour coding}\\
		\hline
		\textup{identity} & \textup{violet}\\ 
		s_1 & \textup{red}\\
		s_1s_2 &\textup{cyan}\\
		s_1s_2s_1 &\textup{orange}\\
	\end{array}
	\]
The KK modules decompose as follows:
\begin{itemize}
\item $K(\lambda,1,\mu) = V_{\lambda+\mu}.$
\item $K(\lambda,s_1,\mu) = V_{\lambda+\mu} \oplus V_{2\varpi_1+2\varpi_2} \oplus V_{3\varpi_2}.$
\item $K(\lambda,s_1s_2,\mu) = V_{\lambda+\mu} \oplus V_{2\varpi_1+2\varpi_2} \oplus V_{3\varpi_2} \oplus V_{4\varpi_1} \oplus V_{2\varpi_1+\varpi_2}.$
\item $K(\lambda,s_1s_2s_1,\mu) = V_{\lambda+\mu} \oplus V_{2\varpi_1+2\varpi_2} \oplus V_{3\varpi_2}\oplus V_{4\varpi_1} \oplus V_{2\varpi_1+\varpi_2} \oplus V_{2\varpi_1+\varpi_2}\oplus V_{\varpi_2} \oplus V_{2\varpi_2}.$
\end{itemize} 
\eexamplenobox
\begin{figure}
\centering
\begin{tikzpicture}
\draw[fill = gray!20]  (6,4.5) -- (12,10.5) -- (12,12) -- (6,12) -- cycle;

\draw[scale = 1.5,style = dotted] (2,2) grid (8,8);
\draw [style = dotted](3,3)--(12,12);
\draw [style = dotted](4.5,3)--(12,10.5)--(10.5,12)--(3,4.5);
\draw [style = dotted](6,3)--(12,9)--(9,12)--(3,6);
\draw [style = dotted](7.5,3)--(12,7.5)--(7.5,12)--(3,7.5);
\draw [style = dotted](9,3)--(12,6)--(6,12)--(3,9);
\draw [style = dotted](10.5,3)--(12,4.5)--(4.5,12)--(3,10.5);

\draw [style = dotted](4.5,3)--(3,4.5);
\draw [style = dotted](6,3)--(3,6);
\draw [style = dotted](7.5,3)--(3,7.5);
\draw [style = dotted](9,3)--(3,9);
\draw [style = dotted](10.5,3)--(3,10.5);
\draw [style = dotted](12,3)--(3,12);



\draw[->,line width = .5 mm] (6,4.5)--(4.5,6);
\node at (4.5,6.2) {$\alpha_2$};
\draw[->,line width = .5 mm] (6,4.5)--(7.5,4.5);
\node at (7.8,4.5) {$\alpha_1$};
\draw[->,line width = .5 mm] (6,4.5)--(6,6);
\node at (6,6.2) {$\varpi_2$};
\draw[->,line width = .5 mm] (6,4.5)--(6.75,5.25);
\node at (6.85,4.9) {$\varpi_1$};

\draw[->,line width = .3 mm,blue] (6,4.5)--(7.5,6);
\node at (7.6,5.8) {$\lambda$};
\draw[->,line width = .3 mm,violet] (7.5,6)--++(64:3.35);
\node at (9,9.2) {$\lambda+\mu$};
\node at (9.5,7.5) {$4\varpi_1$};
\node [rotate = 45] at (7.5,7.8) {$2\varpi_1+\varpi_2$};
\node  at (5.5,7.8) {$2\varpi_2$};
\node [rotate = 45]at (7.5,9.2) {$2\varpi_2+2\varpi_1$};
\node [rotate = 45]at (6,9.2) {$3\varpi_2$};
\draw[->,line width = .3 mm,red] (7.5,6)--++(116:3.35/2)--++(64:3.35/2);
\draw[->,line width = .3 mm,red] (7.5,6)--++(116:3.35);


\draw[->,line width = .3 mm,cyan] (7.5,6)--++(154:3.35/4)--++(26:3*3.35/4);
\draw[->,line width = .3 mm,cyan] (7.5,6)--++(154:3.35/2)--++(26:3.35/2);


\draw[->,line width = .3 mm,style,orange] (7.5,6)--++(206:3.35/4)--++(-26:3.35/3-3.35/4)--++(116:2*3.35/3);
\draw[->,line width = .3 mm,orange] (7.5,6)--++(206:3.35/4)--++(-26:3.35/3-3.35/4)--++(116:3.35/2-3.35/3)--++(64:3.35/2);
\draw[->,line width = .3 mm,orange] (7.5,6)--++(206:3.35/2)--++(-26:2*3.35/3-3.35/2)--++(116:3.35/3);



\end{tikzpicture}
\caption{Decomposition of KK modules of~$V_{2\varpi_1}\tensor V_{2\varpi_1+\varpi_2}$ for $\lieg$ of type~$B_2$} \label{f:b2}
\end{figure}

\mysection{A path model for Kostant-Kumar (KK) modules}\mylabel{s:pathkos}
\noindent
We deduce a path model for KK modules by combining the decomposition rule (Theorem~\ref{t:decompkos}) with the invariance under the root operators (Proposition~\ref{p:wind})   of the association~\eqref{e:defmain} of the Weyl group element~$\weyl(\pi\concat\pi')$ to a concatenation~$\pi\concat\pi'$ of two LS paths.     The restriction on~$\lieg$ (namely, that it be of finite type or symmetric) in the theorem is only because of the use of the decomposition rule,  and is possibly not required.
\begin{theorem}\mylabel{t:pathkos}
	Let $\lieg$ be a symmetrizable Kac-Moody algebra that is either of finite type or symmetric (as in Theorem~\ref{t:decompkos}).
	Let $\lamdba$, $\mu$ be dominant integral weights and $w$ an element of the Weyl group.  Let $\pathsl$ and $\pathsm$ respectively be the sets of Lakshmibai-Seshadri paths of shapes $\lambda$ and $\mu$.   For $\pi\in\pathsl$ and $\pi\in\pathsm$,  let $\wppprime$ be the Weyl group element associated as in~\eqref{e:defmain} in~\S\ref{ss:defkosset} to the concatenated path $\pi\concat\pi'$.    Then the KK set
	\[
		\plcpmw =\{\pi\concat\pi'\st \pi\in\pathsl, \pi'\in\pathsm, \wppprime\leq w\}
	\]
	is a {\em path model\/} for the KK module $\klwm$ in the sense that
	\begin{equation}\label{e:pathkos}
	\character{\klmw}=\sum_{\eta\in\plcpmw}\exp{\eta(1)}
	\end{equation}
\end{theorem}
\bmyproof
From~Theorem~\ref{t:decompkos},  we have:
\begin{equation}\label{e:t:pathkos:1} \character{\klmw}=\sum_{\pi\in\pathsmlw}\character{V_{\lambda+\pi(1)}} \end{equation}
where $\pathsmlw$ is the set of $\lambda$-dominant LS paths of shape $\mu$ with initial direction~$\leq w$.    For $\pi\in\pathsmlw$,  let $\pathsplp$ be the equivalence class in $\pathsl\concat\pathsm$ containing $\pilp$ (under the equivalence relation~$\sim$ defined by the root operators---see \S\ref{sss:sim}),   where $\pil$ denotes the straight line path from the origin to $\lambda$.   Since $\pilp$ lies entirely in the dominant Weyl chamber (this is what it means for $\pi$ to be $\lambda$-dominant),   it follows from the ``Isomorphism Theorem'' in~\cite[Theorem~7.1]{litt:ann} that
\begin{equation}\label{e:t:pathkos:2} \sum_{\sigma\in\pathsplp}\exp{\sigma(1)} \ =\ \sum_{\sigma\in\paths_{\lambda+\pi(1)}}\exp{\sigma(1)} 
\end{equation}
	(where of course $\paths_{\lambda+\pi(1)}$ denotes the set of LS paths of shape $\lambda+\pi(1)$).    By the ``Character formula'' \cite[page~330]{litt:inv}, the right hand side of~\eqref{e:t:pathkos:2} equals $\character{V_{\lambda+\pi(1)}}$, so putting together \eqref{e:t:pathkos:1} and \eqref{e:t:pathkos:2} gives
	\begin{equation}\label{e:t:pathkos:3}
		\character{\klmw}\ = \ \sum_{\pi\in\pathsmlw}\sum_{\sigma\in\pathsplp}\exp{\sigma(1)} 
	\end{equation}
	Thus, for the proof of the theorem,  it suffices to show the following:
	\begin{equation}\label{e:t:pathkos:4}
		\plcpmw\ =\ \sqcup_{\pi\in\pathsmlw}\pathsplp \quad\quad \textup{(disjoint union)}
	\end{equation}
	\newcommand\teta{\tilde{\eta}}
	To prove~\eqref{e:t:pathkos:4},  first let $\pi\in\pathsmlw$.  Let $u$ be an element of the Weyl group such that $u\mu$ is the initial direction of~$\pi$. From our assumption that $uW_\mu\leq wW_\mu$,  it follows that $W_\lambda uW_\mu\leq W_\lambda wW_\mu$ (see Corollary~\ref{c:brcoset}) and $\weyl(\pi_\lambda\concat\pi)\leq w$ (evidently $\weyl(\pi_\lambda\concat\pi)$ is the minimal element in $W_\lambda uW_\mu$).  By Proposition~\ref{p:wind}, it follows that the Weyl group elements associated via~$\weyl$ to elements of~$\pathsplp$ are all the same.   This proves $\plcpmw\supseteq\pathsplp$.

	Now let $\varphi$ be an element in $\plcpmw$. Apply Proposition~\ref{p:dompath} to~$\varphi$ and let $\eta$ be as in the conclusion.   Then $\eta=\pi_\lambda\concat\pi$ for some $\pi\in\pathsmlw$, and the containment $\subseteq$ is proved.  
	\ignore{   
	It follows easily that there is an element---let us say $\teta$---in the equivalence class containing $\eta$ (of the set $\pathsl\concat\pathsm$) that is killed by the root operators~$e_\alpha$ for all simple roots~$\alpha$:   we may arrive at such an element by applying successively the operators $e_\alpha$;  application of $e_\alpha$ shifts the end point up by~$\alpha$, and we cannot keep going ad infinitum since $\pathsl\concat\pathsm$ is a finite set.    From the definition of $e_\alpha$ and item~\eqref{i:hslmconcat} in~\S\ref{ss:stabkosfilt}
it follows that $\teta$ lies entirely in the dominant Weyl chamber.  Writing $\teta$ as a concatenation $\varphi\concat\pi$ with $\varphi\in\pathsl$ and $\pi\in\pathsm$,   it follows that $\varphi$ lies entirely in the dominant Weyl chamber. Thus $\varphi=\pi_\lambda$ (see item~\eqref{i:unique} in~\S\ref{ss:stabkosfilt}) and $\pi$ is $\lambda$-dominant.  
Moreover,  since the association of the Weyl group element~$\weyl$ to paths in~$\pathsl\concat\pathsm$ is invariant under the root operators (Proposition~\ref{p:wind}), we have:
	\[\textup{the initial direction of~$\pi$} = \weyl(\pi_\lambda\concat\pi)=\weyl(\teta)=w(\eta)\leq w\]
	so $\pi$ belongs to $\pathsmlw$.   Evidently,  $\eta$ belongs to $\pathsplp$,  which proves the containment $\subseteq$ in~\eqref{e:t:pathkos:4}.
}

	That the union on the right hand side of~\eqref{e:t:pathkos:4} is disjoint follows from the uniqueness of~$\eta$ in Proposition~\ref{p:dompath}  (which in turn rests on~\cite[Corollary~1\,(b) of~\S7]{litt:ann}):  $\pi_\lambda\concat\pi$ is the unique path in $\pathsplp$ on which $e_\alpha$ vanishes for all simple roots~$\alpha$.
	\emyproof

\mysection{PRV components and generalised PRV components in KK modules}\mylabel{s:randgkprv}
\noindent
We show how the decomposition rule (Theorem~\ref{t:decompkos}) leads easily to results about the existence of PRV components (Theorem~\ref{t:rkprvforkos})  and generalised PRV components (Theorem~\ref{t:gkprv}) in KK modules.    The arguments are well known: see e.g. those by Joseph in~\cite[\S2.7]{joseph:hc}.   In fact,   Theorem~\ref{t:rkprvforkos} for the finite case follows from items~(i) and~(iii) of the theorem in~\cite[\S2.7]{joseph:hc}.

Theorem~\ref{t:rkprvforkos} is at once a generalisation of two results:  the so called refined PRV and KPRV theorems:
\begin{itemize}
	\item
Its special case when $\lieg$ is of finite type and $w=w_0$ (the longest element of the Weyl group) is due to Kumar~\cite[Theorem~1.2]{kumar:refineprv}, who refers to his result as  ``a refinement of the PRV conjecture'' and says that it was conjectured by D.-N.~Verma.
\item
		The special case when $\sigma=w$ was proved by Kumar~\cite[page~117]{kumar:prv} and independently Mathieu~\cite[Corollaire~3]{mathieu}.  Kumar calls it ``the strengthened PRV conjecture (due to Kostant)''. We have called it ``KPRV'' following Khare~\cite{khare}.
\end{itemize}
Theorem~\ref{t:gkprv} is a KK version of Montagard's result~\cite[Theorem~3.1]{montagard} about generalised PRV components.

\ignore{  
The association $\wbar\leftrightarrow w\mu$  (for $w\in W$) gives a bijection of the coset space~$W/W_\mu$ with the set $W\mu$ of $W$-conjugates of~$\mu$.  We identify the sets $W/W_\mu$ and $W\mu$ via this bijection.   The double coset space~$\wlwwm$ may then be identified with the set of $\wl$-orbits of $W\mu$.
\mysubsection{Geometric interpretation of $\wlwwm$}\mylabel{ss:wlwwm}
\bprop\mylabel{p:wlwwm}  Every $\wl$-orbit of the set $W\mu$ of $W$-conjugates of~$\mu$ contains a unique element $w\mu$ such that $\lambda+tw\mu$ is dominant for some real number $t>0$.
\eprop
\bmyproof  
Each such orbit contains a unique $w\mu$ that is $\wl$-dominant.  The $\wl$-dominance means precisely that $\langle w\mu,\alpha^\vee\rangle\geq0$ for every simple root $\alpha$ in~$\wl$.   It is easily verified that 
$w\mu$ has the desired property.    Conversely,  if $w\mu$  is not $\wl$-dominant,  then $\langle w\mu,\alpha^\vee\rangle <0$ for some simple root~$\alpha$ in~$\wl$,  and so $\langle\lambda+tw\mu,\alpha^\vee\rangle=t\langle w\mu,\alpha^\vee\rangle<0$ for $t>0$.
\emyproof
The double coset space $\wlwwm$ may thus be identified with the set of those Weyl conjugates $w\mu$ of~$\mu$ such that $\lambda+tw\mu$ is dominant for some positive $t$.
\bprop\mylabel{p:wlwwm2} Given a double coset in~$\wlwwm$,  let $w$ be the unique minimal element in it with respect to the Bruhat order (as guaranteed by Corollary~\ref{c:cpparpar}).  Then $w\mu$ is such that $\lambda+tw\mu$ is dominant for all small positive~$t$.
\eprop
\bmyproof  From the proof of Proposition~\ref{p:wlwwm},  it is enough to show that $w\mu$ is $\wl$-dominant.    Suppose that this is not so.   Then there exists simple root $\alpha$ with $s_\alpha$ in~$\wl$ such that $\langle w\mu,\alpha^\vee\rangle<0$.   We then have $s_\alpha w<w$, which contradicts the hypothesis that $w$ is the minimal element in its double coset.
\emyproof



}
\mysubsection{The map $\mapphi$}\mylabel{ss:eta}
Let $\lieg$ be a symmetrizable Kac-Moody algebra,  and fix dominant integral weights $\lambda$ and $\mu$.  Let $\wl$ and $\wm$ denote respectively the stabilizers in the Weyl group~$W$ of $\lambda$ and $\mu$.

Consider the map from the Weyl group $W$ to the set $\Domint$ of dominant integral weights given by $\sigma\mapsto \barlsm$,  where $\barlsm$ denotes the dominant Weyl conjugate of the weight $\lambda+\sigma\mu$.  This map factors through the natural quotient map from $W$ to~$\wlwwm$.  We denote by~$\mapphi$ the map $\wlwwm\to\Domint$ given by $\wl\sigma\wm\mapsto\barlsm$.
\mysubsection{PRV components in KK modules}\mylabel{ss:kprv}
\noindent
The restrictive hypothesis on~$\lieg$ in the following theorem (as also in Theorem~\ref{t:gkprv}), namely that it be either of finite type or symmetric,  is inherited from the decomposition rule (Theorem~\ref{t:decompkos}) and is possibly not required.
\begin{theorem}\mylabel{t:rkprvforkos}
	Let $\lieg$ be a symmetrizable Kac Moody algebra that is either of finite type or symmetric (as in Theorem~\ref{t:decompkos}).
	Let $\lamdba$, $\mu$ be dominant integral weights and $w$, $\sigma$ be elements of the Weyl group.   Let $\nu$ be the dominant Weyl conjugate of the weight $\lambda+\sigma\mu$.   Then the irreducible $\lieg$-module~$\Vnu$ occurs in the decomposition into irreducible $\lieg$-modules of the KK module $\klwm$ at least as many times as
	there are elements $\tau\in\wlwwm$ such that $\tau\leq \wl w\wm$ and  $\mapphi(\tau)=\nu$.
\end{theorem}
\bmyproof  
We describe a map $\tmapphi$ from  $\wlwwm$ to the set $\pathsml$ of $\lambda$-dominant LS paths of shape~$\mu$.  Given $\tau\in\wlwwm$, let $v$ be the unique minimal element in~$\tau$.  Consider the path $\varphi=\pi_\lambda\concat\pi_{v\mu}$ in $\paths\concat\pathsm$,  where $\pi_\lambda$ and $\pi_{v\mu}$ are the straight line paths from the origin to $\lambda$ and $v\mu$ respectively.  Note that $\weyl(\varphi)=v$.  Apply Proposition~\ref{p:dompath} to~$\varphi$ and let $\eta$ be as in its conclusion.   
Then  $\eta=\pil\concat\pi$ for some $\pi\in\pathsml$ and $\weyl(\eta)=v$.   We define $\tmapphi(\tau):=\pi$.  Since $\pi$ determines $\eta$ from which we can recover $v$ and in turn~$\tau$,  it follows that $\tmapphi$ is injective.

It follows from Proposition~\ref{p:extremal} that $\varphi$ as above is extremal, which means that $\eta(1)=\overline{\varphi(1)}=\overline{\lambda+v\mu}=\mapphi(\tau)$.    Thus $\tmapphi$ is a ``lift'' to $\pathsml$ of $\mapphi$,  meaning that $\mapphi(\tau)$ is the end point of $\tmapphi(\tau)$ shifted by $\lambda$ for any $\tau$ in~$\wlwwm$.  Combining this fact and the injectivity of~$\tmapphi$ with the decomposition rule (Theorem~\ref{t:decompkos}), we immediately obtain the theorem.
%
\emyproof
\mysubsection{KPRV recovered}\mylabel{ss:kprv}\noindent
It follows immediately from the theorem that $\Vnu$ occurs at least once in~$\klsm$.    To show that it occurs at most once,  we repeat the following elementary argument from~\cite[\S2.7]{kumar:prv}.  
Indeed in any $\lieg$-homomorphism from $\klsm$ module to $\Vnu$,   the vector $\vlambda\tensor v_{\sigma\mu}$ has to map to an element of weight $\lambda+\sigma\mu$.     But the dimension of the $\lambda+\sigma\mu$-weight space in $\Vnu$ is clearly one,  since $\lambda+\sigma\mu$ is a Weyl conjugate of~$\nu$.   Thus the space of $\lieg$-homomorphisms from $\klsm$ to $\Vnu$ is one dimensional.   We thus have:
\bcor\mylabel{c:kprv}
Let $\lieg$, $\lambda$, $\mu$, $\sigma$, and $\nu$ be as in Theorem~\ref{t:rkprvforkos}.    Then the irreducible $\lieg$-module $\Vnu$ occurs exactly once in the decomposition into irreducible $\lieg$-modules of the KK module $\klsm$.
\ecor
\mysubsection{Generalised PRV components in KK modules}\noindent
Importing to our context a result of Montagard~\cite[Theorem~3.1]{montagard}, we prove the following:
\begin{theorem}\mylabel{t:gkprv} Let $\lieg$, $\lambda$, and $\mu$ be as in~Theorem~\ref{t:rkprvforkos}.
	Let $v$, $u$ be elements in the Weyl group and $\beta$ a positive root such that either $v^{-1}\beta$ or $u^{-1}\beta$ is a simple root.   Let $k$ be an integer such that $0\leq k\leq \min\{\langle v\lambda,\beta^\vee\rangle, \langle u\mu,\beta^\vee\rangle\}$ and the integral weight $\nu=v\lambda+u\mu-k\beta$ is dominant.      Then the irreducible $\lieg$-module~$V_\nu$ occurs in the decomposition of the KK module $\klwm$ into irreducibles where $w=v^{-1}s_\beta u$.
\end{theorem}
\bmyproof   First suppose that $v^{-1}\beta$ is simple.    Since $k\leq \langle v^{-1}u\mu,v^{-1}\beta^\vee\rangle$,  it follows that $f^k_{v^{-1}\beta} \pi_{v^{-1}u\mu}$ does not vanish.  Consider the path $\varphi:=\pi_\lamdba\concat f^k_{v^{-1}\beta}\pi_{v^{-1}u\mu}$ in~$\pathsl\concat\pathsm$. 
As is easily verified, the dominant Weyl conjugate of $\varphi(1)$ is $\nu$ and $\weyl(\varphi)$ is either $\brmin{W_\lambda v^{-1}s_\beta u W_\mu}$ or $\brmin{W_\lambda v^{-1}uW_\mu}$ depending upon whether $k>0$ or $k=0$.   An easy verification (using the fact that $\langle u\mu,\beta^\vee\rangle\geq0$ under our assumptions) shows that $w=v^{-1}s_\beta u\geq v^{-1}u$.  Thus $w\geq\weyl(\varphi)$ in either case, and $\varphi\in\plcpmw$.

Apply Proposition~\ref{p:dompath} to the path~$\varphi$ and let $\eta$ be as in its conclusion.  Then $\eta$ is of the form $\pi_\lambda\concat\pi$ with $\pi\in\pathsmlw$.  By the decomposition rule (Theorem~\ref{t:decompkos}),   $V_{\eta(1)}$ occurs in $\klmw$.   But Montagard~\cite[Proof of Theorem~3.1]{montagard} shows that $\varphi$ is extremal, which means that $\nu=\overline{\varphi(1)}=\eta(1)$, and the proof is done in this case.

\begin{figure}[ht!]\label{fig3}
 \centering
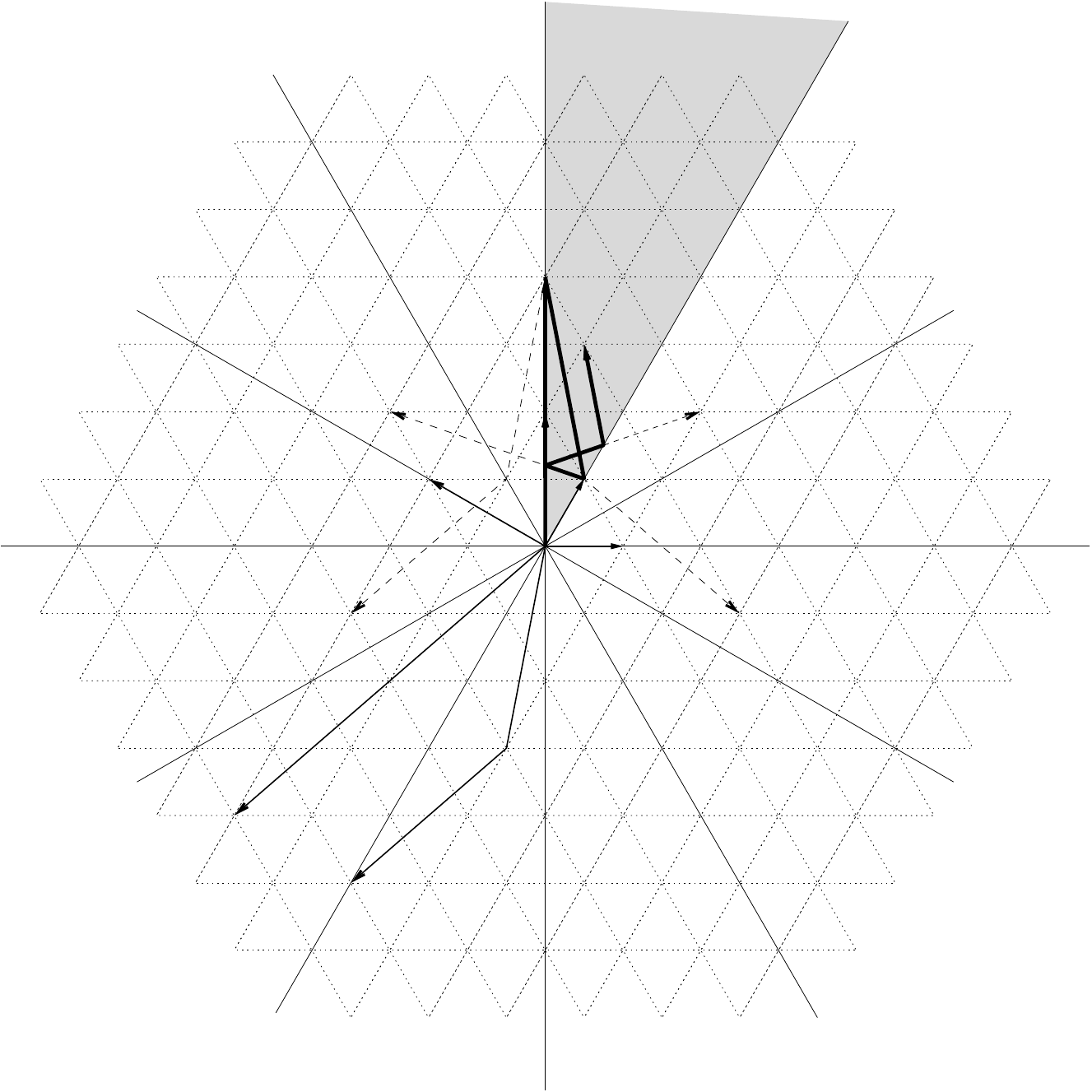 
\caption{$V_\nu \in K(\lambda, w, \mu)$}
	\label{f:g2}
\end{figure}

Now suppose that $u^{-1}\beta$ is simple.   Then, applying the result in the previous case,  we conclude that $V_\nu$ occurs in the KK submodule~$\kmwinvl$ of $V_\mu\tensor V_\lambda$.   But under the $\lieg$-isomorphism $a\tensor b\leftrightarrow b\tensor a$ of $V_\mu\tensor V_\lambda$ with $V_\lambda\tensor V_\mu$,  the submodules $\kmwinvl$ and $\klwm$ map isomorphically to each other.
\emyproof
\bremark
In the set up of the theorem,  suppose we assume that $\beta$ is simple rather than that either $v^{-1}\beta$ or $u^{-1}\beta$ is simple.   In that case too~\cite[Theorem~3.1]{montagard} says that $V_\nu$ occurs in~$V_\lambda\tensor V_\mu$.   We have not handled that case.  
\eremark
\bexample\mylabel{x:g2} We illustrate the result of Theorem~\ref{t:gkprv} and also the idea behind its proof by means of an example borrowed from Montagard~\cite{montagard}.   Let the root system of~$\lieg$ be $G_2$.    Let $e_1$ and $e_2$ be the standard basis vectors of $\mathbb{R}^2$ with its standard inner product.  We may take $\alpha_1:=e_1$ and $\alpha_2:=-\frac{3}{2}e_1+\frac{\sqrt{3}}{2}e_2$ to be the simple roots.  The set of all positive roots is:
\[ \{\alpha_2,\alpha_1, \alpha_2+\alpha_1, \alpha_2+2\alpha_1, \alpha_2+3\alpha_1, 2\alpha_2+3\alpha_1\} \]
The dominant integral weights $\lambda$, $\mu$, $\nu$,   the Weyl group elements $u$, $v$, $w$,   the root~$\beta$, and the integer $k$ are all as shown in Figure~\ref{f:g2}.   The path $\eta$ ending at $\nu$ appears in bold.
\eexample

\mysection{Tableau decomposition rule for Kostant-Kumar (KK) modules}\mylabel{s:tabkk}
\noindent
Fix an integer $d\geq2$. Let~$\lieg=\sld$, the simple Lie algebra of traceless complex $d\times d$ matrices.
There is, in this special case,  the classical Littlewood-Richardson (LR for short) rule (see e.g. \cite{mcld, fulton:yt}) that gives, in terms of tableaux, the decomposition into irreducibles of the tensor product of two finite dimensional irreducible representations of~$\lieg$.     The multiplicities of the irreducibles in this rule are called ``LR coefficients'' and they count certain ``LR tableaux''.
Our purpose in this section is to deduce, from the general decomposition rule (Theorem~\ref{t:decompkos}), 
a version of this classical rule, which we call the ``refined LR rule'', for decomposing as a direct sum of irreducibles any KK submodule of the tensor product: see~\S\ref{ss:tabkk} for the statement.     We call the multiplicities of the irreducibles in this refined rule the ``refined LR coefficients''.     

The refined LR coefficients also count certain LR tableaux.   The identification of the set of LR tableaux to be counted is based upon the association of a permutation to each LR tableau (\S\ref{ss:st-to-w}).     There is a very natural association of a semi-standard Young tableau (SSYT) to an LR tableau (\S\ref{sss:st-to-ssyt}) and the permutation is just 
the initial element of the minimal standard lift when the SSYT is interpreted as a standard concatenation of LS paths (\S\ref{ss:slla}).

In the light of this last mentioned fact, it is noteworthy that the procedure we give for determining the permutation (\S\ref{sss:ssyt-to-w}) from the SSYT is not a repeated application of Deodhar's lemma (Proposition~\ref{p:deodhar}):   it seems to be more efficient than that.  
Lascoux and Sch\"utzenberger~\cite{ls:keys} associate to each SSYT a ``right key'' (which by definition is another SSYT) from which the permutation can be read off.
Willis~\cite{willis} gives an alternative method--``the scanning method''---for finding the right key of an SSYT.   Our procedure is different from those in~\cite{ls:keys,willis}.

\mysubsection{Preliminaries}\mylabel{ss:tabprelims}
The choices involved (Cartan subalgebra $\lieh$, Borel subalgebra~$\lieb$, etc.)\ 
are fixed as usual:  the subalgebra of diagonal (respectively, upper triangular) traceless complex $d\times d$ matrices is taken to be $\lieh$ (respectively, $\lieb$). 
We denote by $\epsilon_j$ the linear functional on $\lieh$ that maps a matrix to its entry in position $(j,j)$.

Recall that a {\em partition\/} is a weakly decreasing sequence $\lambda_1\geq\lambda_2\geq\ldots$ (sometimes also written $\lambda_1+\lambda_2+\cdots$) of non-negative integers that is eventually zero.  
The non-zero elements of the sequence are called the {\em parts}. 
We tacitly identify partitions with their (Young) shapes. 
To a partition $\lambda: \lambda_1\geq\ldots\geq\lambda_d\geq0\geq\ldots$ with at most $d$ parts, we attach the dominant integral weight $\lambda_1\epsilon_1+\cdots+\lambda_d\epsilon_d$.     A second such partition $\lambda'_1\geq\ldots\geq\lambda'_d\geq0\geq\ldots$ corresponds to the same weight as $\lambda$ if and only if $\lambda_1-\lambda'_1=\cdots=\lambda_d-\lambda'_d$ (since $\epsilon_1+\ldots+\epsilon_d=0$ is evidently the only linear dependence relation up to scaling on $\epsilon_1$, \ldots, $\epsilon_d$).   Thus partitions with less than $d$ parts are in one-to-one correspondence with dominant integral weights.   We will abuse notation and use the same symbol for both a partition with less than $d$ parts and the corresponding dominant integral weight.

Let $[j]:=\{1,\ldots,j\}$ for any integer $j\geq1$.    The Weyl group is identified with the group of permutations of the set~$[d]$.  The one line notation for a permutation~$w$ of~$[d]$ is $w_1\ldots w_d$,  where $w_j:=w(j)$ (for $1\leq j\leq d$).
\mysubsection{Semi-standard Skew tableaux (SSST for short)}\mylabel{ss:skewtab}
Let $\nu$ and $\lambda$ be two partitions with the shape of $\nu$ containing the shape of~$\lambda$.  A {\em (semi-standard) skew tableau\/}, {\em SSST\/} for short, of shape $\nu/\lambda$ is a filling up by positive integers of those boxes that are in the shape of~$\nu$ but not in the shape of~$\lambda$ such that
the entries in each row are weakly increasing rightward and those in each column are strictly increasing downward. Here are two examples with $\nu= 7+5+4+3+1$ and $\lambda=4+4+1+1$:
\begin{equation}\label{e:skewtab} 
\begin{tabular}{|c|c|c|c|c|c|c|c|c|c|c|}
\hline
	$\cdot$ &$\cdot$ &$\cdot$ & $\cdot$ & 2 & 2  &3 \\
	\cline{1-7}
$\cdot$ & $\cdot$ &$\cdot$ &$\cdot$ &7 \\
	\cline{1-5}
$\cdot$ & 3 & 3 & 5 \\
	\cline{1-4}
$\cdot$ & 4 & 6\\
	\cline{1-3}
1\\
	\cline{1-1}
\end{tabular}
\ \ \ \ \ \ \ \ \hskip2cm
\begin{tabular}{|c|c|c|c|c|c|c|c|c|c|c|}
\hline
	$\cdot$ &$\cdot$ &$\cdot$ &$\cdot$ & 1 & 1  &1 \\
	\cline{1-7}
$\cdot$ & $\cdot$ &$\cdot$ &$\cdot$ &2 \\
\cline{1-5}
$\cdot$ & 1 & 2 & 3 \\
	\cline{1-4}
$\cdot$ & 3 & 4\\
	\cline{1-3}
1\\
	\cline{1-1}
\end{tabular}
\end{equation}
\ignore{
\mysubsubsection{Bar separated reverse reading words}\mylabel{sss:bsrw}\mylabel{ss:bsrw}
We associate to $T$ its {\em bar separated reverse reading word\/} as follows.  Imagine identical vertical bars,  one less than the number of parts in~$\nu$,   placed on a horizontal line,   like so:
\[ \vert\quad\vert\quad\quad\cdots\quad\quad\vert \quad\vert \]
The spaces between two consecutive vertical bars,  as also the space to the left of the leftmost vertical bar,  and the space to the right of the rightmost vertical bar are called {\em blocks\/}.
We thus have $\nu'_1$ blocks,  where $\nu'_1$ is the number of parts in~$\nu$.   The blocks are numbered starting from $1$,  from left to right.   
For every $i$, $1\leq i\leq \nu'_1$,  we populate the $i^\textup{th}$ block with the entries in the $i^\textup{th}$ row of the skew tableau~$T$ but in reverse order:   the entries in every block will thus be weakly decreasing from left to right.

The bar separated reverse reading words of the two skew tableaux in the display above are respectively:
\begin{equation}\label{e:bsrw}
	3\ 2\ 2\ |\ 7\ |\ 5\ 3\ 3\ |\ 6\ 4\ |\ 1\ 
\ \ \ \ \ \ \ \ \hskip2cm
	1\ 1\ 1\ |\ 2\ |\ 3\ 2\ 1\ |\ 4\ 3\ |\ 1\ 
\end{equation}

	Given the bar separated reverse reading word of $T$ and either $\nu$ or $\lambda$,  we can recover $T$.
\mysubsubsection{Reverse reading words and ballot sequences}\mylabel{ss:rwballotseq} \mylabel{sss:rwballotseq} 
	The {\em reverse reading word\/} of~$T$ is just the bar separated reverse reading word of~$T$ with the vertical bars removed.  It is said to be a {\em ballot sequence\/} if for any integer $j\geq1$  the number of times $j$ occurs up to any point in the word (while scanning it from left to right) is at least the number of times $j+1$ occurs up to that point.   In the examples above,  the reverse reading word on the left is not a ballot sequence but the one on the right is.

If the reverse reading word (of a skew tableau~$T$) is a ballot sequence,  then its {\em type\/} is the partition $\mu=\mu_1+\mu_2+\cdots$,  where $\mu_j$ is the number of times $j$ occurs in the reverse reading word.  For example, the type of the word on the right in~\eqref{e:bsrw} is $\mu=5+2+2+1$.
}
\mysubsubsection{Reverse reading words and ballot sequences}\mylabel{ss:rwballotseq} \mylabel{sss:rwballotseq} 
Let $T$ be a SSST of shape $\nu/\lambda$.
Its {\em reverse reading word\/}, denoted $\wordr(T)$, is defined as follows:  read the entries of $T$ from right to left in every row,  scanning the rows from top to bottom.  
For the two SSSTs in the display above, the reverse reading words respectively are:
\begin{equation}\label{e:bsrw}
	3\ 2\ 2\ 7\ 5\ 3\ 3\ 6\ 4\ 1\ 
\ \ \ \ \ \ \ \ \hskip2cm
	1\ 1\ 1\ 2\ 3\ 2\ 1\ 4\ 3\ 1\ 
\end{equation}
	The word~$\wordr(T)$ (or, more generally, any word in the positive integers) is said to be a {\em ballot sequence\/} if for any integer $j\geq1$  the number of times $j$ occurs up to any point in the word (while scanning it from left to right) is at least the number of times $j+1$ occurs up to that point.   In~\eqref{e:bsrw},  the word on the left is not a ballot sequence but the one on the right~is.
	\mysubsubsection{Type and weight of a word and of a SSST}\mylabel{sss:type}
	The {\em type\/} of any word $\word$ in the positive integers is the sequence $\mu$: $\mu_1$, $\mu_2$, \ldots,  where $\mu_j$ denotes the number of occurrences of $j$ in~$\word$.    The type of the word on the left in~\eqref{e:bsrw} is~$1$, $2$, $3$, $1$, $1$, $1$, $1$, $0$, $0$, \ldots.   Evidently, permuting the letters of a word does not change its type.   
If $\word$ is a word in~$[d]$, then we may further associate to it the integral weight $\mu_1\epsilon_1+\cdots+\mu_d\epsilon_d$ of~$\lieg=\sld$.  This is called the {\em weight\/} of the word and denoted~$\weight(\word)$. 

	The {\em type\/} and {\em weight\/} of a SSST~$T$ are defined respectively to be the type and weight of its reverse reading word~$\wordr(T)$.\footnote{Later on we will introduce the ``column word'' $\wordc(T)$ of~$T$,  which being a permutation of~$\wordr(T)$ shares its type and weight.}

If the word~$\word$ is a ballot sequence,  then its {\em type\/} is a partition: $\mu_1\geq\mu_2\geq\ldots$,  and in this case we use the notation for partitions to denote types.  For example, the type of the word on the right in~\eqref{e:bsrw} is $\mu=5+2+2+1$.    The weight of such a word in~$[d]$ is dominant.
\mysubsection{Littlewood-Richardson (LR for short) tableaux and coefficients}\mylabel{ss:lrcoeffs}
An {\em LR tableau\/} (LR is short for {\em Littlewood-Richardson\/}) is a SSST~$T$ whose reverse reading word~$\wordr(T)$ is a ballot sequence.
Let $\lambda$ and $\mu$ be partitions.
Let $\skewtabml$ denote the set of LR tableau of shape $\nu/\lambda$ and type~$\mu$---here $\nu$ is allowed to vary.
If $T$ in~$\skewtabml$ has shape $\nu/\lambda$,  we write $\nu(T)$ for~$\nu$.
As is well-known,  $\skewtabml$ has representation theoretic and geometric significance.
For example (see e.g.~\cite{mcld, fulton:yt})
	$s_\lambda s_\mu=\sum_{T\in\skewtabml}  s_{\nu(T)}$,
	where $s_\tau$ denotes the Schur function associated to a partition~ $\tau$.

	For a fixed partition $\nu$,  the number of $T$ in~$\skewtabml$ with $\nu(T)=\nu$ is usually denoted $c_{\lambda\mu}^\nu$.   The numbers $c_{\lambda\mu}^\nu$ are called {\em LR coefficients\/}.    In terms of these,  we may write the the above rule for multiplication of Schur functions as 
$s_\lambda s_\mu=\sum_\nu c_{\lambda\mu}^\nu s_\nu$.
\mysubsubsection{Bruhat order on permutations}\mylabel{sss:brperm}   Any permutation $u$ of $[j]$ (for some integer $j\geq1$) can naturally be considered as a permutation of $[k]$, for any integer $k\geq j$.    Given two permutations $u$ and $u'$ (of $[j]$ and $[j']$ respectively),   we write $u\leq u'$ if that is so in the Bruhat order on permutations of $[k]$ for some $k\geq$ both $j$ and $j'$.   If $u\leq u'$ for one such $k$,  then it is so for all such $k$.
\mysubsubsection{Refined Littlewood-Richardson coefficients: their definition}\mylabel{ss:rlr}\mylabel{sss:rlr}
\noindent    
In~\S\ref{ss:st-to-w} below, we specify a procedure that assigns a permutation~$\procw$ to a given SSST~$T$.%
\footnote{It is easy to associate to~$T$ a SSYT~$\ssyt$ of shape~$\mu$---see~\S\ref{sss:st-to-ssyt}. Interpreting $\ssyt$ as a standard concatenation of LS paths in the sense of Proposition~\ref{p:stdjust} in the appendix,   the associated permutation~$\procw$ is just the initial element of the minimal standard lift of~$\ssyt$, as will be proved in \S\ref{s:proc-proof}.  Observe that,  if as in~\S\ref{ss:tabkk} the number of parts in $\nu$ is at most $d$,  then the entries in $\ssyt$ and the number of parts in~$\mu$ are also bounded above by $d$,  so that the interpretation of~$\ssyt$ as a concatenation of LS paths associated to $\lieg=\sld$ is possible,  and $\procw$ is a permutation of~$[d]$.}   
Fix a permutation~$w$ and let $\skewtabmlw$ denote the subset of~$\skewtabml$ consisting of those elements for which the associated permutation $\procw$ satisfies $\procw\leq w$ (in the Bruhat order as defined in~\S\ref{sss:brperm} above).
The result~\eqref{e:kktab} below ascribes representation theoretic meaning to~$\skewtabmlw$.

	For a fixed partition $\nu$,  we denote by $c_{\lambda\mu}^\nu(w)$ the number of $T$ in~$\skewtabmlw$ with $\nu(T)=\nu$.   We call the numbers $c_{\lambda\mu}^\nu(w)$ {\em refined LR coefficients\/}.
\mysubsection{Tableau decomposition rule for KK modules}\mylabel{ss:tabkk}
\noindent
Suppose that $\lambda$, $\mu$ are partitions with less than $d$ parts (or, equivalently, dominant integral weights for $\lieg=\sld$)
and that $w$ is a permutation of~$[d]$ (or, equivalently, an element of the Weyl group). 
Then the decomposition of the Kostant-Kumar module $\koslwm$ (defined in~\S\ref{ss:kosfilt}) as a direct sum of irreducible $\lieg$-modules is given by:
\begin{equation}\label{e:kktab}
	\boxed{	\koslwm = \bigoplus_{T\in\skewtabml(w)}V_{\nu(T)}}
\end{equation}
where $V_{\nu(T)}$ is interpreted to be zero in case $\nu(T)$ has more than $d$ parts.  (Recall from~\S\ref{ss:tabprelims} that to any partition with at most $d$ parts there is associated a dominant integral weight of~$\lieg$.)

Here is an alternative way to express the above decomposition rule:
\begin{equation}\label{e:kktabalt}
	\boxed{	\koslwm = \bigoplus_{\nubar}V_{\nubar}^{\oplus c_{\lambda\mu}^\nu(w)}}
\end{equation}
where the sum runs over all partitions $\nubar$ with less than $d$ parts, and $\nu$ depending on~$\nubar$ denotes the unique partition with at most $d$ parts (if it exists) such that \[\textup{ $\nubar_j=\nu_j-\nu_d$ for $1\leq j<d$ \quad and \quad $\nu_1+\cdots+\nu_d=(\lambda_1+\cdots+\lambda_{d-1})+(\mu_1+\cdots+\mu_{d-1})$}\]

The proof of~\eqref{e:kktab} will be given below in~\S\ref{ss:pfkktab}.
\mysubsubsection{The statement for polynomial representations of~$GL_d(\mathbb{C})$}\mylabel{sss:poly}   For convenience of reference,  we now state, without proof,  a version of the decomposition rule~\eqref{e:kktab} for polynomial representations of the general linear group~$GL_d(\mathbb{C})$.
Suppose that $\lambda$, $\mu$ are partitions with at most $d$ parts and $\Vlambda$, $\Vmu$ the corresponding irreducible polynomial representations.
Let $w$ be a permutation of~$[d]$.  Then the decomposition of the Kostant-Kumar module $\koslwm$ (defined similarly as in~\S\ref{ss:kosfilt}) as a direct sum of irreducible polynomial representations is given by:
\begin{equation}\label{e:kktabpoly}
	\boxed{	\koslwm = \bigoplus_{\nu}V_{\nu}^{\oplus c_{\lambda\mu}^\nu(w)}}
\end{equation}
where the sum runs over all partitions $\nu$ with at most $d$ parts.
\mysubsubsection{An example}\mylabel{sss:decompex}   Here is a simple example illustrating the rules~\eqref{e:kktab} and~\eqref{e:kktabalt}.
Let $d=3$, $\lambda=2+1$, and $\mu=3+1$.   As the reader can readily verify,   there are $7$ elements~$T$ in~$\skewtabml$ with $\nu(T)$ having at most $3$ parts.  These are listed below along with the permutations of~$[3]$ attached to them (as in~\S\ref{ss:st-to-w}):
\begin{gather*}
	\begin{array}{|c|c|c|c|c|}
		\hline
		\cdot &\cdot  & 1 & 1 & 1 \\
		\hline
		\cdot & 2\\
		\cline{1-2}
	\end{array}
	\leftrightarrow 123
	\quad\quad
	\begin{array}{|c|c|c|c|}
		\hline
		\cdot &\cdot  & 1 & 1  \\
		\hline
		\cdot & 1 & 2\\
		\cline{1-3}
	\end{array}
	\leftrightarrow 213 
	\quad\quad
	\begin{array}{|c|c|c|c|c|}
		\hline
		\cdot &\cdot  & 1 & 1 & 1 \\
		\hline
		\cdot \\
		\cline{1-1}
		2 \\
		\cline{1-1}
	\end{array}
	\leftrightarrow 132
\end{gather*}
\begin{gather*}
	\begin{array}{|c|c|c|c|}
		\hline
		\cdot &\cdot  & 1 & 1  \\
		\hline
		 \cdot & 1 \\
		\cline{1-2}
		2 \\
		\cline{1-1}
	\end{array}
	\leftrightarrow 231 
	\quad\quad
	\begin{array}{|c|c|c|c|}
		\hline
		\cdot &\cdot  & 1 & 1  \\
		\hline
		 \cdot & 2 \\
		\cline{1-2}
		1 \\
		\cline{1-1}
	\end{array}
	\leftrightarrow 312 
	\quad
	\quad\quad
	\begin{array}{|c|c|c|}
		\hline
		\cdot &\cdot  & 1   \\
		\hline
		 \cdot &  1 & 2 \\
		 \hline
		1 \\
		\cline{1-1}
	\end{array}
	\leftrightarrow 312 
	\quad\quad
	\begin{array}{|c|c|c|}
		\hline
		\cdot &\cdot  & 1   \\
		\hline
		 \cdot &  1  \\
		 \cline{1-2}
		1  & 2 \\
		\cline{1-2}
	\end{array}
	\leftrightarrow 321 
\end{gather*}
And so we have:
\begin{gather*}
	K(\lambda, 123, \mu) = V_{5+2}\quad\quad
	K(\lambda, 213, \mu) = V_{4+3}\oplus V_{5+2}\quad\quad
	K(\lambda, 132, \mu) = V_{4}\oplus V_{5+2}\\
	K(\lambda, 231, \mu) =  V_{3+1}\oplus V_{4}\oplus V_{4+3}\oplus V_{5+2}\quad
	K(\lambda, 312, \mu) = V_{2+2}\oplus V_{3+1}\oplus V_{4}\oplus V_{4+3}\oplus V_{5+2}\\
	K(\lambda, 321, \mu) = V_{1}\oplus V_{2+2}\oplus V_{3+1}\oplus V_{3+1}\oplus V_{4}\oplus V_{4+3}\oplus V_{5+2}
\end{gather*}
\mysubsection{SSYT and permutations attached to them}\mylabel{ss:ssyt-and-perm}
Let $\mu$ be a partition.
A {\em semi-standard Young tableau\/}, {\em SSYT\/} for short, {\em of shape~$\mu$\/} is just a (semi-standard) skew tableau of shape $\mu/\textup{empty}$ in the sense of~\S\ref{ss:skewtab}.
Here is an example of a SSYT of shape $\mu=4+2+1$:
\begin{equation}\label{e:ssytt'}\begin{array}{|c|c|c|c|}\hline 1 & 3 & 6 & 8 \\
	\hline
	2 & 4 \\
	\cline{1-2}
	7 \\
	\cline{1-1}
\end{array} \end{equation}
\mysubsubsection{Associating a permutation to a given SSYT}\mylabel{sss:ssyt-to-w}  Let $\ssyt$ be a SSYT of shape~$\mu$ and let $m$ be the largest entry of~$\ssyt$.    We associate to~$\ssyt$ a permutation~$\procw$ of~$[m]$, as follows.  
Let $\mpart$ be the number of parts in~$\mu$.  Observe that $m\geq \mpart$ since the entries in every column of~$\ssyt$ are strictly increasing downwards.

Let $\procw_1\procw_2\ldots \procw_m$ be the one-line notation for $\procw$. 
We will describe below an inductive procedure to produce the sequence $\procw_1$, \ldots, $\procw_{\mpart}$.    As for $\procw_{\mpart+1}$, \ldots, $\procw_m$,   we take these to be just the elements of $[m]\setminus\{\procw_1,\ldots,\procw_\mpart\}$ arranged in increasing order.

It is easy to produce~$\procw_1$: it is just the largest (right most) entry in the first row of~$\ssyt$.   
Suppose that $\procw_1$, $\procw_2$, \ldots, $\procw_{\tint-1}$ have been produced (with $1<\tint\leq \mpart$).   
%
%
We now describe a procedure to determine~$\procw_\tint$.

Let $\boxb$ be a box in $\ssyt$.  
Suppose that a box $\boxbp$ in~$\ssyt$ is weakly to the Northeast of~$\boxb$ and has an entry that is less than that of $\boxb$.   Then we write $\boxb\cgeq\boxbp$.
For example,   in the SSYT of~\eqref{e:ssytt'},  if $\boxb$ is the one with entry $7$,  then $\boxbp$ could be any of those containing $1$, $2$, $3$, $4$, or $6$;   if $\boxb$ is the one with entry $4$,  then $\boxbp$ could only be the one containing $3$.

The {\em $\boxb$-depth\/} of such a box $\boxbp$ is defined to be the largest $\delta$ such that there is a chain $\boxb\cgeq\boxb_1\cgeq\boxb_2\cgeq\ldots\cgeq\boxb_\delta=\boxbp$.  The $\boxb$-depth of $\boxb$ itself is defined to be $0$.

Let $\boxbtint$ denote the right most box in row $\tint$.  We write {\em $\tint$-depth\/} for $\boxbtint$-depth.    For $1\leq j\leq \tint$,  we let $y_j$ be maximal possible entry in a box whose $\tint$-depth is~$\tint-j$.  (The box in row $j$ in the same column as $\boxbtint$ has $\tint$-depth $\tint-j$,  so $y_j$ exists.)    By definition, $y_{\tint}$ is the entry in the box $\boxbtint$.    As is easily seen, $y_1<\ldots<y_{\tint}$.   We call this the {\em $\tint$-depth sequence of~$\ssyt$\/}.

Let $a_1<\ldots<a_{\tint-1}$ be the elements $\procw_1$, \ldots, $\procw_{\tint-1}$ arranged in increasing order.
Let $k$, $1\leq k\leq p$, be the largest such that $a_{k-1}<y_k$ ($a_0=-\infty$ by convention).   Take~$\procw_\tint$ to be~$y_k$.
\ignore{   
The recursive procedure to produce~$w_{r+1}$ involves a subset $X$ of $[m]$ whose cardinality decreases by $1$ at each step as the procedure progresses and a strictly decreasing sequence $\yseq$ of elements of $[m]$ whose cardinality at all times is one more than that of $X$.    The initial value of $X$ is $\{w_1,\ldots,w_r\}$,  and that of~$\yseq$ is $y_{r+1}>y_r>\ldots>y_1$.  

The recursive step runs as follows. Let $y$ be the largest element of $\yseq$.   Compare $y$ with the elements of~$X$.   If $y$ is larger than all the elements of $X$,  then we put $w_{r+1}=y$ and the procedure ends.    If not we reset the values of $X$ and $\yseq$ to $X\setminus\{x\}$ and $\yseq\setminus\{y\}$ and repeat,   where $x$ is the smallest element in $X$ not smaller than~$y$.
This recursion will end in at most $r+1$ steps: $X$ will be empty in step $r+1$ and hence the procedure will definitely end at that step (if it gets there at all).

Observe that $w_{r+1}$ is distinct from $w_1$, \ldots, $w_r$.    Indeed,  suppose that the recursion to find it ended in step~$i$ (for some $i$, $1\leq i\leq r+1$).   The value of $y$ during this step (which equals $y_{r+2-i}$ and also $w_{r+1}$) is larger than all the elements of $X$ in this step.     As to other members of $\{w_1, \ldots, w_r\}$,   these must have been removed from $X$ in the earlier steps.  Thus each of them is at least the value of $y$ in some previous step,   and since the value of $y$ keeps strictly decreasing as the steps progress,   it follows that the removed values are all strictly larger than $w_{r+1}$.
\bprop\mylabel{p:procnew}   With notation as above,  suppose that $w_{r+1}=y_s$.   Then if $w^r_1<\ldots<w^r_r$ are the elements $w_1$, \ldots, $w_r$ arranged in increasing order,   then $w^r_1<\ldots<w^r_{s-1}<y_{s}<w^r_s<\ldots<w^r_{r}$.
\eprop
} 
\bprop\mylabel{p:procnew}   With notation as above,  we evidently have:
\begin{itemize}
	\item $a_1<\ldots<a_{k-1}<y_k<a_k<\ldots<a_{\tint-1}$.
	\item $\procw_\tint$ is distinct from $\procw_1$, \ldots, $\procw_{\tint-1}$. \hfill $\Box$
\end{itemize}
\eprop
\bremark  The element $y_j$ in the $\tint$-depth sequence of~$\ssyt$ is just the entry in the right most box of $\tint$-depth $\tint-j$:  ``right most box'' means box in the right most column;   since no two boxes in the same column have the same $\tint$-depth,  this is well defined.  Indeed let $\boxb$ be that box and $e$ its entry.   Clearly $e\leq y_j$.   To show $y_j\leq e$,  first observe that no column to the right of the one containing $\boxb$ has a box of $\tint$-depth $\tint-j$ (by choice of~$\boxb$); secondly that $e$ dominates the entry in any box that is weakly to the Northwest of~$\boxb$ (since $\ssyt$ is a SSYT);  and finally that any box of $\tint$-depth $\tint-j$ strictly South and weakly West of~$\boxb$ can only have an entry that is at most $e$ (for otherwise the $\tint$-depth of~$\boxb$ would exceed~$\tint-j$).
\eremark
\mysubsubsection{Illustration of the procedure above}\mylabel{sss:xbsrrw}
Let $\ssyt$ be the SSYT in~\eqref{e:ssytt'}.
The permutation associated to it is  $83612457$ in one-line notation.  Evidently $\procw_1=8$ and $\mu'_1=3$;   the $2$-depth sequence is $3<4$ and $\procw_2=3$;   the $3$-depth sequence is $3<6<7$ and $\procw_3=6$.
\mysubsubsection{A technical result that will be used later}\mylabel{sss:tlater}   The following lemma will be invoked later on, in Example~\ref{x:deodharcont}.
\blemma\mylabel{l:tlater}  Let $\ssyt$ be a SSYT and $q$ the number of boxes in its right most column.   Let $\ssytp$ be the SSYT obtained from~$\ssyt$ by deleting its last column.     Fix $\tint>q$.    If in the procedure for producing $\procw_p$ (where $\procw$ is the permutation associated to $\ssyt$),  we use the $\tint$-sequence of $\ssytp$ in place of that of~$\ssyt$,   it makes no difference (that is, we still get the same $\procw_p$).
\elemma
\bmyproof   Let $y_1<\ldots<y_\tint$ and $y'_1<\ldots<y'_\tint$ be the $\tint$-depth sequences of~$\ssyt$ and $\ssytp$ respectively, and suppose that $\procw_p=y_k$.
Since the entries in the last column of~$\ssyt$ all belong to $\{\procw_1,\ldots,\procw_{\tint-1}\}$ but, by Proposition~\ref{p:procnew}, $y_k$ does not belong to that set,   it follows that any box of~$\ssyt$ with $y_k$ as its entry belongs to~$\ssytp$.   Thus $y_k=y'_k$.

On the other hand, $y'_j\leq y_j\leq a_{j-1}$ for all $j>k$ (where $a_1<\ldots<a_{\tint-1}$ is the arrangement in increasing order of~$\procw_1$, \ldots, $\procw_{\tint-1}$),  so $k$ is the largest such that $a_{k-1}<y'_k$.~\emyproof
\mysubsection{Association of permutations to LR tableaux}\mylabel{ss:st-to-w}
Recall that the definition in~\S\ref{sss:rlr} of refined LR coefficients refers to a certain association of permutations to LR tableaux.   We describe this association now, after first associating SSYTs to LR tableaux.

Let $T$ be an LR tableau of shape $\nu/\lambda$ and type~$\mu$.
%
If $\nu$ has at most $d$ parts, then so has~$\mu$,  for each entry on row $j$ of $T$ is at most $j$ (for all $j\geq1$).
\mysubsubsection{The SSYT associated to~$T$}\mylabel{sss:st-to-ssyt}  
We associate to~$T$ a SSYT $\ssyt$ of shape $\mu$ as follows.
The entries in row $j$ of $\ssyt$ from left to right are just the row numbers of~$T$ in which the entry~$j$ appears, counted with multiplicity and arranged in weakly increasing order.   That the entries in every column of~$\ssyt$ are strictly increasing downward follows readily from the assumption that the reverse reading word of~$T$ is a ballot sequence:   indeed,  for integers $k\geq1$ and $j\geq2$, if the $k^\textup{th}$ appearance of~$j$ (as we read the reverse reading word from left to right) is in row $r$,    then the $k^\textup{th}$ appearance of $j-1$ is in some row strictly above the $r^\textup{th}$.
\mysubsubsection{The permutation associated to~$T$}\mylabel{sss:skewtab-to-w}
Consider the permutation~$\procw$ associated as in~\S\ref{sss:ssyt-to-w} to the SSYT~$\ssyt$.   We associate~$\procw$ to $T$ itself.
For example, for the skew tableau on the right in~\eqref{e:skewtab}, 
the associated SSYT is the one shown below and the associated permutation is~$51324$:
\begin{equation}\label{e:ssyt}\begin{array}{|c|c|c|c|c|}
	\hline
	1 & 1 & 1 & 3 & 5 \\
	\hline
	2 & 3 \\
	\cline{1-2}
	3 & 4 \\
	\cline{1-2}
	4\\
	\cline{1-1}
\end{array}
	\end{equation}
\mysubsection{$p$-dominance of words}\mylabel{ss:pdom} 
Let $p$: $p_1\geq p_2\geq\ldots$ be a partition.  We denote by $\word(p)$ the word (in the positive integers) that has $p_1$ ones, $p_2$ twos, \ldots\  in succession:    this is just the reverse reading word of the SSYT of shape $p$ all of whose entries in row~$j$ are $j$ (for all $j$).   Note that $\word(p)$ is a ballot sequence.

A word (in the positive integers) is said to be {\em $p$-dominant\/} if when preceded by $\word(p)$ the resulting word is a ballot sequence.  
\bprop\mylabel{p:smallestp} For a given word~$\word$ there is a unique smallest partition $p_\word$ such that $\word$ is $p_\word$-dominant ($p_\word$ is the smallest in the sense that its shape is contained in the shape of any partition~$p$ for which $\word$ is~$p$-dominant).\eprop
\bmyproof  A letter $e>1$ of the given word~$\word$ is said to be a {\em violator\/} if the number of $e-1$ occurring before it does not exceed the number of~$e$ occurring before it.    For $j$ a positive integer,  let $p_j$ be the number of violators in~$\word$ that exceed $j$.  (For example,  $p_1$ is the total number of violators.)
It is elementary to see that the partition $p_1\geq p_2\geq\ldots$ is the unique smallest one for which $\word$ is $p$-dominant.
\emyproof
\ignore{ 
\mysubsubsection{Words attached to a SSYT}\mylabel{sss:ssytw} Let $\ssyt$ be a SSYT.   We denote by $\wordr(\ssyt)$ its reverse reading word (see~\S\ref{sss:rwballotseq}),  and by $\wordc(\ssyt)$ its {\em reverse column word\/} defined as follows:   we read the entries top to bottom in every column beginning with the right most column and ending with the left most.    For the SSYT of~\eqref{e:ssyt},    these words respectively are:
\[ 5311132434 \quad\quad\quad 5311341234   \]
The partition $p_\word$ for either word is the same, namely $3+3+1+1$.  In fact, we have:
\bprop\mylabel{p:pdom}  Let $\ssyt$ be a SSYT and $p$ a partition.  Then $\wordr(\ssyt)$ is $p$-dominant if and only if~$\wordc(\ssyt)$ is so.  \eprop
\bmyproof
For boxes $\boxb_1$ and $\boxb_2$ of~$T$,  the phrase $\boxb_1$ ``occurs before'' $\boxb_2$ in either $\wordc(T)$ or $\wordr(T)$ has the obvious meaning.
Consider a box~$\boxb$ of~$\ssyt$ and let $e$ be its entry.   When we speak 
%
First suppose that $\wordc(\ssyt)$ is $p$-dominant.   To show that $\wordr(\ssyt)$ is $p$-dominant,  it is enough to show the following:  
\begin{itemize}
	\item If a box~$\boxbp$ with entry~$e'$ occurs before $\boxb$ in~$\wordc(\ssyt)$ but not in~$\wordr(\ssyt)$, then $e'\geq e$.
	\item If a box~$\boxbp'$ with entry~$e''$ occurs before $\boxb$ in~$\wordr(\ssyt)$ but not in~$\wordc(\ssyt)$, then $e''< e$.
\end{itemize}
But these two are clear:  any such box~$\boxbp$ is both strictly to the South and strictly to the East of~$\boxb$ and hence $e'>e$; 
 and any such box~$\boxbp'$ is both strictly to the North and strictly to the West of~$\boxb$ and hence $e''<e$. 

 Now suppose that $\wordr(\ssyt)$ is $p$-dominant.    By way of contradiction,  suppose that $\wordc(\ssyt)$ is not so, and that $\boxb$ is a ``violating'' one:  that is,   the number of 
\emyproof
} 
\mysubsubsection{Weights of words in~$[d]$}\mylabel{sss:wtword}   Let $w$ be a word in~$[d]$.    The {\em weight of~$w$\/},   denoted $\weight(w)$,   is defined to be the weight 
\mysubsubsection{The words $\wordr$ and $\wordc$ attached to a SSST}\mylabel{sss:wst} Let $T$ be a SSST. 
We have already defined its reverse reading word~$\wordr(T)$ in~\S\ref{sss:rwballotseq}.
We now define its {\em reverse column word\/}, denoted $\wordc(T)$, as follows:   we read the entries top to bottom in every column beginning with the right most column and ending with the left most.  
For the SSST in~\eqref{e:skewtab},    the reverse column words respectively are
$3227536341$ and $1112324131$.

For the SSST on the left in~\eqref{e:skewtab},  the partitions~$p_\word$ attached (as in Proposition~\ref{p:smallestp}) to its words $\wordr$ and~$\wordc$ turn out to be the same, namely $5+3+2+2+1+1$.  For the SSST on the right in~\eqref{e:skewtab},  both $\wordr$ and $\wordc$ are ballot sequences (so $p_\word$ is empty for both).   Indeed we have:
\bprop\mylabel{p:pdom}  Let $T$ be a SSST and $p$ a partition.  Then $\wordr(T)$ is $p$-dominant if and only if~$\wordc(T)$ is so.\eprop
\bremark\mylabel{r:p:pdom}     This statement is well known at least in the case of a SSYT (see, e.g.,  \cite[Exercise~5.2.4]{lothaire}).   A proof from first principles is given below for the sake of completeness.
\eremark

\bmyproof 
For boxes $\boxb_1$ and $\boxb_2$ of~$T$,  the phrase $\boxb_1$ ``occurs before'' $\boxb_2$ in $\wordc(T)$ (respectively $\wordr(T)$) has the obvious meaning.
We let $\boxb$ be an arbitrary box in~$T$.  Its position is denoted by $(r,c)$ and entry by~$e$.
\begin{enumerate}\item  Let $\boxbp$ be a box that occurs before~$\boxb$ in~$\wordc(T)$ but not in~$\wordr(T)$.  Let its position be denoted by~$(r',c')$ and entry by~$e'$.     Then  $r<r'$, $c<c'$ and, since $T$ is semi-standard, $e<e'$.
\item  Let $\boxb''$ be a box that occurs before~$\boxb$ in~$\wordr(T)$ but not in~$\wordc(T)$.  Let its position be denoted by~$(r'',c'')$ and entry by~$e''$.     Then $r''<r$, $c''<c$ and, since $T$ is semi-standard, $e''<e$.
\end{enumerate}
The following figure depicts the situation:
\begin{center}
	\begin{tikzpicture}[scale=0.5]
	\filldraw[color=red!60, fill=red!5, very thick](0,5)--(0,1)--(-4,1);
	\filldraw[color=blue!60, fill=blue!5, very thick](5,0)--(1,0)--(1,-4);
	\draw[thick] (0,0)--(1,0)--(1,1)--(0,1)--(0,0);
	\node at (0.5,0.5){$\boxb$};
		\node at (4,-1){\textup{Region of~$\boxb'$}};
		\node at (-3,2){\textup{Region of~$\boxb''$}};
\end{tikzpicture}
\end{center}

Suppose first that $\wordc(T)$ is $p$-dominant.     Consider the contributions to the words $\wordr(T)$ and $\wordc(T)$ of an arbitrarily fixed box~$\boxb$ in~$T$.   With notation as above,   observe that no box $\boxb'$ has $e-1$ as an entry  and no box~$\boxb''$ has $e$ as an entry.   Thus,  letting $m_r$ and $n_r$ (respectively $m_c$ and $n_c$) denote respectively the number of occurrences of~$e$ and $e-1$ (strictly) before $\boxb$ in $\wordr(T)$ (respectively $\wordc(T)$),  we have $m_r\leq m_c$ and $n_c\leq n_r$.  Since $m_c\leq n_c$ by $p$-dominance of~$\wordc(T)$,  we have $m_r\leq m_c\leq n_c\leq n_r$,  so $\wordr(T)$ is $p$-dominant too.

Now suppose that $\wordr(T)$ is $p$-dominant.     By way of contradiction,   suppose that $\wordc(T)$ is not $p$-dominant.  Choose a box $\boxb$ in~$T$ which ``violates'' the $p$-dominance of~$\wordc(T)$,   meaning that (with notation as above) $n_c<m_c$.   Since no box of type $\boxb'$ or $\boxb''$ can have an entry equal to~$e$---we have $e''<e<e'$---it follows that $m_c=m_r$.  

Consider a box of type $\boxb''$ with entry equal to~$e-1$.   Let us denote by $\boxb''_1$ any such box and 
suppose that there are $k$ such boxes.
Then $n_r=n_c+k$, since $e'>e$.
The entry in the box just below a box~$\boxb''_1$ must be~$e$ (since such a box is weakly North and strictly West of~$\boxb$ on the one hand,  but on the other hand its entry must be strictly larger than $e-1$).
Thus all the $k$ boxes $\boxb''_1$ must occur in row $r-1$,   
and $T$ looks like:
\begin{center}
	\begin{tikzpicture}[scale=0.5]
	\filldraw[color=red!60, fill=red!5, very thick](0,4)--(0,1)--(-4,1);
	\filldraw[color=blue!60, fill=blue!5, very thick](5,0)--(1,0)--(1,-3);
	\draw[thick] (0,0)--(1,0)--(1,1)--(0,1)--(0,0);
	\node at (0.5,0.5){$e$};
	\draw[thick] (-1,0)--(0,0)--(0,1)--(-1,1)--(-1,0);
	\node at (-0.5,0.5){$e$};
		\node at (-1.5,0.5){\ldots};
	\draw[thick] (-3,0)--(-2,0)--(-2,1)--(-3,1)--(-3,0);
	\node at (-2.5,0.5){$e$};
		\node at (-3.5,-1.0){$\boxb_1$};
		\draw[thick, ->](-3.5,-0.5)--(-2.5,0);
	\draw[thick] (-1,1)--(0,1)--(0,2)--(-1,2)--(-1,1);
	\node at (-0.5,1.5){$f$};
		\node at (-1.5,1.5){\ldots};
	\draw[thick] (-3,1)--(-2,1)--(-2,2)--(-3,2)--(-3,1);
	\node at (-2.5,1.5){$f$};
	\draw[thick] (-4,1)--(-3,1)--(-3,2)--(-4,2)--(-4,1);
	\node at (-3.5,1.5){$g$};
		\node at (4,-1){\textup{Region of~$\boxb'$}};
		\node at (-3,3){\textup{Region of~$\boxb''$}};
	\node at (4,3){$f:=e-1$};
	\node at (4,2){$g<f$};
\end{tikzpicture}
\end{center}

Now let $\boxb_1$ be the box in~$T$ that is $k$ boxes to the left of~$\boxb$.
Let us count the number of entries equal to $e$ (respectively $e-1$) that occur before $\boxb_1$ in~$\wordr(T)$.     This count equals $m_r+k=m_c+k$ (respectively, $n_r=n_c+k$).     We have $n_c+k<m_c+k$ (since $n_c<m_c$ by choice of~$\boxb$).   But this means that the box $\boxb_1$ violates the $p$-dominance of~$\wordr(T)$,  a contradiction.
%
%
\emyproof
\mysubsection{Deconstructing a SSST}\mylabel{ss:dssst}
Let $T$ be a SSST of shape~$\nu/\lambda$.   As before,  we think of~$\lambda$ as being fixed and $\nu$ as varying.
For $k$ a positive integer: 
\begin{itemize}
	\item  Let $n_r(k)$ denote the number of times $k$ appears in row $r$.
	\item  Consider the boxes of $T$ belonging to $\lambda$ and those with entries not exceeding~$k$.   Together they form a Young shape.  Denote by $\lambda^k$ this shape as well as the corresponding partition.   It is convenient to set $\lambda^0=\lambda$.    Observe that \begin{equation}\label{e:lchain}\lambda=\lambda^0\subseteq\lambda^1\subseteq\lambda^2\subseteq\ldots\end{equation}
where $\subseteq$ between shapes means that the former is contained in the latter.    We have $\lambda^k_r-\lambda^{k-1}_r=n_r(k)$.
	\item Denote by $\wordk(T)$  the word comprising the row numbers of~$T$ in which $k$ appears,    listed with multiplicity and in weakly decreasing order.  
		In terms of the integers~$n_r(k)$,  we have $\wordk(T)=   \ldots 2^{n_2(k)}1^{n_1(k)}$.
\end{itemize}
The hypothesis that $T$ is semistandard puts a constraint on the sequence of shapes that can possibly arise as~\eqref{e:lchain}.  Indeed, the fact that an of entry of~$T$ is strictly larger than the one vertically just above it (if the latter happens to exist) means precisely that no two boxes in $\lambda^k\setminus\lambda^{k-1}$ are in the same column, or, in other words:
\begin{equation}\label{e:ldom0} \lambda^k_r\leq\lambda^{k-1}_{r-1} \quad\quad \forall r>1\quad \forall k\geq1  \end{equation}
	In terms of $\lambda$ and $n_r(k)$,  this can also be expressed as the following set of conditions:
\begin{equation}\label{e:ldom}
	\lambda_r+n_r(1)+\cdots+n_r(k)\leq \lambda_{r-1}+n_{r-1}(1)+\cdots+n_{r-1}(k-1)\quad\quad \forall r>1\quad \forall k\geq 1
\end{equation}
\mysubsubsection{The position word~$\wordp$ and its $\lambda$-dominance}\mylabel{sss:wpldom}
To see what~\eqref{e:ldom} translates to in terms of the words~$\wordk(T)$,   let us define the {\em position word\/} of~$T$, denoted $\wordp(T)$, to be the concatenation $w_1(k)w_2(k)\ldots$.  For example,  the position words of the SSST in~\eqref{e:skewtab} are, respectively, $5111334342$ and $5311132434$.   It is readily seen that~\eqref{e:ldom} is equivalent to the $\lambda$-dominance of the word~$\wordp(T)$  (in the sense of~\S\ref{ss:pdom}).
\mysubsubsection{Recovering the SSST~$T$}\mylabel{sss:trcvr}  Evidently the SSST $T$ can be recovered from the collection of integers~$n_r(k)$ (presuming knowledge of the fixed partition~$\lambda$).    Thus it can be recovered either from the sequence~\eqref{e:lchain} of increasing shapes or from the sequence $\word_1(T)$, $\word_2(T)$, \ldots\ of words.   Moreover,  if either the sequence~\eqref{e:lchain} satisfies the constraint~\eqref{e:ldom0} or, equivalently, if the sequence $\word_1(T)$, $\word_2(T)$, \ldots\ is such that~$\wordp(T)$ is $\lambda$-dominant,  then there exists a corresponding~$T$.
\mysubsubsection{Bijection between $\skewtabmlwd$ and $\ssytsetmlwd$}\mylabel{sss:isoTS}  As preparation for the proof in~\S\ref{ss:pfkktab} below of the tableau version of the decomposition rule~\eqref{e:kktab} of KK modules,    we apply the observations above to the case when $T$ is LR.  

Fix notation as in~\S\ref{ss:tabkk}.
Let $\skewtabmld$ denote the subset of $\skewtabml$ consisting of those elements~$T$ such that $\nu(T)$ has at most $d$ parts.
Let $\ssytsetml$ denote those SSYT of shape $\mu$ whose column word is~$\lambda$-dominant (in the sense of~\S\ref{ss:pdom}),
and let $\ssytsetmlw$ be the subset of those elements of~$\ssytsetml$ for which the associated permutation~$\procw$ (as in~\S\ref{sss:ssyt-to-w}) satisfies~$\procw\leq w$.
Put:
\[ \ssytsetmld:=\{S\in\ssytsetml\st \textup{no entry of $S$ exceeds~$d$}\}
\quad\quad\quad
 \ssytsetmlwd:=\ssytsetmlw\cap\ssytsetmld
	\]
%
	The weight of a SSYT with entries from~$[d]$ is its weight thought of as a SSST (see~\ref{sss:type}).
\bprop\mylabel{p:bijTS}   Let $T$ be an element of~$\skewtabml$ and $\ssyt$ the SSYT attached to~$T$ as in~\S\ref{sss:st-to-ssyt}.   The association $T\mapsto \ssyt$ gives a bijection between  $\skewtabml$ and $\ssytsetml$, under which  $\nu(T)=\lambda+\weight(\ssyt)$,   and which also restricts to a bijection between the pairs $\skewtabmlw$, $\ssytsetmlw$  and $\skewtabmlwd$, $\ssytsetmlwd$.
\eprop
\bmyproof
We first show that $T\mapsto \ssyt$ gives a bijection between $\skewtabml$ and $\ssytsetml$.     From Proposition~\ref{p:pdom} it follows that the $\lambda$-dominance of~$\wordr(\ssyt)$ and $\wordc(\ssyt)$ are equivalent,  so
\[\ssytsetml =\{\textup{$\ssyt$ is a SSYT of shape~$\mu$}\st \textup{$\wordr(\ssyt)$ is $\lambda$-dominant}\}
	\]
It is easy to see from their definitions that the words $\wordr(\ssyt)$ and $\wordp(T)$ are the same.  
Thus, from \S\ref{sss:trcvr}, we conclude: 
\begin{itemize}
	\item $\wordp(T)=\wordr(\ssyt)$ is~$\lambda$-dominant,  so $\ssyt$ belongs to~$\ssytsetml$.
	\item  
The sequence $\word_1(T)$, $\word_2(T)$, \ldots\ defined in~\S\ref{ss:dssst} and hence $T$ itself can be recovered readily from $S$ by reading the entries in every row of~$S$ from right to left.   This shows that $T\mapsto S$ is one-to-one.
	\item 
		Given $\ssyt'$ in~$\ssytsetml$,  the $\lambda$-dominance of $\wordr(\ssyt')$ means that there exists a skew tableau $T'$ of shape $\nu/\lambda$ (for some~$\nu$) that corresponds to it (in the sense of~\S\ref{sss:trcvr}).  The fact that the entries along any column of~$\ssyt'$ are strictly increasing downwards translates to the fact that the corresponding $T'$ as above is LR,  so $T'$ belongs to $\skewtabml$ and $T'\mapsto S'$.    This shows that $T\mapsto S$ is surjective. 
\end{itemize}
\noindent This finishes the proof that $T\mapsto\ssyt$ gives a bijection from~$\skewtabml$ to $\ssytsetml$.

It is clear from the description of the association~$T\mapsto\ssyt$ that $\ssyt$ has type~$\mu$ and that $\lambda+\nu(T)=\weight(\ssyt)$.

The association of a permutation to an LR tableau proceeds via the SSYT attached to it, so it immediately follows that $T\mapsto\ssyt$ gives a bijection from $\skewtabmlw$ to $\ssytsetmlw$.

Finally,  the number of parts of~$\nu(T)$ on the one hand and the maximum value of an entry in~$\ssyt$ on the other are upper bounds for each other under~$T\mapsto S$,   so we get a bijection between $\skewtabmlwd$ and $\ssytsetmlwd$.
\emyproof
\ignore{
Let $T$ be a SSST of shape~$\nu/\lambda$.   For the purposes of this subsection,  it is convenient to think of~$\lambda$ as being fixed but $\nu$ as varying.
To $T$ we attach its {\em position word\/} $\wordp(T)$,  a word in the positive integers, as follows.  For $k$ a positive integer,
let $\wordk(T)$  denote the word comprising the row numbers of~$T$ in which $k$ appears,    listed with multiplicity and in weakly decreasing order.   We define $\wordp(T)$ to be the concatenation $\word_1(T)\word_2(T)\ldots$.     
The following claims are obvious:
\begin{enumerate}
	\item\label{i:wordp:k} The number of times $k$ appears in~$\wordp(T)$ equals $\nu_k-\lambda_k$ (for any positive integer~$k$).
	\item\label{i:wordp:lrs}   If $T$ is LR and $\ssyt$ the SSYT attached to it as in~\S\ref{sss:st-to-ssyt},   then $\wordp(T)$ is just the reverse reading word of~$\ssyt$.
	\item\label{i:wordp:d}  If $k$ is the number of parts in~$\nu$, then $\wordp(T)$ is a word in~$[k]$.
\end{enumerate}
\bprop\mylabel{p:wordp:ldom}
	$\wordp(T)$ is $\lambda$-dominant (in the sense of~\S\ref{ss:pdom}).
\eprop
\bmyproof   Let $n_r(k)$ denote the number of entries that are equal to~$k$ in row~$r$ of~$T$.    The semi-standardness of $T$ is equivalent to the following:
\begin{equation}\label{e:ssT}
	abc
\end{equation}
\emyproof
\bremark\mylabel{r:p:wp}   \begin{enumerate}
	\item \label{i:r:wp:any}  As the proof shows,  if the entries of each $\wordk(T)$ are permuted among themselves,   the resulting word is still~$\lambda$-dominant.
\end{enumerate}
\eremark


Let $\nu/\lambda$ be the shape of the skew tableau~$T$,  with $\lambda$ having at most $n$ parts (so $\lambda$ represents also a dominant integral weight).     
\mysubsubsection{$\lambda$-dominance}\mylabel{sss:ldom}  Let $\nu/\lambda$ be the shape of the skew tableau~$T$,  with $\lambda$ having at most $n$ parts (so $\lambda$ represents also a dominant integral weight).     
An SSYT is said to be {\em $\lambda$-dominant\/} if interpreted as a concatenation of LS paths (see~\ref{ss:pathtab}) it is $\lambda$-dominant.
\bprop\mylabel{p:ldom}  The SSYT $T'$ associated to~$T$ in~\S\ref{sss:st-to-ssyt} is $\lambda$-dominant.   Conversely, given a $\lambda$-dominant SSYT $T'$ of shape~$\mu$,   there exists a skew tableau $T$ of shape $\nu/\lambda$ (where $\nu$ is allowed to vary) to which $T'$ is associated.
\eprop
}
\mysubsection{Proof of the tableau KK decomposition rule of~\S\ref{ss:tabkk}}\mylabel{ss:pfkktab}    The decomposition rule~\eqref{e:kktab} in terms of tableaux can be derived, as we now show,  from the general decomposition rule~\eqref{e:decompkos} for KK-modules in~\S\ref{s:decompkos}.  The derivation consists of stringing together three bijections that preserve invariants.   

The first of these is the bijection between~$\skewtabmld$ and $\ssytsetmld$ of Proposition~\ref{p:bijTS}.   The second and third bijections are from the appendix:  
by Corollary~\ref{c:bijSPstd},   we may identify $\ssytsetmd$,  the set of SSYT of shape~$\mu$ with entries from~$[d]$,  with $\pathstd$, the set of standard concatenations of LS paths as in~\S\ref{ss:slla};  and, finally,   there is the crystal isomorphism~$\cryiso$ of~\S\ref{ss:crystal} between the set $\pathsm$ of LS paths of shape~$\mu$ and~$\pathstd$.

In the subsection below, the good properties required of the second bijection are established.  For the first bijection, this was done in Proposition~\ref{p:bijTS}. 
As for the crystal isomorphism~$\cryiso$,  it preserves end points and $\lambda$-dominance as shown in Proposition~\ref{p:cryiso0};  and the
minimal element in the initial direction
of~$\pi$ in~$\pathsm$ is the
initial element of the minimal standard lift
of~$\cryiso\pi$  as shown in Proposition~\ref{p:cryiso}. 

The final upshot is a bijection $T\leftrightarrow\pi$ between $\skewtabmld$ on the one hand and $\pathsml$ on the other such that (a)~$\nu(T)$ equals the end point $\pi(1)$ and (b)~the permutation $\procw$ attached to~$T$ as in~\S\ref{sss:skewtab-to-w} equals the minimal element in the initial direction of~$\pi$.   This will finish the proof of the tableau decomposition rule~\eqref{e:kktab}.

%
\mysubsubsection{Good properties of the bijection of Corollary~\ref{c:bijSPstd}}\mylabel{sss:bijSPgood}
\bprop\mylabel{p:SPstd}
Under the identification between~$\ssytsetmd$ and $\pathstd$ of Corollary~\ref{c:bijSPstd},  let $\ssyt$ in~$\ssytsetmd$ correspond to~$\theta$ in~$\pathstd$.  Then:
\begin{enumerate}
	\item\label{i:wt=end} 
		The weight $\weight(\ssyt)$ of~$\ssyt$ equals the end point $\theta(1)$ of the path~$\theta$.
	\item\label{i:ueqv} The permutation~$\procw$ associated to~$\ssyt$ by the procedure of~\S\ref{sss:ssyt-to-w} equals the initial element of the minimal standard lift of~$\theta$.
	\item\label{i:ldomeq} The column word $\wordc(\ssyt)$ of~$\ssyt$ is $\lambda$-dominant (in the sense of~\S\ref{ss:pdom}) if and only if the path $\theta$ is $\lamdba$-dominant. 
\end{enumerate}\eprop
\bmyproof Item~\eqref{i:wt=end} is immediate from the definitions.  As for item~\eqref{i:ueqv},  the whole of~\S\ref{s:proc-proof} is devoted to its proof.  

Turning to item~\eqref{i:ldomeq},  we first prove the ``only if part''. Let $c$ denote the number of columns in the shape of~$S$,  let $r_j$ denote the number of boxes in column~$j$ of~$S$ (for $1\leq j\leq c$), and let $\word'_j$ denote the word $s_{1j}\ldots s_{r_jj}$  (where, as for a matrix, $s_{ab}$ denotes the entry of~$\ssyt$ in row~$a$ and column~$b$).  The word $\wordc(\ssyt)$ is,  by definition, $\word'_c\word'_{c-1}\cdots\word'_1$.   Its $\lambda$-dominance clearly implies that of any left subword of it,   in particular that of the subwords $\word'_c$, $\word'_c\word'_{c-1}$, \ldots, $\word'_c\word'_{c-1}\cdots\word'_2$, and $\word'_c\word'_{c-1}\cdots\word'_2\word'_1=\wordc(\ssyt)$.  This in turn implies that the weights $\lambda+\weight(\word'_c)$, $\lambda+\weight(\word'_c\word'_{c-1})$, \ldots, $\lambda+\weight(\word'_c\word'_{c-1}\cdots\word'_2)$, and $\lambda+\weight(\word'_c\word'_{c-1}\cdots\word'_2\word'_1)=\lambda+\weight(\wordc(\ssyt))$ are all dominant.  
But the dominance of these $c$ weights is, as is readily seen, precisely equivalent to the $\lambda$-dominance of~$\theta$.

For the ``if part'',  we first make an observation (whose elementary proof we skip).    Suppose that a word $\word$ in~$[d]$ is a concatenation $\word_1\word_2$ of words~$\word_1$ and $\word_2$ such that $\word_1$ is $\lambda$-dominant, $\word_2$ is weakly increasing (left to right), and $\lambda+\weight(\word)$ is dominant.   Then $\word$ is $\lambda$-dominant.

The $\lambda$-dominance of $\theta$ implies that
$\lambda+\weight(\word'_c)$, $\lambda+\weight(\word'_c\word'_{c-1})$, \ldots, $\lambda+\weight(\word'_c\word'_{c-1}\cdots\word'_2)$, and $\lambda+\weight(\word'_c\word'_{c-1}\cdots\word'_2\word'_1)=\lambda+\weight(\wordc(\ssyt))$ are all dominant.   Since each $\word'_j$ is strictly increasing,  we conclude using the observation that $\wordc(\ssyt)$ is $\lambda$-dominant.
\emyproof

\mysection{An important property of the procedure of~\S\ref{sss:ssyt-to-w}}\mylabel{s:proc-proof}
\noindent
Let $\ssyt$ be an SSYT (see~\S\ref{ss:ssyt-and-perm}) none of whose entries exceeds $d$, and $\procw$ the permutation of~$[d]$ obtained by application to~$\ssyt$ of the procedure of \S\ref{sss:ssyt-to-w}.    As explained in~\S\ref{ss:slla} (see, in particular, Corollary~\ref{c:bijSPstd})  such SSYTs may be identified as certain standard concatenations of LS paths whose shapes are fundamental weights (for~$\lieg=\sld$).    In what follows,  we will use the notation for an SSYT to denote also the corresponding standard concatenation of paths. Let $\imsl$ be the initial element of the minimal standard lift of~$\ssyt$ (\S\ref{ss:standard}). 

The purpose of this section is to show that $\procw=\imsl$.   The proof is given in~\S\ref{ss:pfvleqw} and~\S\ref{ss:uleqv} after preparations in the earlier subsections.

The procedure of~\S\ref{sss:ssyt-to-w} seems to be quite different from and more efficient than a repeated application of Deodhar's lemma (Example~\ref{x:deodhar}) to compute the initial element of the minimal standard lift~$\imsl$.    Besides,   the justification we give in Example~\ref{x:deodharcont} of the recipe of~Example~\ref{x:deodhar}  is itself based on the result of this section (that $\procw=\imsl$).
\mysubsection{Notation relating to permutations}\mylabel{ss:notn}
Let $x$ be a permutation and let $x_1x_2\ldots$ denote its one-line notation.

We call $\{i\st x_i>x_{i+1}\}$ the {\em descent set\/} of $x$.
We say that $x$ has {\em only $r$ significant elements\/} if its descent set is contained in~$[r]$, or, in other words, if the sequence $x_{r+1}x_{r+2}\ldots$
is increasing.   E.g.,  the only permutation that has zero significant elements is the identity.      

For $s$ an integer,  let $\vecx^s$ denote the sequence $x^s_1<\ldots<x^s_s$ of the first $s$ elements of~$x$ (namely $x_1$, \ldots, $x_s$) arranged in increasing order. 
\mysubsubsection{On the tableau criterion for Bruhat order}\mylabel{sss:tabbruhat}
Recall the following ``tableau criterion'' for comparability in Bruhat order of two permutations:    
		$x\leq z$ if and only if $\vecx^s\leq \vecz^s$
		for all~$s$,  where $\vecx^s\leq \vecz^s$ is short hand for $x^s_j\leq z^s_j$ for all $1\leq j\leq s$.
		\blemma\mylabel{l:bb} \textup{(\cite[Corollary~(5)]{bb})}  For $x\leq z$,  it suffices that $\vecx^s\leq \vecz^s$ holds for either (a)~all $s$ in the descent set of $x$,
		or (b)~all $s$ not in the descent set of~$z$.
		\elemma
\mysubsection{An example}\mylabel{ss:deodharcont}
For $x$ a permutation,    we denote by $\ltrunc{x}{r}$ the permutation obtained from~$x$ by rearranging the first $r$ elements in its one-line notation in increasing order.   In other words,  $\ltrunc{x}{r}$ is the permutation whose one-line notation is $x^r_1\ldots x^r_r x_{r+1}x_{r+2}\ldots$.
\blemma \mylabel{l:ltrunc} $\ltrunc{x}{s}\leq\ltrunc{x}{r}\leq \ltrunc{x}{1}=x$ for $s\geq r\geq 1$.
\elemma
\bmyproof  Put $y=\ltrunc{x}{s}$ and $z=\ltrunc{x}{r}$.  The descent set of $y$ is contained in $\{s,s+1,\ldots\}$.   For any $t\geq s$,  we have $\vecy^t=\vecz^t$, so it follows from Lemma~\ref{l:bb} that $y\leq z$.  Observe that $\ltrunc{x}{1}=x$.
\emyproof

\noindent
Given a permutation $x$ of~$[n]$ and an integer~$r\leq n$,  we let $\ssyt(r,x)$ denote the SSYT constructed as follows:    it has $n-r+1$ columns;  column~$j$ (counting from the left) has $n+1-j$ boxes and its entries are  the first $n+1-j$ entries of $x$ arranged in increasing order. E.g.,  if $x$ is the permutation of~$[5]$ with one-line notation $45312$,  then $\ssyt(x,3)$ is:
\[
	\begin{array}{|c|c|c|}
		\hline
		1 & 1 & 3 \\
		\hline
		2 & 3 & 4 \\
		\hline
		3 & 4 & 5 \\
		\hline
		4 & 5\\
		\cline{1-2}
		5 \\
		\cline{1-1}
	\end{array}
	\]
\blemma \mylabel{l:ssytrx} The initial element of the minimal standard lift of~$\ssyt(r,x)$ is $\ltrunc{x}{r}$.
\elemma
\bmyproof  Let $a$ be this initial element. By an induction argument,  we may assume that $\ltrunc{x}{s}$ is the initial element of the minimal standard lift of~$\ssyt(s,x)$ for $s>r$.    Thus $z:=\ltrunc{x}{r+1}\leq a$.  Since the first $r$ elements of $a$ match the respective ones of $y:=\ltrunc{x}{r}$, it follows in particular that $\vecy^r\leq \veca^r$.   Since the descent set of $y$ is contained in~$\{r, r+1, \ldots\}$,  and $\vecy^s=\vecz^s\leq \veca^s$ for $s>r$,   it follows from Lemma~\ref{l:bb} that $y\leq a$.

On the other hand, evidently $\ltrunc{x}{n}\leq\ltrunc{x}{n-1}\leq\ldots\leq\ltrunc{x}{r}$ is a standard lift of~$\ssyt(r,x)$,  so $a\leq\ltrunc{x}{r}=y$.   
\emyproof
\bexample\mylabel{x:deodharcont} Let notation be fixed as in Example~\ref{x:deodhar}.    We described there a procedure for determining $\tau:=J_{\sigma W_r}(w)$ without however providing a justification for it.     We now provide such a justification as an application of the main result of this section ($\procw=\imsl$).

Let $\ssytp$ denote the SSYT $\ssyt(r,w)$ and
$\ssyt$ the SSYT obtained by attaching to~$S'$ on the right a column with $r$ boxes whose entries from top to bottom are $\sigma_1$, \ldots, $\sigma_r$.
For the values $n=6$, $r=3$, $\sigma=246135$, and $w=145362$ used as an illustration in Example~\ref{x:deodhar},  $\ssyt$ is:
\begin{equation}\label{e:ssytjsigmaw}
	\begin{array}{|c|c|c|c|c|}
		\hline
		1 & 1 & 1 & 1 & 2\\
		\hline
		1 & 3 & 3 & 4 & 4 \\
		\hline
		3 & 4 & 4 & 5 & 6 \\
		\hline
		4 & 5 & 5\\
		\cline{1-3}
		5 & 6 \\
		\cline{1-2}
		6 \\
		\cline{1-1}
	\end{array}
\end{equation}
By Lemma~\ref{l:ssytrx},   the initial element of the minimal standard lift of~$\ssytp$ is $\ltrunc{w}{r}=w$,  so the initial element~$\imsl$ of the minimal standard lift of~$\ssyt$ is the least element having the following two properties:  $w\leq \imsl$ and the first $r$ elements of~$\imsl$ (in its one-line notation) are $\sigma_1$, \ldots, $\sigma_r$, in that order.

Now, $\tau$ is the least element having the two properties:  $w\leq \tau$ and $\tau W_r=\sigma W_r$.  
Evidently $\ltrunc{\tau}{r}W_r=\tau W_r=\sigma W_r$,  and, by
Lemma~\ref{l:ltrunc}, $w=\ltrunc{w}{r}\leq\ltrunc{\tau}{r}\leq\tau$.
So $\ltrunc{\tau}{r}=\tau$ and the first $r$ elements of $\tau$ are $\sigma_1$, \ldots, $\sigma_r$, in that order.
This means $\tau=\imsl$.

Thus, by the main result of this section,  the element $\procw$ obtained by applying the procedure of~\S\ref{sss:ssyt-to-w} to~$\ssyt$ equals~$\tau$.     It is easily seen that $\procw_j=\sigma_j$, for $1\leq j\leq r$.    For $j>r$, to determine $\procw_j$,   we may use, by Lemma~\ref{l:tlater}, the $j$-depth sequence $y'_1<\ldots<y'_j$ of~$\ssyt'$ instead of that of~$\ssyt$.    The entries in the column with $j$ boxes of~$\ssyt'$ being $w^j_1<\ldots<w^j_j$,  it is clear that $w^j_i\leq y'_i$ for $1\leq i\leq j$.    On the other hand,  since each $y'_i$ must be an entry in one of the columns of~$\ssyt'$ with at most $j$ boxes, it follows that $w^j_i=y'_i$ for every $i$.  

This completes the justification of the recipe of~Exercise~\ref{x:deodhar} to compute~$J_{\sigma W_r}(w)$.
\eexample
\mysubsection{Truncations of permutations and SSYTs}\mylabel{ss:trunc}
For $r$ an integer,  let $x^{(r)}$ denote the permutation obtained from $x$ by rearranging its elements in position $r+1$ and beyond in increasing order.    We call $x^{(r)}$ the {\em $r$-truncation of~$x$\/}. Evidently $x^{(r)}$ has only $r$ significant elements.     As an easy consequence of Lemma~\ref{l:bb}, we have: 
\blemma\mylabel{l:procbasic} Suppose that $x\leq z$.   Then $x^{(r)}\leq z^{(r)}\leq z$.
\elemma

For $r$ an integer, let $\ssytr$ denote the SSYT obtained by taking the first $r$ rows of~$\ssyt$:   if $\ssyt$ has at most $r$ rows, then $\ssytr$ is all of $\ssyt$.       We call $\ssytr$ the {\em $r$-truncation of~$\ssyt$\/}.   Let $\imslsr$ denote the initial element of the minimal standard lift of~$\ssytr$.
\bprop\mylabel{p:procbasic1} 
		Every permutation in the minimal standard lift of $\ssytr$ has only $r$ significant elements.  In particular, if $\ssyt$ has at most $r$ rows, then $\imsl$ has only $r$ significant elements: $\imsltruncr=\imsl$.
		\eprop
\bmyproof 
We use Lemma~\ref{l:procbasic} to observe that the $r$-truncation of any standard lift of~$\ssytr$ continues to be a standard lift,  and moreover that the $r$-truncation of the minimal standard lift is itself.
\emyproof
		\bprop\mylabel{p:procbasic2} $\imslsr=\imsltruncr$.\eprop
			\bmyproof
Using Lemma~\ref{l:procbasic} again, we observe that the $r$-truncation of any standard lift of~$\ssyt$ gives a standard lift of $\ssytr$.  Thus $\imslsr\leq\imsltruncr$.

Let $\sigma_1\leq\ldots\leq \sigma_k=\imslsr$ be the minimal standard lift of~$\ssytr$.    By Proposition~\ref{p:procbasic1},   there are only $r$ significant elements in every $\sigma_j$.   We will construct a standard lift $\sigmat_1\leq\ldots\leq \sigmat_k$ of~$\ssyt$ whose $r$-truncation is $\sigma_1\leq\ldots\leq\sigma_k$.   It will then follow that $\imsl\leq\sigmat_k$,  and so,  by Lemma~\ref{l:procbasic},  $\imsltruncr\leq \sigmat_k^{(r)}=\sigma_k=\imslsr$.

To construct $\sigmat_j$,   proceed as follows.   The first $r$ elements of $\sigmat_j$ are the same as those of~$\sigma_j$.    The first entries of $\sigmat_j$ also match the entries top-downwards in column $j$ of~$\ssyt$ (until the latter entries are exhausted).   The remaining entries of $\sigmat_j$ are arranged in decreasing order.   Criterion~(b) of Lemma~\ref{l:bb} is useful to verify $\sigmat_j\leq\sigmat_{j+1}$.
\emyproof
\mysubsection{The main part of the proof (that $\procw=\imsl$)}\mylabel{ss:pfmain}
Let $k$ be the number of columns in $S$.  For $i$, $1\leq i\leq k$, let $\ssyt[i]$ denote the SSYT consisting only of the first $i$ columns (from the left) of $\ssyt$,     and $\procw[i]$ the permutation obtained by running the procedure of~\S\ref{sss:ssyt-to-w} on~$\ssyt[i]$.  By the description of the procedure, it is clear that running the procedure on $\ssyt[i]^{(\tint)}$ yields $\procw[i]^{(\tint)}$.

Let $p$ denote a positive integer.  We proceed by induction on~$\tint$ to show the following three assertions.\footnote{It is only assertion~\eqref{i:wseq} that we are really interested in.   Once we have it,   it follows rather easily that $\imsl\leq\procw$ (see~\S\ref{ss:pfvleqw}).    The other two assertions are technical devices that facilitate the proof of~\eqref{i:wseq}.}

\begin{enumerate}[label=(a), ref=a]
	\item\label{i:wseq} $\procw[1]^{(\tint)}\leq\procw[2]^{(\tint)}\leq \ldots\leq\procw[k]^{(\tint)}=\procw^{(\tint)}$.
\end{enumerate}
Consider a rectangular grid of boxes with $\tint$ boxes in every column and $k$ boxes in every row.   Suppose we fill the boxes in column $i$ of this grid by the first $\tint$ entries of $\procw[i]$ 
in increasing order.    It is clear from item~\eqref{i:wseq} above that we then get a SSYT (see Lemma~\ref{l:bb}).   Let $\ssyttint$ be the SSYT whose first $\tint$ rows are this rectangular SSYT and whose rows $\tint+1$ and beyond are  the same as the corresponding ones of~$\ssyt$.  For $i$, $1\leq i\leq k$,  we let $\ssyttint[i]$ denote the SSYT consisting of only the first $i$ columns of $\ssyttint$.

For example, on the left in the following display is shown $\ssyt_2$ for $\ssyt$ as in~\eqref{e:ssytt'};
and on the right is shown $\ssyt_4$ for $\ssyt$ as in~\eqref{e:ssytjsigmaw}:
\[
	\begin{array}{|c|c|c|c|}
		\hline
		1 & 3 & 3 & 3 \\
		\hline
		2 & 4 & 6 & 8\\
		\hline
		7\\
		\cline{1-1}
	\end{array}
	\quad\quad\quad\quad
	\begin{array}{|c|c|c|c|c|}
		\hline
		1 & 1 & 1 & 1 & 2\\ 
		\hline
		2 & 3 & 3 & 3 & 3 \\
		\hline
		3 & 4 & 4 & 4 & 4 \\
		\hline
		4 & 5 & 5 & 5 & 6 \\
		\hline
		5 & 6\\
		\cline{1-2}
		6\\
		\cline{1-1}
	\end{array}
	\]
		\begin{enumerate}[label=(b), ref=b]
			\item\label{i:ssytptint} Fix $i$, $1\leq i\leq k$. For any  $r>p$,   the $r$-depth of sequence of $\ssyt[i]$ equals the $r$-depth sequence of $\ssyttint[i]$.
		\end{enumerate}
		\begin{enumerate}[label=(c), ref=c]
			\item
				\label{i:depth}   Fix $i$, $1\leq i\leq k$.   Let $a_1<\ldots<a_{\tint}$ be the the first $\tint$ entries arranged in increasing order of~$\procw[i]$. 
		Let $y_1<\ldots<y_{\tint+1}$ be the $(\tint+1)$-depth sequence of~$\ssyt[i]$.   Let $s$, $1\leq s\leq \tint+1$, be such that $a_1<\ldots<a_{s-1}<y_s<a_s<\ldots<a_\tint$ is the arrangement in increasing order of the first $\tint+1$ elements of~$\procw[i]$ (see Proposition~\ref{p:procnew} and the sentence preceding it).
				Then $y_1=a_1$, \ldots, $y_{s-1}=a_{s-1}$,  and $y_s$ occurs in row $s$ of~$\ssyttint[i]$ and in a column weakly to the right of that in which $\boxbtintpone$ occurs ($\boxbtintpone$ is the right most box in row $\tint+1$ of~$\ssyt$).
\end{enumerate}
\mysubsubsection{Base case of the induction}\mylabel{sss:base}
The assertions are easily verified in case $\tint=1$.  Indeed, for every $i$, $1\leq i\leq k$,  $\procw[i]^{(1)}$
has only one significant element and its first element is the entry in the first row in column~$i$ of~$\ssyt$.  This proves~\eqref{i:wseq}.    
Assertion~\eqref{i:ssytptint} is immediate since $\ssyt_1=\ssyt$.
Assertion~\eqref{i:depth} is vacuous in case~$s=1$.  In case $s=2$,  we have $a_1<y_2$,  where $a_1$ is the entry in the first row and column~$i$ of~$\ssyt[i]$ and $y_2$ is the right most entry in row~$2$ of~$\ssyt[i]$.  It follows that the box in row~$1$ and column~$i$ has $2$-depth~$2$,  so $a_1\leq y_1$.   Since $a_1$ is the largest entry in the first row  and $y_1$ occurs as an entry in the first row,  it follows that $y_1\leq a_1$.   Thus $y_1=a_1$.
\mysubsubsection{Proof of assertion~\eqref{i:wseq}}\mylabel{sss:pfwseq}
To simplify notation, write $g$ and $h$ for $\procw[i]$ and $\procw[i+1]$ respectively.    We need to prove that $\vecg^j\leq\vech^j$ for all $j\leq \tint$ (Lemma~\ref{l:bb}).   By the induction hypothesis,  we know this to be true for $j<\tint$,  so it remains to be proved only for $j=\tint$.    Let us write $a_1<\ldots<a_{\tint-1}$ for $\vecg^{p-1}$ and $b_1<\ldots<b_{\tint-1}$ for $\vech^{p-1}$.



Let $e_1<\ldots<e_p$ and  $f_1<\ldots<f_p$ be the $p$-depth sequences of $\ssyt[i]$ and $\ssyt[i+1]$ respectively. We have, evidently, $e_j\leq f_j$.  Let $s$ and~$t$, $1\leq s,t\leq p$, be such that
\[a_1<\ldots<a_{s-1}<e_s<a_{s}<\ldots<a_{\tint-1}\quad\textup{and}\quad b_1<\ldots<b_{t-1}<f_t<b_{t}<\ldots<b_{\tint-1} \]
are the sequences $\vecg^p$ and $\vech^p$. 

In the case $s\leq t$,\footnote{The case $s<t$ never actually occurs,  but that does not concern us here.}   the desired conclusion $\vecg^\tint\leq\vech^\tint$ follows rather easily from $\vecg^{\tint-1}\leq\vech^{\tint-1}$.  Indeed we have, in the case $s<t$:
\[
	\begin{array}{ll}
		g^\tint_j=a_j\leq b_j=h^\tint_j   &\textup {for $1\leq j\leq s-1$}\\
g^\tint_s=e_s<a_s\leq b_s=h^\tint_s \\
		g^\tint_j=a_{j-1}\leq b_{j-1}=h^\tint_{j-1}<h^\tint_j &\textup{for $s+1\leq j< t$}\\
		g^\tint_t=a_{t-1}\leq b_{t-1}<f_t=h^\tint_t \\ 
		g^\tint_j=a_{j-1}\leq b_{j-1}=h^\tint_j&\textup{for $t<j\leq p$}
	\end{array}
	\]
	\noindent For $s=t$,  the three middle lines in the display above should be replaced by $g^\tint_s=e_s\leq f_s=h^\tint_s$.

So let us assume $s>t$.   The cases $1\leq j<t$ and $s<j\leq p$,    are similar respectively to the first and last cases above.    It is for $j$ in the range $t\leq j\leq s$ that we need some care.      The key observation here is that $a_j=e_j$ for $j<s$.    
This follows from assertion~\eqref{i:depth} with $\tint$ replaced by $\tint-1$ (which we may assume to be true by induction).
Indeed, using this,  we are done as follows:
\[
	\begin{array}{ll}
		g^\tint_j=a_j =e_j\leq f_j\leq b_{j-1}=h^\tint_j   &\textup {for $t< j\leq s$}\\
		g^\tint_t=a_{t-1}=e_t<f_t=h^\tint_t \hfill \Box\\ 
	\end{array}
	\]
\mysubsubsection{Proof of assertion~\eqref{i:ssytptint}}\mylabel{sss:pffour} Fix $r>\tint$.   
By induction,  we know the statement for $\tint-1$ in place~$\tint$,  so the $r$-depth sequences of $\ssyt[i]$ and $\ssyttintmone[i]$ are the same.   It is therefore enough to prove that  the $r$-depth sequences of~$\ssyttintmone[i]$ and $\ssyttint[i]$ are the same.  
It is convenient to omit the ``$[i]$'' and just write $\ssyttintmone$ and $\ssyttint$ for $\ssyttintmone[i]$ and $\ssyttint[i]$ respectively.  
Assertion~\eqref{i:ssytptint} follows immediately from Corollaries~\ref{c:p:oldbox} and~\ref{c:p:2:crucial} below.

We denote by $\colmax$ the column number in which the right most box~$\boxbtint$ in row~$\tint$ of~$\ssyttintmone$ (equivalently~$\ssyt$) occurs.
\ignore{ 
A box in row $j$, $1\leq j\leq \tint$, of~$\ssyttintmone$ is said to be of {\em maximal depth\/} if it has $\tint$-depth $\tint-j$:  note that any box in row $j$ has $\tint$-depth at most $\tint-j$,  which justifies the terminology.
\bprop\mylabel{p:maxdepth}  Recall that $\boxbtint$ denotes the right most box in row~$\tint$ of~$\ssyttintmone$ (equivalently, of~$\ssyt$).
\begin{enumerate}
		\item\label{i:omd} If a box has maximal depth, then so has any box on top of it in the same column, and also so has any box to the left of it in the same row.
	\item\label{i:omd}  The box $\boxbtint$ and the boxes vertically above it all have maximal depth.
\end{enumerate}
\eprop
} 
The entries in any column of~$\ssyttintmone$ are also entries in that same column of~$\ssyttint$.  For every box $\boxb$ of~$\ssyttintmone$,   we denote by $\boxb'$ the (unique) box of~$\ssyttint$ in the same column as $\boxb$ and having the same entry.  The association $\boxb\mapsto\boxb'$ is evidently one-to-one.
Either $\boxb'$ is in the same row as $\boxb$ or in the next lower row. 

We classify boxes of~$\ssyttint$ as follows:  
\begin{itemize}
	\item {\em Old\/} boxes are those that are in the image of the above map $\boxb\mapsto\boxb'$.  {\em New\/} boxes are those that are not old.
	\item An {\em unmoved\/} box is an old box~$\boxbp$ that is in the same row as its preimage~$\boxb$. We write $\boxb'=\boxb$ in this case.
		A {\em moved\/} box is an old one that is not unmoved, or, in other words,  an old one that is in a row one lower than its preimage.
\end{itemize}
A box~$\boxb$ of~$\ssyttintmone$ is {\em moved\/} or {\em unmoved\/} accordingly as its image $\boxb'$ in~$\ssyttint$ is so.

As an illustration,   shown on the left in the display below is $\ssyt_3$ and on the right 
is $\ssyt_4$ in a particular case: $\tint=3$ and $\colmax=2$.     The entries in all the new boxes are in bold and underlined; those in unmoved boxes are in red;   those in moved boxes are in blue.  The $4$-depth sequences for $\ssyt_3[3]$ and $\ssyt_3[4]$
respectively are $7$, $5$, $2$, $1$ and $7$, $5$, $4$, $2$;    those for $\ssyt_3[l]$ for $l\geq5$ are all $7$, $6$, $4$, $3$.
\ignore{ 
\begin{equation}\label{e:oldnewbox}
	\begin{array}{|c|c|c|c|c|c|c|c|}
		\hline
		1 & 1 & 1 & 1 & 2 & 4 & 6 & 7 \\
		\hline
		2 & 2 & 2 & 2 & 3 & 5 & 8 & 8 \\
		\hline
		3 & 3 & 3 & 5 & 6 & 8 & 9 & 9 \\
		\hline
		5 & 6 & 7\\
		\cline{1-3}
		6 & 7 \\
		\cline{1-2}
		8 & 9\\
		\cline{1-2}
		9\\
		\cline{1-1}
	\end{array}
		\quad\quad\quad\quad
	\begin{array}{|c|c|c|c|c|c|c|c|}
		\hline
		1 & 1 & 1 & 1 & 2 & 4 & \nbem{4} & \nbem{4} \\
		\hline
		2 & 2 & 2 & 2 & 3 & 5 & 6 & 7 \\
		\hline
		3 & 3 & 3 & 5 & 6 & \nbem{6} & 8 & 8 \\
		\hline
		5 & 6 & 7 & \nbem{7} & \nbem{7} & 8 & 9 & 9 \\
		\hline
		6 & 7 \\
		\cline{1-2}
		8 & 9\\
		\cline{1-2}
		9\\
		\cline{1-1}
	\end{array}
\end{equation}
} 
\begin{equation}\label{e:oldnewbox}
	\begin{array}{|c|c|c|c|c|c|c|c|}
		\hline
		1 & 1 & 1 & 2 & 3 & 3 & 4 & 7 \\
		\hline
		2 & 2 & 2 & 4 & 6 & 7 & 8 & 8 \\
		\hline
		3 & 3 & 5 & 5 & 8 & 8 & 9 & 9 \\
		\hline
		6 & 7\\
		\cline{1-2}
		7 \\
		\cline{1-1}
		9\\
		\cline{1-1}
	\end{array}
		\quad\quad
	\begin{array}{|c|c|c|c|c|c|c|c|}
		\hline
		\textcolor{red}{1} & \textcolor{red}{1} & \textcolor{red}{1} & \textcolor{red}{2} & \textcolor{red}{3} & \textcolor{red}{3} & \nbem{3} & \nbem{3} \\
		\hline
		\textcolor{red}{2} & \textcolor{red}{2} & \textcolor{red}{2} & \textcolor{red}{4} & \nbem{4} & \nbem{4} & \textcolor{blue}{4} & \textcolor{blue}{7} \\
		\hline
		\textcolor{red}{3} & \textcolor{red}{3} & \textcolor{red}{5} & \textcolor{red}{5} & \textcolor{blue}{6} & \textcolor{blue}{7} & \textcolor{blue}{8} & \textcolor{blue}{8} \\
		\hline
		\textcolor{red}{6} & \textcolor{red}{7} & \nbem{7} & \nbem{7} & \textcolor{blue}{8} & \textcolor{blue}{8} & \textcolor{blue}{9} & \textcolor{blue}{9} \\
		\hline
		\textcolor{red}{7} \\
		\cline{1-1}
		\textcolor{red}{9}\\
		\cline{1-1}
	\end{array}
\end{equation}

\bprop\mylabel{p:boxtype} 
\begin{enumerate}
	\item\label{i:unmoved} In any column $i$ of $\ssyttintmone$ (respectively $\ssyttint$) with $i\leq\colmax$,  all boxes are unmoved (respectively old and unmoved).  In particular,  the right most box~$\boxbr$ in row~$r$ (with $r>\tint$) of~$\ssyttint$ is old and unmoved.
	\item\label{i:nb} Let $\newboxn$ be a new box.  Then, to the left of $\newboxn$ and in the same row,  in a column with number $i\geq\colmax$,  there is an old box carrying the same entry as $\newboxn$.\footnote{Any such box is actually unmoved,  but we don't need that bit of detail.}
\end{enumerate}
\eprop
\bmyproof Item~\eqref{i:unmoved} is clear.   Indeed $\ssyttintmone$ and $\ssyttint$ are identical in columns $i\leq\colmax$. 

To prove item~\eqref{i:nb}, suppose that $\newboxn$ occurs in column~$c$ of~$\ssyttint$.   
Then $\ssyttintmone[c]$ has $\tint-1$ boxes in its last column.   Let $a_1<\ldots<a_{\tint-1}$ be the entries in that column (top to bottom).   Let $y_1<\ldots<y_\tint$ be the $\tint$-depth sequence of~$\ssyt[c]$ (or, what amounts to the same by the induction hypothesis, of~$\ssyttintmone[c]$) and let $s$, $1\leq s\leq \tint$, be such that $a_1<\ldots<a_{s-1}<y_s<a_s<\ldots<a_{\tint-1}$ are the entries in the last column of~$\ssyttint[c]$.   The box with $y_s$ as its entry is~$\newboxn$, and $\newboxn$ occurs in row~$s$ of~$\ssyttint[c]$.  We may assume by induction that assertion~\eqref{i:depth} of~\S\ref{ss:pfmain} is true with $\tint-1$ in place of~$\tint$,  and conclude that
$y_s$ appears as an entry in row~$s$ of~$\ssyttintmone[c]$ in a column with number $i\geq\colmax$.
\emyproof
\bprop\mylabel{p:oldbox}    Let $\boxb_1$ and $\boxb_2$ be boxes of $\ssyttintmone$.    Then $\boxb_1\cgeq\boxb_2$ if and only if~$\boxb'_1\cgeq\boxb'_2$. 
\eprop
\bmyproof   Neither the entry nor the column number changes on passage from $\boxb$ to $\boxb'$.    While the row number could increase by at most~$1$ on this passage,   consider the facts that both $\ssyttintmone$ and $\ssyttint$ are SSYTs and that $\boxb_2$ (respectively $\boxb'_2$) is weakly to the East of~$\boxb_1$ (respectively $\boxb'_1$) and carries an entry which is strictly less. Together these imply that $\boxb_2$ (respectively $\boxb'_2$) occurs in a higher row than $\boxb_1$ (respectively $\boxb'_1$).
\emyproof
\bcor\mylabel{c:p:oldbox} Fix $r> p$.
Suppose that $\boxbr=\boxb_0\cgeq\ldots\cgeq\boxbdelta$ is a chain of boxes in~$\ssyttintmone$.
Then $\boxbr=\boxb'_0\cgeq\ldots\cgeq\boxbpdelta$ is a chain of boxes in~$\ssyttint$.  In particular,   the $r$-depth sequence of $\ssyttintmone$ is term for term dominated by the $r$-depth sequence of~$\ssyttint$.
\ecor
\bmyproof  That we get a chain on passing from $\boxb$ to~$\boxb'$ is clear from Proposition~\ref{p:oldbox}.   That ${\boxbr}'=\boxbr$ follows from Proposition~\ref{p:boxtype}~\eqref{i:unmoved}.    It is clear from the description of the association $\boxb\mapsto\boxb'$ that the entries of $\boxbdelta$ and $\boxbpdelta$ are the same.
\emyproof
\bprop\mylabel{p:2:crucial}   Fix $r>p$.  Given a chain $\boxbr=\boxbt_0\cgeq\ldots\cgeq\boxbtdelta$ of boxes in~$\ssyttint$,   there exists a chain 
$\boxbr=\boxb_0\cgeq\ldots\cgeq\boxbdelta$ of boxes in~$\ssyttintmone$ with $\boxbdelta$ having the same entry as $\boxbtdelta$ and being weakly to the Northwest of it (meaning,  the row and column numbers of $\boxbdelta$ each is at most that of the corresponding number of~$\boxbtdelta$).
\eprop
\bmyproof   Proceed by induction on~$\delta$.    For $\delta=0$,  the statement is easily seen to be true since $\boxbr={\boxbr}'$ is old and unmoved (Proposition~\ref{p:boxtype}~\eqref{i:unmoved}).  Suppose that $\delta\geq1$.   

First suppose that $\boxbtdelta$ is an old box.  Let $\boxbdelta$ be the unique box in~$\ssyttintmone$ such that $\boxbpdelta=\boxbtdelta$.  Note that $\boxbdelta$ shares its entry and column number with~$\boxbtdelta$ and is weakly to the North of it.    By induction, choose $\boxbr=\boxb_0\cgeq\ldots\cgeq\boxb_{\delta-1}$ with $\boxb_{\delta-1}$ being weakly to the Northwest of~$\boxbt_{\delta-1}$ and having the same entry.  Since $\boxb_{\delta-1}$ is weakly to the West of~$\boxbdelta$ with a strictly larger entry,  it follows that it is on a strictly lower row, and so $\boxb_{\delta-1}\cgeq\boxbdelta$.

Now suppose that $\boxbtdelta$ is a new box.   Using Proposition~\ref{p:boxtype}~\eqref{i:nb},  replace it by an old and unmoved box having the same entry and being to the left in the same row.  Suppose that the new $\boxbtdelta$ is in column~$c$.  If any $\boxbt_j$ for $j<\delta$ has a column number higher than~$c$,   replace it by the one in the same row in column number $c$.   We now get a chain with $\boxbtdelta$ being old,  so we are reduced to the case settled in the previous paragraph.
\emyproof
\bcor\mylabel{c:p:2:crucial} Fix $r>p$. The $r$-depth sequence of $\ssyttintmone$ dominates term for term the $r$-depth sequence of~$\ssyttint$.
\ecor
\ignore{
Given a relation $\boxb\cgeq\boxb_1$ among boxes in $\ssyttintmone$,   we have $\boxb'\cgeq\boxb_1'$  (the box $\boxb_1'$ has to be in some row above that of $\boxb'$,  since its entry is strictly less than that of~$\boxb'$ and it lies strictly to the East of~$\boxb'$).   If $\boxb$ of~$\ssyttintmone$ is in some column with $\tint$ or more boxes,  then $\boxb=\boxb'$ (meaning that $\boxb'$ occurs in the same row and column as $\boxb$).   In particular,  this is true of~$\boxbr$, the right most box in row~$r$.  It now follows that the $r$-depth sequence of~$\ssyttint$ is term for term at least that of~$\ssyttintmone$.

To prove that the $r$-depth sequence of~$\ssyttintmone$ is term for term at least that of~$\ssyttint$,  which is the harder direction, we introduce a definition and prove a proposition.  A box of $\ssyttint$ that is not of the form~$\boxb'$ for some $\boxb$ of~$\ssyttintmone$ is called a {\em new box\/}.     Each column of~$\ssyttint$ has at most one new box.   A column of~$\ssyttint$ has a new box if and only if that column of~$\ssyttintmone$ has $\tint-1$ boxes.    Let $\colmax$ be the largest column number such that $\ssyttintmone$ (equivalently $\ssyt$) has at least $\tint$~boxes in that column.  In other words,  $\colmax$ is the number of the column in which the right most box~$\boxbtint$ in row~$\tint$ of~$\ssyttintmone$ (equivalently $\ssyt$) occurs.
\bprop\mylabel{p:newbox}
Let $\newboxn$ be a new box of~$\ssyttint$.    Then, to the left of~$\newboxn$ (not necessarily immediately), there is a box that is not new, has the same entry as~$\newboxn$,  and lies in a column  whose number is at least $\colmax$. 
\eprop

Let $y_k$ be the maximal entry in a box of~$\boxbr$-depth $r-k$ of~$\ssyttint$.   Let $\boxbr=\boxc_0\cgeq \boxc_1\cgeq\ldots\cgeq\boxc_k$  be a chain of boxes in~$\ssyttint$ with the entry in~$\boxc_k$ being $y_k$.  If $\boxc_k$ is a new box,  then,  thanks to the proposition,  we replace it by a box $\boxc_k'$ in the same row that is not new and has the same entry.  Now, by induction,  we start replacing the boxes.
}

\mysubsubsection{Proof of assertion~\eqref{i:depth}}\mylabel{sss:pfdepth}
By assertion~\eqref{i:ssytptint}, we may take $y_1<\ldots<y_{\tint+1}$ to be the $(\tint+1)$-depth sequence of~$\ssyttint[i]$.  
Fix $j<s$.     We would like to show that $y_j=a_j$.    In what follows,  we write just ``depth'' to mean ``$(\tint+1)$-depth''.
Recall that, by definition, $y_j$ is the maximal entry in a box of depth $\tint+1-j$;  and $a_j$ occurs as the entry in row $j$ and column~$i$ of~$\ssyttint[i]$.
Any box of depth $\tint+1-j$ occurs in row~$j$ or above,  and $a_j$ dominates all the entries in those rows.
Thus it is enough to show that the box in row $j$ and column~$i$ of~$\ssyttint[i]$ has depth $\tint+1-j$.
Further, since any box in row~$j$ has depth at most $\tint+1-j$,   it is enough to show that the depth of that box is at least $\tint+1-j$.
Further,  it is enough to show this for $j=s-1$, since it follows then for the other $j<s$ as well.


By definition, $y_s$ occurs as an entry in a box $\boxb$ of~$\ssyttint[i]$ of depth $\tint+1-s$.      Such a~$\boxb$ can only appear in row $s$ or above.    But since $a_{s-1}<y_s$,  it follows that $\boxb$ cannot occur in row $s-1$ or above.   So it appears in row $s$, and so $\boxb\cgeq\boxc$ where $\boxc$ is the box in row~$s-1$ and column~$i$ of~$\ssyttint[i]$,  which means that $\boxc$ has depth $\tint+2-s$.\hfill$\Box$
\mysubsection{Proof that $\imsl\leq \procw$}\mylabel{ss:pfvleqw} 
It follows from assertion~\eqref{i:wseq} that $\procw[1]^{(\tint)}\leq \ldots\leq\procw[k]^{(\tint)}$ is a standard lift of~$\ssyt^{(\tint)}$.   Since $\imsl^{(\tint)}=\imsl(\tint)$ is the initial element of the minimal lift of~$\ssyt^{(\tint)}$ (by Proposition~\ref{p:procbasic2}),  it follows that $\imsl^{(\tint)}\leq\procw[k]^{(\tint)} =\procw^{(\tint)}$. Since $\imsl=\imsl^{(\tint)}$ and $\procw=\procw^{(\tint)}$ for large $p$,  it follows that $\imsl\leq\procw$.~\hfill$\Box$
\mysubsection{A technical lemma (that is invoked in~\S\ref{ss:uleqv})}\mylabel{ss:tlemma}
\ignore{ 
\bmyproof This is clear by an induction on the number $c$ of columns of~$\ssyt$.    When $c=1$,  the assertion is evident.  Let $\ssyt'$ be obtained from $\ssyt$ by deleting its last column.   Then, by induction, the initial element $v'$ of the minimal standard lift of~$\ssyt'$ has only $r$ significant elements.  By definition, $v$ is the minimal element subject to two constraints, namely, (a)~$v'\leq v$ and (b)~the initial elements of $v$ are the same as those that appear in the last column of~$S$.   Both these constraints will still be satisfied if the elements of $v$ in positions $r+1$ and beyond are rearranged in increasing order (see Lemma~\ref{l:procbasic} above).  Thus $v$ has only $r$ significant elements.   \emyproof
\blemma\mylabel{l:ony} Let $\sigma_1\leq\ldots\leq\sigma_k$ be a standard lift of~$\ssyt$.     Let $c$ be the serial number (from the left) of a column of~$\ssyt$ that contains a box of $\boxbrpone$-depth $\delta$.   Let $y$ be the entry in that box.   Then,  among the first $r+1$ elements of $\sigma_c$ (in its one-line notation),  there are at least $\delta+1$ that are at least $y$.
\elemma
\bmyproof Proceed by induction on~$\delta$.  Suppose first that $\delta=0$.    The only box with $\boxbrpone$-depth $0$ is the box $\boxbrpone$ itself.  Since $\boxbrpone$ occurs on row $r+1$,  the conclusion is easily verified to be true.

Now suppose that $\delta\geq1$.     From the hypothesis (and the definition of $\boxbrpone$-depth),  there exists, for some $c'<c$,  a box in the column~$c'$ of $\ssyt$ with entry $y'>y$ and of $\boxbrpone$-depth $\delta-1$.      By the induction hypothesis,   there exist,  among the first $r+1$ elements of $\sigma_{c'}$,  $\delta$ that are at least $y'$.   Since $\sigma_{c'}\leq\sigma_c$,   the same assertion holds with $\sigma_{c'}$ replaced by $\sigma_c$.     Now, $y$ too occurs in the first $r+1$ elements of $\sigma_c$.    Thus there are at least $\delta+1$ among the first $r+1$ elements of $\sigma_c$ that are at least $y$.  \emyproof
\bcor\mylabel{c:ony}   Let $y_{r+1}>\ldots>y_1$ the sequence as in the procedure in~\S\ref{sss:ssyt-to-w} to determine $\procwrpone$.   Then,  for every $j$, $1\leq j\leq r+1$,   among the first $r+1$ elements of $\imsl$,   there occur at least $r+2-j$ elements that are at least $y_j$.
\ecor
\bmyproof   By definition, $y_j$ occurs as an entry in some box of~$\ssyt$ of $\boxbrpone$-depth $r+1-j$.    Suppose $c$ is the column number in which such a box appears.   Choose the standard lift in the lemma above to be the minimal one.  Then,  by the lemma,    among the first $r+1$ elements of $\sigma_c$, there occur at least $r+2-j$ that are at least $y_j$.   Since $\sigma_c\leq v$,  the same assertion holds with $v$ in place of~$\sigma_c$.
\emyproof
} 
\blemma\mylabel{l:ony} Let $\sigma_1\leq\ldots\leq\sigma_k$ be a standard lift of~$\ssyt$.   Consider any box of $\boxbtint$-depth $\delta$ in $\ssyt$, for some positive integer $p$.  (Recall that $\boxbtint$ denotes the right most box in row $\tint$ of~$\ssyt$.)  Let $y$ be the entry in that box and $c$ be the serial number (from the left) of the column in which that box appears.  Then,  among the first $p$ elements of $\sigma_c$ (in its one-line notation),  there are at least $\delta+1$ that are at least $y$.
\elemma
\bmyproof Proceed by induction on~$\delta$.  Suppose first that $\delta=0$.    The only box with $\boxbtint$-depth~$0$ is the box $\boxbtint$ itself.  Since $\boxbtint$ occurs on row $\tint$,  the conclusion is easily verified to be true.

Now suppose that $\delta\geq1$.     From the hypothesis (and the definition of $\boxbtint$-depth),  there exists, for some $c'<c$,  a box in the column~$c'$ of $\ssyt$ with entry $y'>y$ and of $\boxbtint$-depth $\delta-1$.      By the induction hypothesis,   there exist,  among the first $\tint$ elements of $\sigma_{c'}$,  $\delta$ that are at least $y'$.   Since $\sigma_{c'}\leq\sigma_c$,   the same assertion holds with $\sigma_{c'}$ replaced by $\sigma_c$.     Now, $y$ too occurs in the first $\tint$ elements of $\sigma_c$.    Thus there are at least $\delta+1$ among the first $\tint$ elements of $\sigma_c$ that are at least $y$.  \emyproof
\bcor\mylabel{c:ony}   Let $y_1<\ldots<y_{\tint}$ be the $\tint$-depth sequence of~$\ssyt$ (this was used in the procedure in~\S\ref{sss:ssyt-to-w} to determine~$\procwtint$).   Then,  for every $j$, $1\leq j\leq \tint$,   among the first $\tint$ elements of $\imsl$,   there occur at least $p+1-j$ elements that are at least $y_j$.
\ecor
\bmyproof   By definition, $y_j$ occurs as an entry in some box of~$\ssyt$ of $\boxbtint$-depth $\tint-j$.    Suppose $c$ is the column number in which such a box appears.   Choose the standard lift in the lemma above to be the minimal one.  Then,  by the lemma,    among the first $\tint$ elements of $\sigma_c$, there occur at least $\tint+1-j$ that are at least $y_j$.   Since $\sigma_c\leq v$,  the same assertion holds with $v$ in place of~$\sigma_c$.
\emyproof
\mysubsection{Proof that $\procw\leq\imsl$}\mylabel{ss:uleqv} For $\tint$ a positive integer,  we prove, by induction on~$\tint$, that $\procw^{(\tint)}\leq \imsl^{(\tint)}$.    Since $\procw^{(\tint)}=\procw$ and $\imsl^{(\tint)}=\imsl$ for large $\tint$,  it will follow that $\procw\leq\imsl$.   First consider the case $\tint=1$. 
Let the right most entry in the first row of~$\ssyt$ be~$a$.  From the description of the procedure to produce $\procw$ in~\S\ref{sss:ssyt-to-w}, it is clear that $\procw_1=a$.  On the other hand, evidently,  the initial element of any standard lift of~$\ssyt$ has $a$ as its first element (in its one-line notation), so in particular $\imsl_1=a$.  This proves $\procw^{(1)}=\imsl^{(1)}$.

Now let $p>1$.  
%
By the induction hypothesis,  we have $\procw^{(\tint-1)}\leq\imsl^{(\tint-1)}$.
It is enough therefore to prove that $\procw_{\tint}\leq v_{\tint}$.  

Since we have proved that $\imsl\leq\procw$ (\S\ref{ss:pfvleqw}),  it follows that that $\imsl^{(\tint-1)}= \procw^{(\tint-1)}$.
Let $y_1<\ldots<y_p$ be the $p$-depth sequence of~$\ssyt$ and let $j$, $1\leq j\leq \tint$, be such that $\procw_\tint=y_j$.  Then there are exactly $\tint-j$ elements among $\procw_1$, \ldots, $\procw_{\tint-1}$  that are at least $y_j$.  
Corollary~\ref{c:ony} guarantees that
	among the first $\tint$ elements of $\imsl$, there are at least $\tint+1-j$ that are at least~$y_j$.   
Since $\procw_j=\imsl_j$ for $j\leq\tint-1$, it follows that $\procw_\tint=y_j\leq\imsl_\tint$, and we are done.~\hfill$\Box$

\begin{appendix}
\mysection{Multiple concatenations of LS paths}\mylabel{s:lsmult}
\noindent
The immediate provocation for this appendix comes from the need to quote its results (Propositions~\ref{p:dompatheqn} and~\ref{p:cryiso}) in the proof of the tableau decomposition rule for KK modules (\S\ref{ss:pfkktab}).  
These results are part of folklore. They are already hinted at by Littelmann in \cite{litt:inv}: see the ``precise combinatorial criterion'' alluded to in the paragraph preceding the theorem in~\S8.1 of that paper.  They are also later stated in~\cite[\S11]{litt_plactic} with a sketch of proofs.   However, we could not find a suitable reference with complete proofs.  This appendix aims to provide precisely such a
reference, presupposing knowledge of
		(a)~Littelmann's basic definitions and results on paths as in~\cite{litt:ann}
		and
		(b)~the results recalled and proved from scratch in~\S\ref{s:extremal} above.  

\ignore{
For $\theta$ a path in~$\paths$, write $\theta=\pi_1\concat\cdots\concat\pi_n$.  For $1\leq j\leq n$,  let $\tau_j$, $\varphi_j$ be respectively the final direction and initial direction of~$\pi_j$,  and let $\taut_j$, $\phit_j$ be arbitrary lifts in~$W$ of $\tau_j$, $\varphi_j$.  For $1\leq j<n$,  define
\begin{equation}\label{e:weylj}
	\weylj(\theta):=\brmin{W_{\lambda_j}\intvl{\taut_j^{-1}}\phit_{j+1}W_{\lambda_{j+1}}}
\end{equation}
\bprop\mylabel{p:invn}
The tuple $(\weyl_{1,2}(\theta),\ldots,\weyl_{n-1,n}(\theta))$ remains the same as $\theta$ varies over its equivalence class in~$\paths$ with respect to~$\sim$.
\eprop
\bmyproof   Lemma~\ref{l:wind}~\eqref{i:epipi'} generalizes to $n$-fold concatenations of LS paths (with the same proof):  if $e_\alpha(\theta)$ does not vanish,  then there exists $j_0$, $1\leq j_0\leq n$, such that  $e_\alpha(\theta)=\pi'_1\concat\cdots\concat\pi'_n$,  where $\pi'_j=\pi_j$ for $j\neq j_0$, $1\leq j\leq n$,  and $\pi'_{j_0}=e_\alpha\pi_{j_0}$.     Suppose that $e_\alpha(\theta)$ does not vanish and fix $j_0$ as above.    It follows easily from the definition of $e_\alpha$ that $e_\alpha(\pi_{j_0-1}\concat \pi_{j_0})=\pi_{j-1}\concat e_\alpha(\pi_{j_0})$ and $e_\alpha(\pi_{j_0}\concat\pi_{j_0+1})=e_\alpha\pi_{j_0}\concat\pi_{j_0+1}$.   The result now follows from Lemma~\ref{l:wind}~\eqref{i:winv}.
\emyproof
}
\mysubsection{Standard concatenations}\mylabel{ss:etastd}   Let $\lieg$ be symmetrizable Kac-Moody algebra.
Let $\lambda_1$, \ldots, $\lambda_n$ be dominant integral weights. 
For $j$, $1\leq j\leq n$,  let $\paths_{j}$ denote the set of LS paths of shape~$\lambda_j$.
Consider the set $\paths:=\pmconcat:=\{\pi_1\concat\cdots\concat\pi_n\st \pi_j\in\paths_{j}\ \textup{for}\ 1\leq j\leq n\}$ of paths.   For paths $\theta$ and $\theta'$ in~$\paths$, let us write $\theta\sim\theta'$ if either $e_\alpha\theta=\theta'$ or $f_\alpha\theta=\theta'$ for some simple root~$\alpha$.  This is a symmetric relation.   Let us continue to denote by $\sim$ the reflexive and transitive closure of this relation on~$\paths$.   

\mysubsubsection{The path $\eta(\theta)$}\mylabel{sss:appeta}
Fix a $\theta=\pi_1\concat\cdots\concat\pi_n$ in~$\paths$.   As in Proposition~\ref{p:dompath}, which is the special case $n=2$ of the present set up, it follows that:
\begin{itemize}
	\item  In the equivalence class of~$\paths$ containing~$\theta$,   there exists a unique path $\eta(\theta)$ that is killed by $e_\alpha$ for every simple root~$\alpha$.
	\item The $\eta(\theta)$ as above lies entirely in the dominant chamber.
\end{itemize}
\mysubsubsection{Standard concatenations}\mylabel{sss:appstd}
We want to characterize those $\theta$ for which $\eta(\theta)=\pi_{\lambda_1}\concat\cdots\concat\pi_{\lambda_n}$,  where as usual $\pi_{\lambda_j}$ denotes the straight line path from the origin to~$\lambda_j$.
Towards this,  put $W_j:=W_{\lambda_j}$, the stabiliser of~$\lambda_j$ in the Weyl group~$W$, and let $\tau_{1,j}>\ldots>\tau_{r_j,j}$ be the chain of elements in~$W/W_{j}$ forming the LS path~$\pi_j$ (for $1\leq j\leq n$).    Consider the tuple
\begin{equation}\label{e:pseq}
	\left(\tau_{1,1},\ldots\tau_{r_1,1},
\quad\ldots\quad,\tau_{1,j},\ldots,\tau_{r_j,j},\quad\ldots\quad,
	\tau_{1,n},\ldots,\tau_{r_n,n}\right)
\end{equation}
	which is an element of
	\begin{equation}\label{e:parseq}
		(W/W_{1})^{\times\ \textup{$r_1$ times}} \times \cdots \times
		(W/W_{j})^{\times\ \textup{$r_j$ times}} \times \cdots \times
		(W/W_{n})^{\times\ \textup{$r_n$ times}}
	\end{equation}
		We call the path $\theta$ {\em standard\/} if the tuple~\eqref{e:pseq} is standard in the sense of~\S\ref{ss:standard}. 
		A standard lift (respectively, a minimal standard lift) in the sense of~\S\ref{ss:standard} of the tuple~\eqref{e:pseq} is called a {\em standard lift\/} (respectively, {\em minimal standard lift\/}) of~$\theta$. 
		We denote by $\weyl(\theta)$ the initial element of the minimal standard lift of~$\theta$.

		We denote by $\pathstd$ the subset of~$\paths$ consisting of standard paths.
\bexamplenobox   The path $\pi_{\lambda_1}\concat\cdots\concat\pi_{\lambda_n}$ is standard, for $(\textup{identity},\ldots,\textup{identity})$ is its minimal standard lift.     Moreover, it is the only standard path in~$\paths$ with identity as the initial element of its minimal standard lift.  Thus:
\begin{equation} \label{e:minstdliftid}
	\{\theta \in \pathstd\st \weyl(\theta)\leq \textup{identity}\} =\{\pi_{\lambda_1}\concat\cdots\concat\pi_{\lambda_n}\}
\end{equation}
\eexamplenobox
Here is a characterization of the paths~$\theta$ in~$\paths$ for which $\eta(\theta)=\pi_{\lambda_1}\concat\cdots\concat\pi_{\lambda_n}$:
\bprop\mylabel{p:dompatheqn} \textup{(see \cite[\S8.1]{litt:inv})}  
$\eta(\theta)=\pi_{\lambda_1}\concat\cdots\concat\pi_{\lambda_n}$ if and only if $\theta$~is standard.
\eprop
\noindent  The proof of this proposition is given in~\S\ref{ss:pfdompatheqn}.
\ignore{
\bexamplenobox  Let $\lieg$ be finite dimensional of type $A_\ell$ and let the fundamental weights $\varpi_1$, \ldots, $\varpi_\ell$ be as numbered as in~\cite[]{bourbaki}.    Let $\lambda_1$, \ldots, $\lambda_n$ be:
\[
	\underbrace{\varpi_1,\ldots,\varpi_1}_{\textup{$m_1$ times}},\ldots,
	\underbrace{\varpi_j,\ldots,\varpi_j}_{\textup{$m_j$ times}},\ldots,
	\underbrace{\varpi_\ell,\ldots,\varpi_\ell}_{\textup{$m_\ell$ times}}
	\]
\noindent
The elements of $W/W_{\varpi_j}$ are parametrized by $j$ element subsets $\{1,\ldots,\ell+1\}$.  Each such subset is written as $\{1\leq i_1<\ldots<i_j\leq \ell+1\}$.   Given two such subsets $\underline{i}=\{1\leq i_1<\ldots<i_j\leq \ell+1\}$ and $\underline{i}'=\{1\leq i'_1<\ldots<i'_j\leq \ell+1\}$,  we have $\underline{i}\leq\underline{i}'$ in the Bruhat order on $W/W_{\varpi_j}$ if and only if $i_1\leq i'_1$, \ldots, $i_{j-1}\leq i'_{j-1}$, and $i_j\leq i'_j$.   For a permutation $\sigma$ of~$\ell+1$ elements whose one line notation is~$\sigma_1\ldots\sigma_{\ell+1}$,  the coset $\sigma W_{\varpi_j}$  corresponds to $\{1\leq i_1<\ldots<i_j\leq\ell+1\}$,  where $i_1$, \ldots, $i_j$ are the elements $\sigma_1$, \ldots, $\sigma_j$ arranged in increasing order.

The paths of shape $\varpi_j$ are all straight lines,  so they too are parametrized by elements of $W/W_{\varpi_j}$.   Thus a path in~$\paths$ can be represented by a ``tableau'',  where a {\em tableau\/} consists of $m_\ell+\ldots+m_1$  top-justified columns of boxes,  with each of the first $m_\ell$ columns (from the left) having $\ell$ boxes, each of the next $m_{\ell-1}$ columns having $\ell-1$ boxes, and so on;   the boxes are filled with numbers between $1$ and $\ell+1$,  the entries in each column being strictly decreasing downwards.

Let, for example, $d=5$, and $\mu=5+3+3+2$. Then $m_1=2$, $m_2=0$, $m_3=1$, $m_4=2$, $m_5=0$.   Then $\lambda_1$, \ldots, $\lambda_n$ equals $\varpi_1$, $\varpi_1$, $\varpi_3$, $\varpi_4$, $\varpi_4$.    The paths in~$\paths$ can then be identified with tableaux consisting of $5$ columns of boxes,  the first two columns having $4$ boxes each, the next column having $3$ boxes, and the last two columns having $1$ box each.\footnote{The reversal of order, which is admittedly annoying, is necessary to preserve entrenched conventions.}
Here are two examples of such tableaux:
\begin{equation}\label{e:tab}
	\begin{array}{|c|c|c|c|c|}
		\hline
		1 & 1 & 2 & 1 & 4\\
		\hline
		3 & 2 & 4 \\
		\cline{1-3}
		4 & 3 & 5\\
		\cline{1-3}
		5 & 5\\
		\cline{1-2}
		\end{array}
		\quad\quad \quad\quad \quad\quad \quad\quad
	\begin{array}{|c|c|c|c|c|}
		\hline
		1 & 1 & 2 & 3 & 4\\
		\hline
		2 & 3 & 4 \\
		\cline{1-3}
		3 & 4 & 5\\
		\cline{1-3}
		5 & 5\\
		\cline{1-2}
		\end{array}
\end{equation}
For $\theta=(\tau_1,\ldots,\tau_5)$ in~$\paths$,  the entries in the first column of the corresponding tableau define $\tau_5$ (which is path of shape $\lambda_5=\varpi_4$),  the entries in the second column define $\tau_4$, and so on, so that the entries in the last column define $\tau_1$.

A tableau as in~\eqref{e:tab} is called {\em standard\/} if the entries in every row are weakly increasing left to right.   The tableau on the right in~\eqref{e:tab} is standard but not the one on the left.    This terminology is justified because:
\bprop\label{e:stdjust}
A tableau is standard if and only if the corresponding path~$\theta$ is standard as in the sense defined earlier in this section.
\eprop
\bmyproof  Let $j$ be such that $1\leq j\leq \ell$,  and let $\underline{h}=\{1\leq h_1<\ldots<h_j\leq\ell+1\}$ be an element in~$W/W_{\varpi_j}$.   Denote by $\tilde{\underline{h}}$ the permutation whose one line notation is $h_1\ldots h_j h'_1\ldots h'_{\ell+1-j}$,  where $h'_1$, \ldots, $h'_{\ell+1-j}$ are the elements of $\{1,\ldots,\ell+1\}\setminus\{h_1,\ldots, h_j\}$ arranged in decreasing order.   Then $\tilde{\underline{h}}W_{\varpi_j}=\underline{h}$,  and $\tilde{\underline{h}}\leq \tilde{\underline{i}}$ where $\tilde{\underline{i}}$ is defined analogously for $\underline{i}=\{1\leq i_1<\ldots<i_k\leq \ell+1\}$ is an element of $W/W_{\varpi_k}$ with $j\geq k$ and $h_1\leq i_1$, \ldots, $h_{k-1}\leq i_{k-1}$, and $h_k\leq i_k$.  This proves the ``only if'' part.

Let $j$, $k$ be integers such that $1\leq k\leq j\leq \ell$.   Let $\sigma$, $\tau$ be permutations of~$\{1,\ldots,\ell+1\}$ whose one line notations respectively are $\sigma_1\ldots\sigma_{\ell+1}$ and $\tau_1\ldots\tau_{\ell+1}$.
Then $\sigma W_{\varpi_j}=\{h_1<\ldots<h_j\}=\underline{h}$ where $h_1$, \ldots, $h_j$ are just $\sigma_1$, \ldots, $\sigma_j$ arranged in increasing order, and $\tau W_{\varpi_k}=\{i_1<\ldots<i_k\}=\underline{i}$ where $i_1$, \ldots, $i_k$ are just $\tau_1$, \ldots, $\tau_k$ arranged in increasing order. 
If $\sigma\leq\tau$,  then $h'_1\leq i_1$, \ldots, $h'_k\leq i_k$,  where $h'_1$, \ldots, $h'_k$ are $\sigma_1$, \ldots, $\sigma_k$ arranged in increasing order.   It follows that $h_1\leq i_1$, \ldots, $h_k\leq i_k$,  since evidently $h_1\leq h'_1$, \ldots, $h_k\leq h'_k$.  This proves the ``if'' part of the assertion.\emyproof

\noindent
A procedure to find the initial element of the minimal standard lift of a standard tableau is given in~\S\ref{ss:tableaux_kk}.
\eexamplenobox
}
\mysubsection{Specializing to a classical case: the case of the special linear Lie algebra}\mylabel{ss:slla}   Preserve the notation of the previous subsection and specialize to the situation of~\S\ref{s:tabkk}: an integer $d\geq2$ is fixed, $\lieg=\sld$, etc.   Let $\mu$ be a dominant integral weight, or, equivalently a partition with less than $d$ parts.   Write $\mu$ as $\mu_1\geq\mu_2\geq\ldots\geq\mu_{d-1}\geq0\geq\ldots$.   Let $\varpi_1=\epsilon_1$, $\varpi_2=\epsilon_1+\epsilon_2$, \ldots, $\varpi_{d-1}=\epsilon_1+\cdots+\epsilon_{d-1}$ be the fundamental weights.   Let $W_{\varpi_j}$, $1\leq j< d$, denote the stabiliser in~$W$ of~$\varpi_j$.

Put $m_1=\mu_1-\mu_2$, \ldots, $m_{d-2}=\mu_{d-2}-\mu_{d-1}$, $m_{d-1}=\mu_{d-1}-\mu_{d}=\mu_{d-1}$,  and $n=m_1+m_2+\cdots+m_{d-1}$  (note that $n=\mu_1$).
Let $\lambda_1$, \ldots, $\lambda_n$ be:
\[
	\underbrace{\varpi_1,\ldots,\varpi_1}_{\textup{$m_1$ times}},\ldots,
	\underbrace{\varpi_j,\ldots,\varpi_j}_{\textup{$m_j$ times}},\ldots,
	\underbrace{\varpi_\ell,\ldots,\varpi_\ell}_{\textup{$m_\ell$ times}}
	\]
\noindent
so that $\mu=\lambda_1+\cdots+\lambda_n$.

The elements of $W/W_{\varpi_j}$ are parametrized by subsets of cardinality $j$ of $[d]$. Each such subset is written as $\{1\leq i_1<\ldots<i_j\leq d\}$.   Given two such subsets $\underline{i}=\{1\leq i_1<\ldots<i_j\leq d \}$ and $\underline{i}'=\{1\leq i'_1<\ldots<i'_j\leq d\}$,  we have $\underline{i}\leq\underline{i}'$ in the Bruhat order on $W/W_{\varpi_j}$ if and only if $i_1\leq i'_1$, \ldots, $i_{j-1}\leq i'_{j-1}$, and $i_j\leq i'_j$.   For a permutation $\sigma$ of~$[d]$ whose one line notation is~$\sigma_1\ldots\sigma_{d}$,  the coset $\sigma W_{\varpi_j}$  corresponds to $\{1\leq i_1<\ldots<i_j\leq d\}$,  where $i_1$, \ldots, $i_j$ are the elements $\sigma_1$, \ldots, $\sigma_j$ arranged in increasing order.

For permutations $\sigma$ and $\tau$ of~$[d]$ with respective one-line notations $\sigma_1\ldots\sigma_d$ and $\tau_1\ldots\tau_d$,    we have $\sigma\leq\tau$ in the Bruhat order if and only if $\sigma W_{\varpi_j}\leq\tau W_{\varpi_j}$ for every $j$, $1\leq j<d$:  see, for example, \cite{bb}.

The LS paths of shape $\varpi_j$ are all straight lines,  so they too are parametrized by elements of $W/W_{\varpi_j}$.   Thus a path in~$\paths$ can be represented by a ``tableau'',  where a {\em tableau\/} consists of $m_{d-1}+\ldots+m_1$  top-justified columns of boxes,  where each of the first $m_{d-1}$ columns (from the left) has $d-1$ boxes, each of the next $m_{d-2}$ columns has $d-2$ boxes, and so on;   the boxes are filled with numbers between $1$ and $d$,  the entries in each column being strictly increasing downwards.\footnote{The reversal of order, which is admittedly annoying, is necessary to preserve entrenched conventions.}

Let, for example, $d=5$ and $\mu=6+3+3+2$.  Then $m_1=3$, $m_2=0$, $m_3=1$, $m_4=2$;  $n=6$, and $\lambda_1$, \ldots, $\lambda_6$ equals $\varpi_1$, $\varpi_1$, $\varpi_1$, $\varpi_3$, $\varpi_4$, $\varpi_4$.    And the paths in~$\paths$ can be identified with tableaux consisting of $6$ top-justified columns of boxes,  the first two columns having $4$ boxes each, the next column having $3$ boxes, and the last three columns having $1$ box each.
Here are two examples of such tableaux:
\begin{equation}\label{e:tab}
	\begin{array}{|c|c|c|c|c|c|}
		\hline
		1 & 1 & 2 & 1 & 4 &3\\
		\hline
		3 & 2 & 4 \\
		\cline{1-3}
		4 & 3 & 5\\
		\cline{1-3}
		5 & 5\\
		\cline{1-2}
		\end{array}
		\quad\quad \quad\quad \quad\quad \quad\quad
	\begin{array}{|c|c|c|c|c|c|}
		\hline
		1 & 1 & 2 & 3 & 4 &5\\
		\hline
		2 & 3 & 4 \\
		\cline{1-3}
		3 & 4 & 5\\
		\cline{1-3}
		5 & 5\\
		\cline{1-2}
		\end{array}
\end{equation}
For $\theta=\tau_1\concat\cdots\concat\tau_6$ in~$\paths$,  the entries in the first column of the corresponding tableau define $\tau_6$ (which is a path of shape $\lambda_5=\varpi_4$),  the entries in the second column define $\tau_5$, and so on, until the entries in the last column define $\tau_1$.

\bprop\mylabel{p:stdjust}
A path $\theta$ in~$\paths$ is standard as defined earlier in this section (\S\ref{sss:appstd}) if and only if the entries in the tableau corresponding to it are weakly increasing in every row from left to right.  
\eprop
\bmyproof   Proposition~\ref{p:deodhar:0}~\eqref{i:deodp03} is relevant here.   In particular,  we could use it prove the if part,   but instead we directly construct an explicit standard lift.
Let $j$ be such that $1\leq j< d$,  and let $\underline{h}=\{1\leq h_1<\ldots<h_j\leq d\}$ be an element in~$W/W_{\varpi_j}$.   Denote by $\tilde{\underline{h}}$ the permutation whose one line notation is $h_1\ldots h_j h'_1\ldots h'_{d-j}$,  where $h'_1$, \ldots, $h'_{d-j}$ are the elements of $[d]\setminus\{h_1,\ldots, h_j\}$ arranged in decreasing order.   Clearly,  $\tilde{\underline{h}}W_{\varpi_j}=\underline{h}$.   Let $k\leq j$ and let
$\underline{i}=\{1\leq i_1<\ldots<i_k\leq d\}$ be an element of $W/W_{\varpi_k}$ such that $h_1\leq i_1$, \ldots, $h_{k-1}\leq i_{k-1}$, and $h_k\leq i_k$.  Let $\tilde{\underline{i}}$ be defined from $\underline{i}$ (as $\tilde{\underline{h}}$ is from $\underline{h}$).     Then,  as is not to hard to see,
$\tilde{\underline{h}}\leq \tilde{\underline{i}}$.  This proves the if part.

Let $j$, $k$ be integers such that $1\leq k\leq j< d$.   Let $\sigma$, $\tau$ be permutations of~$[d]$ with respective one line notations $\sigma_1\ldots\sigma_{d}$ and $\tau_1\ldots\tau_{d}$.
Then $\sigma W_{\varpi_j}=\{h_1<\ldots<h_j\}=\underline{h}$ where $h_1$, \ldots, $h_j$ are just $\sigma_1$, \ldots, $\sigma_j$ arranged in increasing order, and $\tau W_{\varpi_k}=\{i_1<\ldots<i_k\}=\underline{i}$ where $i_1$, \ldots, $i_k$ are just $\tau_1$, \ldots, $\tau_k$ arranged in increasing order. 
Suppose that $\sigma\leq\tau$.  Then $h'_1\leq i_1$, \ldots, $h'_k\leq i_k$,  where $h'_1$, \ldots, $h'_k$ are $\sigma_1$, \ldots, $\sigma_k$ arranged in increasing order.   It follows that $h_1\leq i_1$, \ldots, $h_k\leq i_k$,  since evidently $h_1\leq h'_1$, \ldots, $h_k\leq h'_k$.  This proves the only if part of the assertion.\emyproof
%
\noindent
\bcor\mylabel{c:bijSPstd}   The set~$\pathstd$ of standard paths in~$\paths$ may be identified with the set $\ssytsetmd$ of SSYT of shape $\mu$ (in the sense of~\S\ref{ss:ssyt-and-perm}) with entries from~$[d]$.
\ecor
\noindent
The path represented by the tableau on the left in~\eqref{e:tab} is not standard whereas the one represented by the tableau on the right is standard:  the tableau on the left is not a SSYT whereas the tableau on the right is.



\noindent
\mysubsection{Proof of Proposition~\ref{p:dompatheqn}}\mylabel{ss:pfdompatheqn}
\noindent
Towards the proof,  we first prove a lemma.
		\blemma\mylabel{l:standard}    If $\theta$ is standard,  then so is every element in the equivalence class of~$\paths$ containing~$\theta$.  \elemma
\bmyproofnobox 
Let $\theta$ be standard and $\alpha$ be a simple root.   We will presently show that $f_\alpha\theta$ is standard in case it does not vanish.  The proof that $e_\alpha\theta$ is also standard, which we omit, is analogous.  This will suffice to prove the lemma.  Let us write $W/W_{\lambda_{i_1}}\times \cdots \times W/W_{\lambda_{i_m}}$ for the Cartesian product~\eqref{e:parseq}, and denote by $(\tau_1,\ldots,\tau_m)$ the tuple~\eqref{e:pseq}.  

Suppose that $f_\alpha\theta$ does not vanish.  From the definition of $f_\alpha$,  it follows that, by increasing $m$ and replacing $\tau_j$ by $\tau_j$, $\tau_j$ for some choices of $j$, $1\leq j\leq m$, as necessary,   we may assume that the tuple as in~\eqref{e:pseq} corresponding to~$f_\alpha\theta$ is $(\tau'_1,\ldots,\tau'_m)$ where for every $j$, $1\leq j\leq m$,  we have
\begin{equation}\label{e:tautp}
	\textup{$\tau'_j$ is either $\tau_j$ or $s_\alpha\tau_j$, depending upon certain conditions.}
\end{equation}
To exploit these conditions,  it is useful to introduce the following terminology.  Let $j$ be an integer, $1\leq j\leq m$.   We call $j$
\begin{equation}\label{e:jdefn}
	\begin{array} {cl}
		\textup{{\em changing}} &  \textup{if $\tau_j'=s_\alpha\tau_j\neq\tau_j$}\\
		\textup{{\em changeable (but not changing)}} & \textup{if $\tau_j'=\tau_j<s_\alpha\tau_j$}\\
		\textup{{\em resisting}} & \textup{if $s_\alpha\tau_j<\tau_j$}\\
		\textup{{\em flat}} & \textup{if $s_\alpha\tau_j=\tau_j$}\\
	\end{array}
\end{equation}

Using this terminology,  we record some simple observations (\eqref{e:jchanging}),~\eqref{e:jresisting1}, and~\eqref{e:jchangeable} below) that we need for the proof. All of these follow readily from the definition of~$f_\alpha$ as in~\cite{litt:ann}.   To begin with:
\begin{equation}\label{e:jchanging}
	\begin{array}{l}
		\textup{$j$ is changing only if $\tau_j<s_\alpha\tau_j$,}\\
		\textup{so the cases in~\eqref{e:jdefn} are exhaustive and mutually exclusive.}
	\end{array}
\end{equation}
In particular this means that $\tau'_j=\tau_j$ if $j$ is resisting.  So we have:
\begin{equation}\label{e:jresisting1}
	\textup{If $j$ is resisting or flat or changeable (but not changing), then $\tau_j'=\tau_j$.}
\end{equation}

We call $j$, $1\leq j\leq m$, {\em unobstructed\/} if there exists $k$, $j\leq k\leq m$, such that $k$ is changing and there does not exist $j'$ with $j'$ resisting and $j\leq j'< k$.    We call $j$ {\em obstructed\/} if it is not unobstructed.    Evidently:
\begin{equation}\label{e:jevident}\textup{$j$ is unobstructed if it is changing, and $j$ is obstructed if it is resisting.}\end{equation} 
	so, from~\eqref{e:jresisting1}:
	\begin{equation}\label{e:jobstructed}
		\textup{If $j$ is obstructed,  then $\tau'_j=\tau_j$.}
	\end{equation}
	We also have (from the definition of the operator~$f_\alpha$):
\begin{gather}
\label{e:jchangeable}
	\textup{If $j$ is changeable (but not changing), then $j$ is obstructed.}  
\end{gather}

Let now $\taut_1\geq\ldots\geq\taut_m$ be a standard lift of~$\theta$.   For $j$, $1\leq j\leq m$,  define $\taut_j'$ by:
\begin{equation}\label{e:taupt1}
	\taut_j' =\left\{
		\begin{array}{cl}
			s_\alpha\taut_j & \textup{if $j$ is changing}\\
			\taut_j & \textup{if $j$ is obstructed}\\
		\end{array}
		\right.
\end{equation}
and, when $j$ is flat and unobstructed,   by a downward induction as required:
\begin{equation}\label{e:taupt2}
	\taut_j':=\textup{the smaller of $\taut_j$ and $s_\alpha\taut_j$ that is larger than or equal to $\taut'_{j+1}$}
\end{equation}
Since $\taut'_{j+1}$  is either $\taut_{j+1}$ or $s_\alpha\taut_{j+1}$ (by downward induction),  it follows 
	(by an application of the basic observation~\eqref{e:parallel} in~\S\ref{ss:rsltsextrml} applied to the hypothesis that $\taut_j\geq\taut_{j+1}$)
that \begin{equation}\label{e:tautpstd} \taut_{j}\join s_\alpha\taut_j\geq\taut'_{j+1}\end{equation}
	so at least one of~$\taut_j$ and $s_\alpha\taut_j$ is larger than or equal to $\taut'_{j+1}$ and~\eqref{e:taupt2} makes sense.

We now argue that $\taut'_j\geq\taut'_{j+1}$ for all $1\leq j<m$.  If $j$ is flat and unobstructed,  then this follows from the definition~\eqref{e:taupt2} of $\taut'_j$.   If $j$ is either changing or resisting, then $\taut'_j=\taut_j\join s_\alpha\taut_j$ (from \eqref{e:jdefn}, \eqref{e:jchanging}, and~\eqref{e:jresisting1}),  so it follows from~\eqref{e:tautpstd}  that~$\taut'_j\geq\taut'_{j+1}$.  By the mutual exclusivity of the cases in~\eqref{e:jdefn} (which follows from~\eqref{e:jchanging} as already remarked) and \eqref{e:jchangeable}, we may assume that $j$ is obstructed but not resisting.    But then $j+1$ is also obstructed, and so $\taut'_{j+1}=\taut_{j+1}$ by~\eqref{e:taupt1},   and $\taut'_j=\taut_j\geq\taut_{j+1}=\taut'_{j+1}$.


We claim that $\taut'_1\geq \ldots\geq\taut'_m$ is a standard lift of~$f_\alpha\theta$.
It remains only to verify that $\taut_j' W_{\lambda_{i_j}}=\tau_j'$  for every $j$, $1\leq  j\leq m$.    This is easily done, as follows:
\begin{itemize}
	\item $j$ changing:  $\taut'_j=s_\alpha\taut_j$ by~\eqref{e:taupt1}, so $\taut'_jW_{\lambda_{i_j}}=s_\alpha\taut W_{\lambda_{i_j}}=s_\alpha\tau_j$.  But $s_\alpha\tau_j=\tau'_j$  by~\eqref{e:jdefn}.
	\item $j$ obstructed:  $\taut'_j=\taut_j$ by~\eqref{e:taupt1}, so $\taut'_jW_{\lambda_{i_j}}=\taut W_{\lambda_{i_j}}=\tau_j$.  But $\tau_j=\tau'_j$  by~\eqref{e:jobstructed}.
	\item $j$ flat:  $\taut'_j$ is either $\taut_j$ or $s_\alpha\taut_j$,   so $\taut'_j W_{\lambda_{i_j}}$ is either $\tau_j$ or $s_\alpha\tau_j$. But $\tau_j=s_\alpha\tau_j=\tau'_j$ by~\eqref{e:jresisting1}.\hfill$\Box$
\end{itemize}
\emyproofnobox
\pubpri{}{\noindent
For the record, here is yet another fact that follows readily from the definition of~$f_\alpha$ (to go along with~\eqref{e:jchanging}, \eqref{e:jresisting1}, and \eqref{e:jchangeable}),  although we have no need to refer to it.
\bprop\mylabel{p:faextra}
With notation and terminology as in the proof of Lemma~\ref{l:standard},   suppose that $1\leq j\leq \ell\leq m$,  with $j$ changing and $\ell$ resisting.  Then there exists $k$ with $j<k<\ell$ such that $k$ is changeable (but not changing). 
\eprop
}
\bmyproofof{Proposition~\ref{p:dompatheqn}}  Write $\eta$ for $\eta(\theta)$. If $\eta=\pi_{\lambda_1}\concat\cdots\concat\pi_{\lambda_n}$,  then $\eta$ is standard and so $\theta$ is standard by the previous lemma.  Now suppose that $\theta$ is standard.  Then so is~$\eta$ by the lemma.  Let us write $W/W_{\lambda_{i_1}}\times \cdots \times W/W_{\lambda_{i_m}}$ for the Cartesian product~\eqref{e:parseq}, and denote by $(\sigma_1,\ldots,\sigma_m)$ the minimal standard lift of $\eta$.   

Let $\alpha$ be any simple root.   We claim that there cannot exist~$k$, $1\leq k\leq m$, such that:
\begin{equation}\label{e:kcond}
	\textup{$s_\alpha\sigma_k <\sigma_k$}\quad\quad\textup{and}\quad\quad\textup{$s_\alpha\sigma_j W_{\lambda_{i_j}}=\sigma_j W_{\lambda_{i_j}}$ for all $1\leq j<k$}
\end{equation}
	To prove the claim, we suppose such a $k$ exists and arrive at a contradiction.   We have $s_\alpha\sigma_k W_{\lambda_{i_k}}\leq \sigma_k W_{\lambda_{i_k}}$.   If strict inequality holds here,  then $e_\alpha \eta$ does not vanish,  a contradiction,     so equality holds.      If $s_\alpha\sigma_{k+1}>\sigma_{k+1}$, then $s_\alpha\sigma_k=\sigma_k\meet s_\alpha\sigma_k>\sigma_{k+1}\meet s_\alpha\sigma_{k+1}=\sigma_{k+1}$,   a contradiction to the hypothesis that $(\sigma_1,\ldots,\sigma_m)$ is a minimal standard lift of~$\eta$ (because then $s_\alpha\sigma_k$ would work as a lift in place of $\sigma_k$).   Thus~\eqref{e:kcond} holds with $k$ replaced by $k+1$.   Repeating these arguments sufficiently many times,  we conclude that $s_\alpha\sigma_m<\sigma_m$ and $s_\alpha\sigma_j W_{\lambda_{i_j}}=\sigma_j W_{\lambda_{i_j}}$ for all $1\leq j\leq m$.    But then $\sigma_m$ is not the minimal element in the coset $\sigma_m W_{\lambda_{i_m}}$,  which contradicts the hypothesis that $(\sigma_1,\ldots,\sigma_m)$ is the minimal standard lift of~$\eta$.

	To show that $\eta=\pi_{\lambda_1}\concat\cdots\concat\pi_{\lambda_n}$,  it suffices to show that $\sigma_1$ is the identity element of the Weyl group~$W$.  If $\sigma_1$ is not the identity element, let $\alpha$ be a simple root such that $s_\alpha\sigma_1<\sigma_1$.   Then~\eqref{e:kcond} holds with $k=1$,   a contradiction.
\emyproofof
\mysubsection{The crystal isomorphism}\mylabel{ss:crystal}
Fix notation as in the beginning of~\S\ref{ss:etastd}.    
Let $\pathstd$ denote the set of all standard paths in~$\paths$.  
By Proposition~\ref{p:dompatheqn},  $\pathstd$ is precisely the set of paths $\theta$ in $\paths$ for which $\eta(\theta)=\pi_{\lambda_1}\concat\cdots\concat\pi_{\lambda_n}$.
Thus, by~\cite[Theorem~7.1]{litt:ann},  there is a (unique) crystal isomorphism\footnote{``Crystal isomorphism'' just means a bijection that commutes with the action of the root operators $f_\alpha$ and $e_\alpha$.} 
 $\cryiso:\pathsl\to\pathstd$,  where $\pathsl$ denotes the set of LS paths of shape $\lambda=\lambda_1+\cdots+\lambda_n$.
 \bprop\mylabel{p:cryiso0}
 The isomorphism~$\cryiso$ has the following properties:
 \begin{itemize}
	 \item 
 The straight line path~$\pi_\lambda$ (from the origin to $\lambda$) is mapped under $\cryiso$ to $\pi_{\lambda_1}\concat\cdots\concat\pi_{\lambda_n}$. 
		 \item The end point of~$\pi$ in $\pathsl$ is the same as that of its image~$\cryiso\pi$.
		 \item $\pi$ is $\lambda$-dominant if and only $\cryiso\pi$ is so.
 \end{itemize}
 \eprop
 \bmyproofnobox
 The first item is because~$\pi_\lambda$ (respectively $\pi_{\lambda_1}\concat\cdots\concat\pi_{\lambda_n}$) is the unique path in~$\pathsl$ (respectively $\pathstd$) on which $e_\alpha$ vanishes for every simple~$\alpha$.   The second is because (a)~$\pi_\lambda$ and $\pi_{\lambda_1}\concat\cdots\concat\pi_{\lambda_n}$ both have $\lambda$ as end point, (b)~every path in $\pathsl$ (respectively $\pathstd$) can be obtained by acting a sequence of~$f_\alpha$ operators on~$\pi_\lambda$ (respectively $\pi_{\lambda_1}\concat\cdots\concat\pi_{\lambda_n}$),  (c)~the first item, and finally (d)~if $f_\alpha$ does not vanish on any path~$\sigma$ (in either $\pathsl$ or $\pathstd$) then $f_\alpha\sigma(1)=\sigma(1)-\alpha$.  As for the third item,  we make two observations from which it follows that $\cryiso$ preserves $\lambda$-dominance:
 \begin{itemize}
	 \item a path $\sigma$ in $\pathsl\cup\pathstd$ is $\lambda$-dominant if and only if $e_\alpha(\pi_\lambda\concat\sigma)$ vanishes for all simple roots $\alpha$.   
		 \item For paths $\pi_1$ and $\pi_2$ in $\pathsl\cup\pathstd$,   $e_\alpha(\pi_1\concat\pi_2)$ equals either $e_\alpha\pi_1\concat\pi_2$ or $\pi_1\concat e_\alpha\pi_2$ depending precisely upon whether or not $i\geq j$ where $i$ (respectively~$j$) is the maximum non-negative integer~$k$ such that $f_\alpha^k \pi_1$ (respectively $e_\alpha^k\pi_2$) does not vanish.~\hfill$\Box$
 \end{itemize}
\emyproofnobox
 
\bprop\mylabel{p:cryiso}  For an LS path $\pi$ of shape $\lambda$,  the minimal element in the initial direction of~$\pi$ equals the initial element $\weyl(\cryiso{\pi})$ of the standard minimal lift of $\cryiso{\pi}$.
\eprop
\noindent
The proposition follows by combining Corollary~\ref{c:l:cryisoprelim:2} with Lemma~\ref{l:cryisoprelim2}.
\ignore{
\begin{equation}\label{e:taupt}
	\taut_j' =\left\{
		\begin{array}{ll}
			s_\alpha\taut_j = \taut_j\join s_\alpha\taut_j & \textup{if $\tau_j'=s_\alpha\tau_j$ and $s_\alpha\tau_j>\tau_j$}\\
			\taut_j & \textup{if $\tau_j'=\tau_j$ and $s_\alpha\tau_j\neq\tau_j$}\\
			\taut_j\join s_\alpha\taut_j & \textup{if $s_\alpha\tau_j=\tau_j$, and there exists $j'<j<j''$, with both $j'$ and $j''$ changing}\\
		\end{array}
		\right.
\end{equation}
\emyproof
}
\mysubsubsection{A useful observation (Corollary~\ref{c:l:cryisoprelim:2})}\mylabel{sss:useful}
Let $\mu$ be a dominant integral weight.  For a Weyl group valued function $\funf:\pathsm\to W$ on the set $\pathsm$ of LS paths of shape $\mu$,  and $v$ an element of~$W$,   put $\pathsmvf:=\{\pi\in\pathsm\st \funf(\pi)\leq v\}$.
\blemma\mylabel{l:cryisoprelim} 
Suppose that the following conditions hold for $\pi$ an arbitrary path in~$\pathsm$ and $w:=\funf(\pi)$:   
\begin{enumerate}
	\item \label{i:l:cryisoprelim:1}If $\alpha$ a simple root with $s_\alpha w<w$,  then $e_\alpha\pi$ does not vanish.
	\item \label{i:l:cryisoprelim:2} Suppose $f_\alpha\pi$ does not vanish.   Then either \textup{(a)}~$\funf(f_\alpha\pi)=w$  or \textup{(b)}~$\funf(f_\alpha\pi)=s_\alpha w>w$ and $e_\alpha\pi$ vanishes.
\end{enumerate}
Then, for $v$ in $W$ and $\beta$ simple such that $s_\beta v<v$:
\[\pathsmvf=\{f_\beta^k\pi\st\pi\in\paths_{\mu,s_\beta v}(\funf),\ k\geq0,\ \textup{$f_\beta^k\pi$ does not vanish}\}\]
\elemma
\bmyproof    
Since $\pathsmbvf\subseteq\pathsmvf$ and, by~\eqref{i:l:cryisoprelim:2},  $\pathsmvf$ is closed under the action of~$f_\alpha$, it follows that the right hand side is contained in~$\pathsmvf$.  
To prove the other containment,  let $\sigma$ be in~$\pathsmvf$.   Let $k\geq0$ be maximal such that $e_\beta^k\sigma$ does not vanish,  and put $\pi:=e_\beta^k\sigma$.    Then $f_\beta^k\pi=\sigma$,  so it is enough to show that $\pi$ is in~$\pathsmbvf$.

Put $w:=\funf(\pi)$.    On the one hand,  since $f_\beta^k\pi=\sigma$, it follows from~\eqref{i:l:cryisoprelim:2} that $\funf(\sigma)$ equals either $w$ or $s_\beta w$, so that $w\leq w\join s_\beta w =\funf(\sigma)\join s_\beta\funf(\sigma)\leq v\join s_\beta v=v$.   But, on the other,  if $s_\beta w<w$,  then $e_\beta\pi$ does not vanish by~\eqref{i:l:cryisoprelim:1},  a contradiction to the maximality of~$k$.    Thus we have $w<s_\beta w$ and $w=w\meet s_\beta w\leq v\meet s_\beta v=s_\beta v$.
\ignore{
To prove the other containment,  first observe that $\pathsmbvf$ is evidently contained in the right hand side.  Let $\sigma\in\pathsmvf$.  If
$\funf(\sigma)\leq s_\beta v$, 
then, by the observation,   $\sigma$ belongs to the right hand side.    So suppose that $\funf(\sigma)\not\leq s_\beta v$.   Then $s_\beta\funf(\sigma)<\funf(\sigma)$ (for,  if not,  then $\funf(\sigma)=\funf(\sigma)\meet s_\beta\funf(\sigma)\leq v\meet s_\beta v=s_\beta v$).

Put $w:=\funf(\sigma)$.     By~\eqref{i:l:cryisoprelim:1},  $e_\beta \sigma$ does not vanish.   Let $k$ be maximum such that $e_\beta^k\sigma$ does not vanish.  We have $f_\beta^k(e_\beta^k\sigma)=\sigma$.      By~\eqref{i:l:cryisoprelim:2},  $\funf(e_\beta^k\sigma)$ is either $w$ or $s_\beta w$.   But if it is~$w$,  then $e_\beta(e_\beta^k\sigma)$ does not vanish by~\eqref{i:l:cryisoprelim:1},  so the maximality of~$k$ is violated.    Thus $\funf(e_\beta^k\sigma)=s_\beta w$,  which means that $e_\beta^k\sigma$ belongs to~$\pathsmbvf$ and $\sigma=f_\beta^k(e_\beta^k\sigma)$ belongs to the right hand side.

To prove the other containment,  let $\pi$ be in~$\pathsmvf$.     If $s_\beta \funf(\pi)>\funf(\pi)$,   then $\funf(\pi)=s_\beta \funf(\pi)\meet \funf(\pi)\leq s_\beta v\meet v=s_\beta v$,  so $\pi$ belongs to $\pathsmbvf$,  which is evidently contained in the right hand side.  

So suppose that $s_\beta\funf(v)<\funf(v)$.   Then $e_\beta \pi$ does not vanish by~\eqref{i:l:cryisoprelim:1}.  Let $k$ be maximal such that $e_\beta^k\pi$ does not vanish.   Since $f_\beta^k(e_\beta^k\pi)=\pi$,   it follows from~\eqref{i:l:cryisoprelim:2} that $\funf(e_\beta^k\pi)$ is either $v$ or $s_\beta v$.   If it equals $v$,  then~\eqref{i:l:cryisoprelim:1} implies that~$e_\beta (e_\beta^k\pi)$ does not vanish,  a contradiction to the choice of~$k$.   Thus $\funf(e_\beta^k\pi)=s_\beta v$,  and so $\pi=f_\beta^k(e_\beta^k\pi)$ belongs to the right hand side.
}
\emyproof
\ignore{
\blemma\mylabel{l:cryisoprelim} 
Suppose that the following conditions hold for $\pi$ an arbitrary path in~$\pathsm$ and $w:=\funf(\pi)$:   
\begin{enumerate}
	\item \label{i:l:cryisoprelim:1}If $\alpha$ a simple root with $s_\alpha w<w$,  then $e_\alpha\pi$ does not vanish.
	\item \label{i:l:cryisoprelim:2} Suppose $f_\alpha\pi$ does not vanish.   Then either \textup{(a)}~$\funf(f_\alpha\pi)=w$  or \textup{(b)}~$\funf(f_\alpha\pi)=s_\alpha w>w$ and $e_\alpha\pi$ vanishes.
\end{enumerate}
Then, for $v$ in $W$ and $\beta$ simple such that $s_\beta v<v$:
\[\pathsmvf=\{f_\beta^k\pi\st\pi\in\paths_{\mu,s_\beta v}(\funf),\ \ k\geq0,\ \ \textup{$f_\beta^k\pi$ does not vanish}\}\]
\elemma
\bmyproof    For $\pi$ in $\pathsmbvf$ or in $\pathsmvf$,  we have $\funf(f_\beta \pi)\geq v$ by~\eqref{i:l:cryisoprelim:2}, so the right hand side is contained in~$\pathsmvf$.    To prove the other containment,  let $\pi$ be in~$\pathsmvf$.     If $s_\beta \funf(\pi)>\funf(\pi)$,   then $\funf(\pi)=s_\beta \funf(\pi)\meet \funf(\pi)\leq s_\beta v\meet v=s_\beta v$,  so $\pi$ belongs to $\pathsmbvf$,  which is evidently contained in the right hand side.  

So suppose that $s_\beta\funf(v)<\funf(v)$.   Then $e_\beta \pi$ does not vanish by~\eqref{i:l:cryisoprelim:1}.  Let $k$ be maximal such that $e_\beta^k\pi$ does not vanish.   Since $f_\beta^k(e_\beta^k\pi)=\pi$,   it follows from~\eqref{i:l:cryisoprelim:2} that $\funf(e_\beta^k\pi)$ is either $v$ or $s_\beta v$.   If it equals $v$,  then~\eqref{i:l:cryisoprelim:1} implies that~$e_\beta (e_\beta^k\pi)$ does not vanish,  a contradiction to the choice of~$k$.   Thus $\funf(e_\beta^k\pi)=s_\beta v$,  and so $\pi=f_\beta^k(e_\beta^k\pi)$ belongs to the right hand side.
\emyproof
}
\bcor\mylabel{c:l:cryisoprelim} Let $\iota:\pathsm\to W$ be the function that maps each path to the minimal element in its initial direction. 
Then, for $v$ in $W$ and $\beta$ simple such that $s_\beta v<v$:
\[\pathsmvi=\{f_\beta^k\pi\st\pi\in\paths_{\mu,s_\beta v}(\iota), k\geq0, \textup{$f_\beta^k\pi$ does not vanish}\}\]
\ecor
\bmyproof   It follows easily from the definition of the operators $e_\alpha$ and $f_\alpha$ that the hypothesis of the lemma are satisfied for the function~$\iota$.   (See also~\cite[Lemma in \S5.3]{litt:inv}.)
\emyproof
\bcor\mylabel{c:l:cryisoprelim:2}   Suppose in addition to the conditions~\eqref{i:l:cryisoprelim:1} and~\eqref{i:l:cryisoprelim:2} of~Lemma~\ref{l:cryisoprelim} the function~$\funf$ satisfies the following:  $\paths_{\mu,\textup{identity}}(\funf)=\{\pi_\lambda\}$. 
Then $\funf=\iota$.
\ecor
\bmyproof
We proceed by induction on $v$ to show that $\pathsmvf=\pathsmvi$.
This will suffice.  If $v=\textup{identity}$, then both sets are equal to $\{\pi_\lambda\}$ and the result holds.   So suppose that $v>\textup{identity}$.   Choose a simple root $\alpha$ such that $w:=s_\alpha v<v$.    By the induction hypothesis, $\pathsmwf=\pathsmwi$.
But, by Lemma~\ref{l:cryisoprelim},  we have \[\pathsmvf =\{f_\alpha^k\pi\st\pi\in\pathsmwf,\ k\geq0,\ f_\alpha^k\pi\textup{ does not vanish}\}\] and, by Corollary~\ref{c:l:cryisoprelim},  we have \[ \pathsmvi=\{f_\alpha^k\pi\st\pi\in\pathsmwi,\ k\geq0,\ f_\alpha^k\pi\textup{ does not vanish}\},\]
so it is clear that $\pathsmvf=\pathsmvi$.
	\emyproof

	\noindent
The above corollary together with the following lemma proves Proposition~\ref{p:cryiso}.  The proof of the lemma occupies~\S\ref{ss:pfcryiso}
	\blemma\mylabel{l:cryisoprelim2}  Fix notation as in the first paragraph of \S\ref{ss:crystal}.  Let $\funf:\pathsl\to W$ be the Weyl group valued function on~$\pathsl$ given by $\funf(\pi):=\weyl(\cryiso(\pi))$.   Then $\paths_{\lambda,\textup{identity}}(\funf)=\{\pi_\lambda\}$ and $\funf$ satisfies the conditions~\eqref{i:l:cryisoprelim:1} and~\eqref{i:l:cryisoprelim:2} of Lemma~\ref{l:cryisoprelim}.
	\elemma
\ignore{
\blemma\mylabel{l:cryisoprelim} Let $\theta$ be a path in~$\pathstd$ and let $w:=\weyl(\theta)$ be the initial element of its minimal standard lift.  Then
\begin{enumerate}
	\item If $\alpha$ a simple root with $s_\alpha w<w$,  then $e_\alpha\theta$ does not vanish.
	\item Suppose $f_\alpha\theta$ does not vanish.   Then either (a)~$\weyl(f_\alpha\theta)=w$  or (b)~$\weyl(f_\alpha\theta)=s_\alpha w>w$ and $e_\alpha\theta$ vanishes.
	\item Suppose $e_\alpha\theta$ does not vanish.   Then either (a)~$\weyl(e_\alpha\theta)=w$  or (b)~$\weyl(e_\alpha\theta)=s_\alpha w<w$ and $e_\alpha^2\theta$ vanishes.
\end{enumerate}
\elemma
}
\mysubsection{Completion of the proof of Proposition~\ref{p:cryiso}:  Proof of Lemma~\ref{l:cryisoprelim2}}\mylabel{ss:pfcryiso}
We first prove:
\ignore{
For $v$ an element of the Weyl group,  put $\pathslv:=\{\pi\in\pathsl\st \indir(\pi)\leq v\}$ and $\pathstdv:=\{\theta\in\pathstd\st\weyl(\theta)\leq v\}$.  We proceed by induction on $v$ to show that $\cryiso$ induces a crystal isomorphism between $\pathslv$ and $\pathstdv$.  This will suffice.

If $v=\textup{identity}$,  then $\pathslv=\{\pi_\lambda\}$ and $\pathstdv=\{\pi_{\lambda_1}\concat\cdots\concat\pi_{\lambda_n}\}$,  and the result holds.   So suppose that $v>\textup{identity}$.   Choose a simple root $\alpha$ such that $w:=s_\alpha v<v$.    By the induction hypothesis,  $\cryiso$ induces a crystal isomorphism from $\pathslw$ to $\pathstdw$.   But by Corollary~\ref{c:l:cryisoprelim}  we have
\( \pathslv=\{f_\alpha^k\pi\st\pi\in\pathslw,\ k\geq0,\ f_\alpha^k\pi\textup{ does not vanish}\}\)
	and by Proposition~\ref{p:cryisoprelim} below
	\(\pathstdv=\{f_\alpha^k\pi\st\pi\in\pathstdw,\ k\geq0,\ f_\alpha^k\pi\textup{ does not vanish}\} \),  so it is clear that $\pathslv$ and $\pathstdv$ are crystal isomorphic via~$\cryiso$.  \hfill$\Box$ \\

\noindent	It remains only to prove the following proposition,  which has been quoted in the proof just above:
}
\blemma\mylabel{l:minstdlift}  With notation as in the statement and proof of Lemma~\ref{l:standard},  suppose that $\taut_1\geq\ldots\geq\taut_m$ be the minimal standard lift of~$\theta$. Then
\begin{enumerate}
	\item\label{i:fafact}  Suppose that $\taut_p>s_\alpha\taut_p$.   Then $p$ is either resisting or flat.    If $p$ is flat,  then there exists $r$, $p<r\leq m$,  with $r$ resisting and every $q$ such that $p<q<r$ is flat.
	\item\label{i:tautpgeqtaut} $\taut'_j\geq\taut_j$ for all $j$, $1\leq j\leq m$.
	\item\label{i:changing} Suppose that $j$ is changing and $j<m$.   Then $\taut'_{j+1}\meet s_\alpha\taut'_{j+1}=\taut_{j+1}$.
	\item\label{i:tautminstd}
$\taut'_1\geq\ldots\geq\taut'_m$ is the minimal standard lift of $f_\alpha\theta$
\end{enumerate}
\elemma
\bmyproof
\eqref{i:fafact}\ \  If $p$ is changing or changeable (but not changing),  then $\tau_p<s_\alpha\tau_p$,  so it would mean that $\taut_p<s_\alpha\taut_p$ (Corollary~\ref{c2:brcoset}).   This proves that $p$ can only be either flat or resisting.   Suppose now that $p$ is flat.  Let $r$ be the least integer, $p<r\leq m$, (if it exists) such that $r$ is not flat. If such an $r$ doesn't exist, put $r=m+1$.  For all $q$, $p\leq q<r$, put $\sigma_q:=\taut_q\meet s_\alpha\taut_q$.  Then $\sigma_p=s_\alpha\taut_p$.  We have $\sigma_p\geq\ldots\geq\sigma_{r-1}$ (by the basic fact~\eqref{e:parallel} in~\S\ref{ss:rsltsextrml}).   If $r<m$ and $r$ is not resisting,  then $\taut_r=\taut_r\meet s_\alpha\taut_r$,  so that $\sigma_{r-1}\geq\taut_r$.   Thus $\sigma_p\geq\ldots\sigma_{r-1}\geq\taut_{r}\geq\ldots\geq\taut_m$ would be a standard lift of $(\tau_p,\ldots,\tau_m)$,  which we could complete to a standard lift of $\theta$.   But then $\sigma_p=s_\alpha\taut_p<\taut_p$,  which contradicts the hypothesis that $\taut_1\geq\ldots\geq\taut_m$ is the minimal standard lift.

\eqref{i:tautpgeqtaut}\ \    Proceed by downward induction on~$j$.   Since $\taut'_j=\taut_j$ in case $j$ is obstructed,  and $\taut'_j=s_\alpha\taut_j>\taut_j$ in case $j$ is changing,  we may assume that $j$ is flat and unobstructed,  so $j<m$ and $\tau'_j=\tau_j$. We have,  by the induction hypothesis, $\taut'_{j+1}\geq\taut_{j+1}$,  and so by~\eqref{i:r:p:deodhar:ww'} of Remark~\ref{r:p:deodhar}:  \[\taut'_j\geq\brmin{J_{\tau'_j}(\taut'_{j+1})}=\brmin{J_{\tau_j}(\taut'_{j+1})}\geq\brmin{J_{\tau_j}(\taut_{j+1})}=\taut_j.\]

\eqref{i:changing}\ \ Since by definition $\taut'_{j+1}$ is either $\taut_{j+1}$ or $s_\alpha\taut_{j+1}$,   and $\taut'_{j+1}\geq\taut_{j+1}$ by item~\eqref{i:tautpgeqtaut},   it is enough to show that $\taut'_{j+1}<s_\alpha\taut'_{j+1}$.  If not,  then, by item~\eqref{i:fafact},  there exists $r$ such that $j<r\leq m$ with $r$ resisting and every $q$ such that $j<q<r$ is flat.   But this cannot happen since $j$ is changing,  by the definition of the operator $f_\alpha$.

\eqref{i:tautminstd}\ \ Let $\taut''_1\geq\ldots\geq\taut''_m$ be another standard lift of~$f_\alpha\theta$.   It suffices to show that $\taut''_j\geq\taut'_j$ for every $j$, $1\leq j\leq m$. 
Proceed by downward induction on~$j$.
It is convenient to put~$\taut_{m+1}=\taut'_{m+1}=\taut''_{m+1}=\textup{identity}$.
By the induction hypothesis, $\taut''_{j+1}\geq\taut'_{j+1}$.

In case $j$ is obstructed,
$\tau'_j=\tau_j$ by~\eqref{e:jobstructed}, and we have  
\[\taut''_j\geq
\brmin{J_{\tau'_j}(\taut''_{j+1})}=
\brmin{J_{\tau_j}(\taut''_{j+1})}\geq
\brmin{J_{\tau_j}(\taut_{j+1})} =\taut_j=\taut'_j\]
Suppose now that $j$ is changing.
Then  $\tau_j<s_\alpha\tau_j=\tau'_j$ by \eqref{e:jdefn}~and~\eqref{e:jchanging}.  
By \eqref{i:r:p:deodhar:ww'} of Remark~\ref{r:p:deodhar}, Corollary~\ref{c2:brcoset} and Lemma~\ref{l:deodhar}~\eqref{i:ldeod2}, and item~\eqref{i:changing} above:
\begin{align*}
	\taut''_j\geq \brmin{J_{s_\alpha\tau_j}(\taut''_{j+1})}
\geq \brmin{J_{s_\alpha\tau_j}(\taut'_{j+1})}
	= &s_\alpha\brmin{J_{\tau_j}(\taut'_{j+1}\meet s_\alpha\tau'_{j+1})}\\
	= &s_\alpha\brmin{J_{\tau_j}(\taut_{j+1})}=
s_\alpha\taut_j
=\taut'_j
\end{align*}

The only remaining case is when $j$ is flat and unobstructed.    We then have $j<m$ and $\tau'_j=\tau_j$.   By the induction hypothesis and item~\eqref{i:tautpgeqtaut} above,  we have $\taut''_{j+1}\geq\taut'_{j+1}\geq\taut_{j+1}$,  so by \eqref{i:r:p:deodhar:ww'} of Remark~\ref{r:p:deodhar}:
\[\taut''_j\geq
\brmin{J_{\tau'_j}(\taut''_{j+1})}
=\brmin{J_{\tau_j}(\taut''_{j+1})}
\geq \brmin{J_{\tau_j}(\taut'_{j+1})}
\geq \brmin{J_{\tau_j}(\taut_{j+1})}
= \taut_{j} \]
This means we would be done in case $\taut_{j}\geq\taut'_{j}$ (which by item~\eqref{i:tautpgeqtaut} is equivalent to $\taut_j=\taut'_j$).  But,  $\taut'_j$ is by definition the smaller of $\taut_j$ and $s_\alpha\taut_j$ that is larger than $\taut'_{j+1}$.   So it only remains to consider the case when $\taut_j<s_\alpha\taut_j$ and $\taut_j\not\geq\taut'_{j+1}$.  In this situation, $\taut_{j+1}<s_\alpha\taut_{j+1}=\taut'_{j+1}$  (for $\taut'_{j+1}$ is by definition either $s_\alpha\taut_{j+1}$ or $\taut_{j+1}$, and $\taut_j\geq\taut_{j+1}$).   This implies by item~\eqref{i:fafact} that $j+1$ is obstructed and therefore $j$ is also obstructed,  a contradiction.
%
%
\emyproof
\bmyproofof{Lemma~\ref{l:cryisoprelim2}} That $\paths_{\lambda,\textup{identity}}(\funf)=\pi_\lambda$ follows from~\eqref{e:minstdliftid}. 

Now Put $\theta = \cryiso{\pi}$ and $\weyl(\theta)=w$.    Let $\theta=(\tau_1,\ldots,\tau_m)$ and let $(\taut_1,\ldots,\taut_m)$ be the minimal standard lift of~$\theta$ (so that $w=\taut_1$).    

Proof of condition~\eqref{i:l:cryisoprelim:1} of Lemma~\ref{l:cryisoprelim}:   Let $\alpha$ be a simple root such that $s_\alpha w<w$.   To show that $e_\alpha\pi$ does not vanish,  it is enough to show that $e_\alpha\theta$ does not vanish,  and for this it is enough to show that there exists $r$, $1\leq r\leq m$, such that $s_\alpha\tau_r<\tau_r$, and $s_\alpha\tau_j=\tau_j$ for all $j$, $1\leq j<r$.      By way of contradiction,   suppose that $s_\alpha\tau_r>\tau_r$ for the least $r$ such that $s_\alpha\tau_r\neq\tau_r$ (the case when $s_\alpha\tau_j=\tau_j$ for all $1\leq j\leq m$ is included in the consideration:  we put $r=m+1$ in this case).     For $j$, $1\leq j<r$,  set $\sigma'_j:=\taut_j\meet s_\alpha\taut_j$.    Observe that $(\sigma'_1,\ldots,\sigma'_{r-1},\taut_r,\ldots,\taut_m)$ is also a standard lift of~$\theta$.    But then $\sigma'_1=s_\alpha w<w=\taut_1$,   which contradicts the choice of $(\taut_1,\ldots,\taut_m)$ as the minimal standard lift of~$\theta$.

Proof of condition~\eqref{i:l:cryisoprelim:2} of Lemma~\ref{l:cryisoprelim}: Suppose that $f_\alpha\pi$ does not vanish.  Then $f_\alpha\theta$ does not vanish either.   By Lemma~\ref{l:standard},  $f_\alpha\theta$ is standard.  Moreover, by Lemma~\ref{l:minstdlift} $(\taut'_1,\ldots,\taut'_m)$ is the minimal standard lift of~$\theta$.  Since $\taut'_1$ is either $\taut_1$ or $s_\alpha\taut_1$ by its definition,   it follows that $\funf(f_\alpha\pi)$ is either $w$ or $s_\alpha w$.   Suppose that $\funf(f_\alpha\pi)\neq\taut_1=w$.  Then, since $\taut'_1\geq\tau_1$ by item~\eqref{i:tautpgeqtaut} of Lemma~\ref{l:minstdlift},   it follows that $\funf(f_\alpha\pi)=s_\alpha w>w$.  Moreover,  this happens only if $1$ is unobstructed,  which means that minimum is~$0$ of the function $t\mapsto\langle\pi(t),\alpha^\vee\rangle$ on the interval $[0,1]$, and so $e_\alpha\pi$ vanishes.
\emyproofof

\end{appendix}

\bibliographystyle{bibsty-final-no-issn-isbn}
\ifthenelse{\equal{\finalized}{no}}{
\bibliography{abbrev,references}
}{
}

\end{document}